\documentclass[preprint,12pt]{elsarticle}

\usepackage{amssymb}
\usepackage{amsmath}
\usepackage{mathtools}
\usepackage{xcolor}
\usepackage{graphicx}
\usepackage{dsfont}
\usepackage{subcaption}
\usepackage[section]{placeins}
\usepackage{hyperref}
\newcommand{\R}{{\mathbb{R}}}
\renewcommand{\d}{{\mathrm d}}
\DeclareMathOperator{\sech}{sech}
\newcommand{\joinR}{\hspace{-.15em}}
\newcommand{\RomanI}{\scalebox{0.6}{\textit{I}}}
\newcommand{\RomanII}{\mbox{\RomanI\joinR\RomanI}}
\newcommand{\RomanIII}{\mbox{\RomanI\joinR\RomanII}}

\usepackage{amsthm}
\newtheorem{Theorem}{Theorem}[section]

\newtheorem{Lemma}[Theorem]{Lemma}

\newtheorem{Assumption}[Theorem]{Assumption}

\numberwithin{equation}{section}

\journal{Physica D: Nonlinear Phenomena}

\begin{document}

\begin{frontmatter}



\title{Long-Timescale Soliton Dynamics in the Korteweg-de Vries Equation with Multiplicative Translation-Invariant Noise}


\author[inst1]{R.W.S. Westdorp}
\ead{r.w.s.westdorp@math.leidenuniv.nl}

\affiliation[inst1]{organization={Mathematical Institute},
            addressline={Niels Bohrweg 1}, 
            city={Leiden},
            postcode={2333 CA}, 
            country={the Netherlands}}

\author[inst1]{H.J. Hupkes}
\ead{hhupkes@math.leidenuniv.nl}


\begin{abstract}
This paper studies the behavior of solitons in the Korteweg-de Vries equation under the influence of multiplicative noise. We introduce stochastic processes that track the amplitude and position of solitons based on a rescaled frame formulation and stability properties of the soliton family. We furthermore construct tractable approximations to the stochastic soliton amplitude and position which reveal their leading-order
drift. We find that the statistical properties predicted by our method agree well with numerical evidence.
\end{abstract}




\end{frontmatter}


\section{Introduction}
\label{sec:intro}
In recent years, stochastic traveling waves and more general stochastic pattern dynamics have become major areas of interest in the field of SPDEs. This paper employs modern stochastic phase-tracking techniques to study traveling waves in stochastic Korteweg-de Vries (KdV) equations with multiplicative noise, such as
\begin{align}
\label{eqn:SKDV}
    \d u=-(\partial_x^3u+2u\partial_x u)\ \d t+\sigma u\cdot \d W^Q_t.
\end{align}
Here $u$ is a real-valued process on $(t,x)\in\R^+ \times \R$, and the scalar parameter $\sigma>0$ encodes the noise strength. The noise $W^Q_t$ is a cylindrical $Q$-Wiener process on a separable Hilbert space $\mathcal{H}$, and $\cdot$ denotes a suitable product between elements of $\mathcal{H}$ and $L^2(\R)$. Both multiplicative space-time noise and multiplicative scalar noise can be treated in this setting. 

In the deterministic case ($\sigma=0$), it is well-known that \eqref{eqn:SKDV} admits soliton solutions $u(t,x)=\phi_c(x-ct)$ of the form 
\begin{align}
\label{eqn:soliton}
    \phi_c(x)=\tfrac{3c}{2} \sech^2(\sqrt{c}x/2), \quad c>0.
\end{align}
 We describe the evolution $u(t,x)$ of such a soliton
 under the influence of the multiplicative stochastic forcing ($\sigma>0$) 
 by using the modulation Ansatz 
\begin{align}
\label{eqn:simpledecomposition}
u(t,x)= \phi_{c(t)}(x-\xi(t))+r(t,x).
\end{align}
Here, $c(t)$ and $\xi(t)$ are stochastic processes that track the amplitude and position of the modulated soliton, respectively, while the perturbation $r$ remains small in a suitable sense. Numerical evidence based on these phase definitions strongly suggests that such solutions remain close to the soliton family for exponentially long times. We furthermore construct tractable approximations to the modulation parameters which reveal their leading-order drift.

\paragraph{Solitons in the Korteweg-de Vries equation}
The family of solitons \eqref{eqn:soliton} has been central to the analysis of the deterministic KdV equation ($\sigma=0$). At the time of its introduction by Boussinesq \cite{boussinesq} and rediscovery by Korteweg and de Vries \cite{korteweg} in the late nineteenth century, the KdV equation was primarily used as a model for shallow water waves along a canal. The equation has since appeared as an amplitude equation to describe a wide variety of physical wave phenomena, such as internal waves in stratified oceans \cite{internalwaves} and acoustic waves in plasmas \cite{plasma}. Let us specifically mention the 
Fermi-Pasta-Ulam-Tsingou (FPUT) chain, whose dynamics in the continuum regime can be described by the KdV equation \cite{friesecke, friesecke2,hong}. 

As a dispersive system, the dynamics generated by the KdV equation spread out localized initial conditions. Yet, due to nonlinear effects, \eqref{eqn:SKDV} with $\sigma=0$ admits the family of traveling wave solutions \eqref{eqn:soliton} that can have arbitrary amplitude,
propagating at the proportional velocity. 
The relation 
$\phi_c(x)=c \phi_1(\sqrt{c}x)$
apparent in \eqref{eqn:soliton} is a direct consequence
of the fact that the KdV equation 
enjoys 
the scaling invariance
\begin{align}\label{eqn:symmetry}
    u(t,x)\mapsto \alpha^{2}u(\alpha^3 t, \alpha x),
\end{align}
in addition to its translational symmetry.

Zabusky and Kruskal observed in numerical experiments \cite{zabusky} that, asymptotically, solutions to the KdV evolution decompose into several solitons of the form \eqref{eqn:soliton} followed by a radiation component, which undergoes a dispersive evolution. See the work of Schuur \cite{schuur} for rigorous results in this direction. Another key property of the deterministic KdV equation is that it is a completely integrable system, and consequently has an infinite number of conserved quantities \cite{tao}. For example, the KdV evolution conserves the energy $\int_{\R}u^2(t,x)\d x $.

\paragraph{Deterministic stability}
In the deterministic setting, the family of solitons in \eqref{eqn:soliton} was shown to be orbitally stable by Bona, Souganidis and Strauss \cite{bona}. Their work asserts that a slightly perturbed soliton remains, upto translations, within an $H^1$-neighborhood of the initial soliton. This result was improved upon by Pego and Weinstein, who established \textit{asymptotic} stability of the soliton family \cite{pegoweinstein}. More precisely, the authors show that if $u(t,x)$ is initially a small perturbation of the soliton $\phi_{c_*}(\cdot-\xi_*)$, then there exists a final velocity $c>0$ and a final phase-shift $\xi\in \R$ such that
\[u(t,\cdot+\xi+ct)-\phi_c \to 0 \quad \text{as} \quad t \uparrow \infty,\]
in the weighted spaces
\begin{align}\label{eqn:weightedspace}
L^2_a:=L^2(\R,e^{2ax}\d x), \quad 0<a<\sqrt{c_*}.
\end{align}
Here, the final velocity $c$ and phase-shift $\xi$ contain small corrections to their starting values $c_*$ and $\xi_*$. The exponential weight ensures that disturbances in the wake of the soliton decay at an exponential rate. The result holds under the assumption that the initial perturbation and its derivative lie in $L^2 \cap L^2_a$. This assumption was later relaxed by Merle and Vega \cite{merle} to accommodate general $L^2$-perturbations, with convergence in $L^2_{loc}$. 

The proof of Pego and Weinstein relies on the spectral stability of the linearization of the KdV evolution around a soliton $\phi_c$, encoded by the operator
\begin{align}\mathcal{L}_c = -\partial_{x}^3+(c-2\phi_c) \partial_{x}-2\partial_x \phi_c = -\partial_{x}^3+c\partial_x -2\partial_x[\phi_c \cdot].\label{eqn:linearop}
\end{align}
As a linear operator on the space $L^2(\R)$, the operator $\mathcal{L}_c$ has spectrum $i\R$. Embedded in the continuous spectrum is a double eigenvalue at 0 with 
an associated two-dimensional generalized kernel spanned by $\partial_x \phi_c$ and $\partial_c \phi_c$. Considering the operator $\mathcal{L}_c$ on a weighted space 
$L^2_a$ with $0<a<\sqrt{c}$ moves the essential spectrum to the left of the imaginary axis.\footnote{The width of the resulting spectral gap is $a(c-a^2)$, which is maximal at $a=\sqrt{c/3}$. In the literature, one often encounters the restriction $0<a<\sqrt{c/3}$. This avoids the use of spaces which are more restrictive without gaining a faster decay rate.} Its (formal) adjoint $\mathcal{L}^*_c$ on the space $L^2_{-a}$ has a generalized kernel spanned by the soliton $ \phi_c$ and the primitive 
\[\zeta_c (x)=\int_{-\infty}^x \partial_c \phi_c(y)\d y \in L^2_{-a}.\]
This primitive $\zeta_c$ is not a localized function, which is evident from the fact that $\zeta_c (x)$ tends to $\int_{\R} \partial_c \phi_c \d x =3c^{-1/2}$ as $x\to\infty$. Pego and Weinstein show that $\mathcal{L}_c$ generates an exponentially stable $C_0$-semigroup $\{e^{\mathcal{L}_{c}t}\}_{t\geq 0}$ on the subspace of $L_a^2$ consisting of functions $v\in L_a^2$ which satisfy the orthogonality conditions\footnote{In expression \eqref{eqn:orthogonality} we slightly abuse notation, as $\zeta_c \notin L^2$. The product $\int_{\R}v \zeta_c \ \d x$ is, however, a well-defined real number, since $e^{-ax}\zeta_c \in L^2$ and $e^{ax}v \in L^2$. }
\begin{align}\label{eqn:orthogonality}
\langle v, \zeta_c \rangle_{L^2}= \langle v, \phi_c \rangle_{L^2}=0.
\end{align}
 
\paragraph{Stochastic KdV equations}
Several stochastic versions of the KdV equation have been introduced in the literature, which incorporate random perturbations that affect the propagation of solitons. In \cite{herman}, Herman derives a KdV equation perturbed by a single Brownian motion to model the propagation of an ion-acoustic soliton in the presence of noise. Since the KdV equation arises as an approximation for more involved models, such as for fluid dynamics or wave dynamics in the FPUT lattice, stochastic KdV equations also serve as starting point for studying the effects of random perturbations in such systems \cite{mcginniswright, arevalo}. 

In \cite{debouardwellposedness}, de Bouard and Debussche establish the well-posedness of \eqref{eqn:SKDV} in the space $H^1(\R)$, in the case that the covariance operator $Q$ is a translation-invariant and non-negative convolution operator on $L^2(\R)$ given by
\begin{align}
\label{eqn:convolution}
    Qf(x)=\int_{\R}q(x-y)f(y)\ \d y,
\end{align}
with a convolution kernel $q$ in $H^1(\R)\cap L^1(\R)$. The same authors also prove the existence of solutions to \eqref{eqn:SKDV} in two cases that approximate a space-time white noise on $L^2(\R)$ \cite{almostwhitenoise, periodicwhitenoise}, corresponding to $Q=I_{L^2}$. 

The multiplicative forcing in \eqref{eqn:SKDV} breaks the conservative nature of the KdV equation, which can readily be seen in the scalar case where $W^Q_t$ is a real-valued Brownian motion $\beta_t$. Indeed, by formally applying It\^{o}'s lemma to the SPDE \[ \d u=-(\partial_x^3u+2u\partial_x u) \ \d t+\sigma u \ \d \beta_t, \]
one finds 
\begin{align}
    \d \langle u, u\rangle_{L^2} =& \bigl[-2 \langle u,\partial_x^3u+2u\partial_x u\rangle_{L^2}+\sigma^2 \langle u, u\rangle_{L^2} \bigr]\ \d t \nonumber\\
    &+2\sigma\langle u, u\rangle_{L^2}\ \d \beta_t\nonumber\\
    =& \sigma^2 \langle u, u\rangle_{L^2}\ \d t+2\sigma\langle u, u\rangle_{L^2}\ \d \beta_t, \label{eqn:energygbm}
\end{align}
which shows that the energy $\int_{\R}u^2(t,x) \d x $ undergoes a geometric Brownian motion with positive drift. See \cite[Proposition 3.1]{debouardwellposedness} for a similar result in the case of space-time noise. As the stochastic forcing slightly perturbs the soliton continually, its aggregated effect produces a stochastic phase-shift and a substantially varying soliton parameter $c(t)$.

The conservative nature of the equation is, however, not entirely lost. One easily verifies that the \textit{average} $L^2$-norm of solutions to the stochastic KdV equation is conserved: \[\mathbb{E}\bigl[\|u(t,\cdot)\|_{L^2}\bigr]=\|u_0\|_{L^2}.\]
Due to its diffusion, the same does not hold for nonzero powers $p\in \R$ of the $L^2$-norm. 
In this context it is worthwhile to point out 
that the soliton parameter $c$ is related to the $L^2$-energy as $\|\phi_c\|_{L^2}^2=6c^{3/2}$. This hints at a relation between the stochastic soliton parameter $c(t)$ and the process 
\begin{align}\|u(t,\cdot)\|_{L^2}^{4/3}=\|u_0\|_{L^2}^{4/3}e^{-\frac{2}{3}\sigma^2 t+\frac{4}{3}\sigma  \beta_t}, \label{eqn:heuristic}
\end{align}
a basic prediction that the results in this paper will reproduce and refine.

\paragraph{Stochastic traveling waves}
The effect of noise on traveling waves has been previously been analyzed in various settings, primarily for equations of reaction-diffusion type. The subject has seen considerable activity in the physics literature, see for instance the works by Garc\'{i}a-Ojalvo and Schimansky-Geier\cite{garcia, geier} which introduce stochastic corrections to traveling waves in bistable RDEs. One well-known example of a noise-induced velocity correction is the Brunet-Derrida conjecture \cite{brunet1997shift,brunet2001effect}, which describes a speed correction for traveling fronts in randomly perturbed Fischer-KPP equations and was proved by Mueller \cite{mueller}. We also refer to the works \cite{lordpaper, lord}, which analyze numerical methods for the simulation of stochastic traveling waves and provide numerous intriguing observations. 

Various contributions have appeared in the mathematics literature to give (further) rigorous meaning to such results. Kr\"{u}ger and Stannat introduced a multiscale expansion of a stochastic phase for traveling waves in bistable RDEs \cite{krugerstannat, stannat2014stability}. This method was later applied to the FitzHugh-Nagumo equation in \cite{gnann} and extended upon in \cite{maclaurinnew}. Hamster and the second author have developed a phase-tracking method in the setting of reaction-diffusion equations that tracks stochastic traveling waves over exponentially long timescales \cite{hamsterscalar, hamster, hamsterstability}. We refer to the review by Kuehn \cite{kuehnreview} for a more detailed overview of results on stochastic traveling waves in reaction-diffusion equations. 

\paragraph{Stochastic KdV waves}
The stochastic dynamics of solitons in a randomly perturbed KdV equation was first considered by Wadati \cite{wadati}, who  derived statistical properties of solitons in the KdV equation with additive scalar noise using a Galilean coordinate transformation. See also the works \cite{wadati2, garnier}, which expand on this method. 

De Bouard and Debussche analyzed the stochastic soliton dynamics produced by the KdV equation \eqref{eqn:SKDV}  with multiplicative space-time noise in \cite{debouardsoliton}. Similar to the approach taken in this work, the authors consider a decomposition of the form \eqref{eqn:simpledecomposition} and formulate modulation equations for the soliton parameters. The authors, moreover, estimate the exit-time of the solution from a neighborhood of the modulated soliton. See also the works \cite{debouardadditive,NLSmodulation}, where the same authors study solitons in the KdV equation with additive noise and a stochastic Gross-Pitaevskii equation.

While the method in \cite{debouardsoliton} yields rigorous stability results, it only allows for a small variation of the soliton parameter $c(t)$. The resulting modulation equations, therefore, do not incorporate the full dynamics of the soliton amplitude and are valid only on timescales of order $O(1/\sigma^2)$, where the parameter $c(t)$ remains close to its starting value $c_*$. Our method allows for large variations in the parameter $c(t)$ by not only translating the solution, but also rescaling the solution in accordance with the natural scaling $u(t,\cdot) \mapsto \alpha^2 u(t,\alpha \cdot)$ of the soliton family \eqref{eqn:soliton}.

In a more recent work by Cartwright and Gottwald \cite{cartwright}, the authors apply a collective coordinate framework developed in \cite{cartwright2019collective} to \eqref{eqn:SKDV} in the scalar case, where $ W_t^Q$ in \eqref{eqn:SKDV} is a scalar Brownian motion $ \beta_t$. One version of their approach modulates the amplitude and width of the soliton independently, whereas another version requires that the modulated soliton remains in the family \eqref{eqn:soliton}. The former method is found to be advantageous in case the perturbations are large. The authors develop an $O(\sigma)$ modulation system based on this method.
Our work here enables the inclusion of higher order effects, which are essential to capture the dominant behaviour of fundamental soliton features such as the width.

\paragraph{Soliton tracking}
In this paper, we introduce a soliton-tracking method that adapts a phase-tracking method developed in \cite{hamster} by Hamster and the second author for traveling waves in reaction-diffusion equations. The method relies on a transformation of \eqref{eqn:SKDV} that translates and rescales a solution $u(t,x)$ to closely match a fixed soliton $\phi_{c_*}$ at the origin. More precisely, 
we introduce the process 
\begin{align}
\label{eqn:decomposition}
    v(t,x)=\alpha^2(t)u\bigl(t,\alpha(t) x+\xi(t)\bigr)-\phi_{c_*}(x),
\end{align}
which is the remaining difference between the translated and rescaled solution, and the fixed soliton $\phi_{c_*}$. From the scaling process $\alpha(t)$, we can recover the effective soliton parameter $c(t)$ of the solution $u(t,x)$ as $c(t)=c_* \alpha^{-2}(t)$. 

The translation process $\xi(t)$ and scaling process $\alpha(t)$ are a-priori not specified. Their drift and martingale components bring about four degrees of freedom. We formulate an SPDE that governs the dynamics of the remainder $v$ produced by translation and rescaling as in \eqref{eqn:decomposition} by noise-driven processes $\xi(t)$ and $\alpha(t)$. We use the four degrees of freedom brought about by $\xi(t)$ and $\alpha(t)$ to ensure that  $v$ satisfies the orthogonality conditions of \eqref{eqn:orthogonality} at all times. Intuitively, this should mean that the remainder $v$ continuously experiences exponential damping, and will remain small. This argument has been made rigorous by Hamster and the second author for traveling waves in reaction-diffusion equations in \cite{hamsterstability}, where the authors establish exit-times on the remainder which are exponentially long with respect to the parameter $1/\sigma$.  

We restrict ourselves in this work to formal arguments to support the construction of the soliton-tracking method. Our method produces a coupled SPDE system that describes the evolution of $v(t), c(t)$ and $\xi(t)$, which we characterize in considerable detail. On the one hand, this is intended to facilitate the re-use of our ideas in other contexts. On the other hand, we view this paper as the basis for ongoing efforts to rigorously establish long-time stability of the KdV soliton under stochastic perturbations. We complement the presentation of our general method with three worked-out examples, which allow us to showcase the strong predictive power of our method, as validated by extensive numerical
simulations.

Although we anticipate that a  meta-stability result based on the methods presented in this work is possible, one needs to control transformations of the perturbation due to the rescaling. It seems that this requires spatial information on the linearized frozen-frame evolution beyond the semigroup bounds obtained in \cite{pegoweinstein}.

\paragraph{Stochastic soliton dynamics}

In order to gain insight into the stochastic soliton dynamics described by our method, we design an approximation procedure for the processes $\alpha(t)$ and $\xi(t)$. 
In contrast to the stochastic wave position studied in \cite{hamsterscalar} and \cite{hamster}, the evolution of the perturbation $v(t)$ does not decouple from the rescaling process $\alpha(t)$. This significantly complicates our analysis, since $\alpha(t)$ develops significant fluctuations on short timescales. In order to account for this, we first expand the perturbation $v(t)$ in terms of the small parameter $\sigma$, using $\alpha$ as an external input. Expanding the modulation equation for $\alpha(t)$ in terms of $v$ subsequently produces SDE approximations for $\alpha(t)$ that have random coefficients. Solving these SDEs then provides the desired final approximations, which can in principle be computed at any desired order in $\sigma$. This non-standard nested approximation procedure is discussed in detail in \S\ref{sec:solitondynamics}.

We find that the soliton velocity primarily follows the amplitude process $c(t)$, with additional correction due to the noise. The amplitude process $c(t)$ experiences an almost linear positive drift which develops on the time-scale $O(\sigma^2 t)$. Numerical simulations of \eqref{eqn:SKDV} show that this drift is captured quite well by an approximation of $c(t)$ that takes into account coupling terms that are quadratic in $\sigma$. 

The amplitude growth rate and velocity correction are determined by statistical properties of the perturbation with respect to the modulated soliton. This further motivates  the use of the frozen-frame formulation \eqref{eqn:decomposition}, which keeps the stochastic traveling wave in a fixed reference frame and facilitates analysis of the perturbation shape. We emphasize that the general approach is not restricted to the KdV equation but can also be used in other settings where (approximate) self-similarity properties hold.

\paragraph{Outlook}
This work extends the phase-tracking results of \cite{hamsterscalar} and \cite{hamster} to a system where stochastic perturbations not only induce a translation, but also a rescaling of the traveling wave. Numerical simulations indicate that our method tracks the KdV soliton over exponentially long times. However, this work does not provide a rigorous stability result to support this claim. We hope that the decompositions and approximations obtained here will provide a pathway towards such a result.

\paragraph{Outline} 
This paper is organized as follows. In \S\ref{sec:tracking}, we provide a detailed derivation of our phase-tracking approach and the resulting modulation system for the soliton parameters. We also introduce the three example systems that we use throughout the work to illustrate our techniques and results.
In \S\ref{sec:solitondynamics}, we formulate an expansion of the SPDE system that governs the dynamics of $\alpha(t)$ and $v(t)$, and from there derive the leading-order behavior of the mean and variance of the effective soliton parameter $c(t)$ and position $\xi(t)$. We verify these findings via numerical simulations throughout \S\ref{sec:solitondynamics}.

\paragraph{Acknowledgements}
The first author acknowledges support from the Netherlands Organization for Scientific Research (NWO) (grant 613.009.137).

\section{Stochastic soliton tracking}
\label{sec:tracking}

In this section, we introduce a system of modulation equations to describe the evolution of the solution $u(t,x)$ to the stochastic KdV equation in a stochastic co-moving frame. In particular, we consider the SPDE
\begin{align}
\label{eqn:skdvgeneral}
    \d u&= -(\partial_x^3 u +2u \partial_x u)\d t+\sigma M(u) [\d W_t^Q],
\end{align}
with initial condition $u(0,x)=\phi_{c_*}(x)$, where $\phi_{c_*}$ is the soliton defined in \eqref{eqn:soliton}. The noise term $M(u)$ is of general multiplicative form and is driven by a translation-invariant noise process $W_t^Q$, which we both describe in more detail in \S\ref{subsec:setup}.

We introduce a position process $\xi$, which (roughly) keeps the stochastically evolving soliton centered at the origin, and a rescaling process $\alpha$, which neutralizes its amplitude fluctuations. More precisely, we introduce the remainder 
\begin{align}\label{eqn:remainder2}
    v(t,x)=\alpha^2(t)u\bigl(t,\alpha(t) x+\xi(t)\bigr)-\phi_{c_*}(x),
\end{align}
which describes the deviation from the soliton $\phi_{c_*}$ in a frame where the solution $u(t,x)$ has been translated and rescaled according to the natural scaling $\phi_c \mapsto \alpha^2 \phi_c(\alpha \cdot)$ of the soliton family. The aim is to choose the processes $\alpha$ and $\xi$ in a fashion that keeps the perturbation $v$ in the space characterized by \eqref{eqn:orthogonality}, where the linearized evolution is stable. 

First, in \S\ref{subsec:setup}, we describe the forms of multiplicative noise that can be treated by our method. We introduce three concrete settings that will appear throughout the paper: scalar noise, translation-invariant colored noise and space-time white noise. In \S\ref{subsec:frozen}, we derive an SPDE that describes the evolution of the perturbation $v$ in the co-moving frame, {which we use in \S\ref{subsec:modulation} to prescribe the dynamics of our processes $\xi$ and $\alpha$.
This results in a modulation system of the form
\begin{align*}
    \d v =& \alpha^{-3}\bigl[\mathcal{L}_{c_*} v + O(v^2)\bigr]\ \d t+\sigma^2 O(1)\ \d t+ \sigma O(1)\ \hat{T}_{\alpha}\d W_t^Q,\\
     \d \alpha =&\ \bigl[\alpha^{-2}O(v^2)+\sigma^2 \overline{\gamma}_d(v,\alpha)\bigr]\ \d t +\sigma \alpha O(1)\ \hat{T}_{\alpha} \d W_t^Q,\\
    \d \xi =&\ \bigl[\alpha^{-2}(c_*+O(v^2))+\sigma^2 \overline{\mu}_d(v,\alpha)\bigr] \ \d t +\sigma  \alpha O(1)\ \hat{T}_{\alpha} \d W_t^Q.
\end{align*}
Here, $\mathcal{L}_{c_*}$ is the linear operator introduced in \eqref{eqn:linearop} and the operator $\hat{T}_{\alpha}$ corrects for the spatial rescaling by $\alpha$. The drift terms $\overline{\gamma}_d(v,\alpha)$ and $\overline{\mu}_d(v,\alpha)$ are both $O(1)$ in $v$, but their dependence on $\alpha$ varies per noise type. We explore these differences by returning to the three example settings.

Finally, we illustrate the effectiveness of our decomposition in \S\ref{subsec:numerics} 
via numerical simulations.
Since our approach keeps $v$ in the stable subspace of $\mathcal{L}_{c_*}$, we expect that the perturbation $v$ only grows logarithmically in time. Our numerical simulations strongly suggest that this is indeed the case.


\subsection{Stochastic set-up}\label{subsec:setup}
Let us outline the setting of \eqref{eqn:skdvgeneral} in more detail. We follow the approach of \cite{daprato} and  \cite{liurockner}, and consider noise from a separable Hilbert space $\mathcal{H}$ with the inner product $\langle \cdot,\cdot \rangle_{\mathcal{H}}$ and an orthonormal basis $\{e_k\}_{k=0}^\infty$. We then pick a covariance operator $Q$ that satisfies the following properties:
\begin{Assumption}\label{ass:Q}
The operator $Q:\mathcal{H}\to\mathcal{H}$ is linear and bounded, and for each $h,g\in \mathcal{H}$ we have
    \begin{itemize}
        \item $\langle Q h, h\rangle_{\mathcal{H}}\geq 0$;\quad (non-negativity)
    \item $\langle Q h, g\rangle_{\mathcal{H}}=\langle  h,Q g\rangle_{\mathcal{H}}$. \quad (symmetry)
    \end{itemize}
\end{Assumption}
With this assumption in place, we let $W^Q_t$ be a $Q$-cylindrical Wiener process on $\mathcal{H}$; cf. \cite[\S4.1.2]{daprato} and \cite[\S2.5.1]{liurockner}. This $\mathcal{H}$-valued Wiener process has the property that for each $h\in \mathcal{H}$, the process $\langle W^Q_t, h\rangle$ defines a real-valued Brownian motion, which satisfies the correlation identity 
\[\mathbb{E}\bigl[\langle W^Q_t, h\rangle \langle W^Q_s, g\rangle\bigr]= (t \wedge s)\langle Q h,g\rangle,\]
for $t,s\geq 0$ and $h,g\in \mathcal{H}$. 

In order to define a stochastic integral with respect to $W_t^Q$, we follow \cite{daprato,liurockner} and introduce the space
$\mathcal{H}_Q:=Q^{1/2}(\mathcal{H})$. Equipped with the inner product 
\[\langle v, w\rangle_{\mathcal{H}_Q}=\langle Q^{-1/2}v,Q^{-1/2}w\rangle_{\mathcal{H}},\]
$\mathcal{H}_Q$ is a separable Hilbert space for which $\{Q^{1/2}e_k\}_{k=0}^\infty$ is an orthonormal basis. Let us furthermore introduce the notation $HS(\mathcal{H}_Q,\mathcal{H})$ for the space of Hilbert-Schmidt operators between $\mathcal{H}_Q$ and $\mathcal{H}$. We recall that $HS(\mathcal{H}_Q,\mathcal{H})$ is a Hilbert space with the inner-product 
\[\langle A,B\rangle_{HS(\mathcal{H}_Q,\mathcal{H})}=\sum_{k=0}^\infty \bigl\langle A[Q^{1/2}e_k],B[Q^{1/2}e_k]\bigr\rangle_{\mathcal{H}}.\]

We refer to \cite[\S4.2.1]{daprato} and \cite[\S2.5.2]{liurockner} for the construction of the stochastic integral 
\begin{align}\label{eqn:stochint}
\int_0^t \Phi(s) \ \d W_s^Q
\end{align}
with respect to $W_t^Q$, which defines an $\mathcal{H}$-valued stochastic process for $HS(\mathcal{H}_Q,\mathcal{H})$-valued integrands $\Phi$ that satisfy the integrability condition
\[\mathbb{E}\int_0^t \|\Phi(s)\|^2_{HS(\mathcal{H}_Q,\mathcal{H})}\ \d s < \infty.\]

For convenience, we also introduce a rescaling and translation transformation $T_{\alpha,\xi}$, which acts on functions $u\in L^2(\R)$ as 
\begin{align}\label{eqn:mappingT}
    T_{\alpha,\xi}u= u(\alpha\cdot+\xi).
\end{align}
This allows us to write \eqref{eqn:remainder2} as
\begin{align}
\label{eqn:remainder}
v(t,x)=\alpha^2(t) T_{\alpha(t),\xi(t)}u(t,x)-\phi_{c_*}(x).
\end{align}
With these preliminaries in place, we impose the following conditions on the noise-term $M(u)$ and its relation to the scaling operators $T_{\alpha,\xi}$.
\begin{Assumption}\label{ass:multiplicative}
For each $u\in L^2(\R)$:
\begin{enumerate}
    \item $M(u)$ defines a Hilbert-Schmidt operator from $\mathcal{H}_Q$ to $L^2(\R)$. If furthermore $u\in H^1(\R)$, then $M(u)$ defines a bounded linear operator from $\mathcal{H}$ to $L^2(\R)$.
\item For each $\beta \in \R$ and $h\in\mathcal{H}$ we have the identity\footnotemark
\[\beta M(u)[h]=M(\beta  u)[h].\]
\item There exists a bounded linear operator $\hat{T}_{\alpha,\xi}$ on $\mathcal{H}$, such that for each $\alpha>0$ and $\xi \in \R$ we have
\[T_{\alpha,\xi}\bigl[M(u)[h]\bigr]=M(T_{\alpha,\xi}u)[\hat{T}_{\alpha,\xi}h].\]
\item There exists a linear operator $\hat{\partial}_x$ on $\mathcal{H}_Q$ such that we have\footnotemark[\value{footnote}]
 \[\partial_x\bigl[M(u)[h]\bigr]=M(u_x)[h]+M(u)[\hat{\partial}_xh],\]
for every $h\in \mathcal{H}_Q$ and $u\in H^1(\R)$.
\item The translation invariance identities
\begin{align*}
\hat{T}_{\alpha,\xi} Q\hat{T}^*_{\alpha,\xi}=\hat{T}_{\alpha} Q\hat{T}^*_{\alpha} \quad \text{and }\quad \hat{T}_{\alpha,\xi} \hat{\partial}_x Q\hat{T}^*_{\alpha,\xi}=\hat{T}_{\alpha}  \hat{\partial}_x Q\hat{T}^*_{\alpha}
\end{align*}
hold for each $\alpha>0$ and $\xi\in \R$, where $\hat{T}_{\alpha}$ is shorthand for $\hat{T}_{\alpha,0}$ and $\hat{T}^*_{\alpha,\xi}$ denotes the $\mathcal{H}$-adjoint of $\hat{T}_{\alpha,\xi}$.\end{enumerate}
\end{Assumption}

\footnotetext{Items 2 and 4 of Assumption~\ref{ass:multiplicative} restrict the setting to linear noise terms, in order to keep our computations tractable. More general noise terms can be treated by applying a function $g:\R \to \R$ point-wise and considering the noise term $M\bigl(g(u)\bigr)$. Items 2 and 4 then generalise to $\beta M\bigl(g(u)\bigr)[h]=M\bigl(\beta g(u)\bigr)[h]$ and $\partial_x\bigl[M\bigl(g(u)\bigr)[h]\bigr]=M\bigl(g^\prime(u)u_x\bigr)[h]+M\bigl(g(u)\bigr)[\hat{\partial}_x h]$.}

The operator $M$ takes on the role of multiplication between $u$ and an element of $\mathcal{H}$. We remark that for each $u\in H^1(\R)$
the operator $M(u)$ has an adjoint $M^*(u):L^2(\R)\to \mathcal{H}$, which by definition satisfies
\begin{align}\label{eqn:adjoint}
\bigl\langle M(u)[h],f\bigr\rangle_{L^2(\R)} =\bigl\langle h, M^*(u)[f]\bigr\rangle_{\mathcal{H}}\end{align}
for each $h\in \mathcal{H}$ and $f\in L^2(\R)$. 

\paragraph{Solution types}
With this assumption in place, we can assign a rigorous meaning to the SPDE \eqref{eqn:skdvgeneral}. Based on the flow $S(t)=e^{-t\partial_x^3}$ generated by the linear equation $u_t=-\partial_x^3 u$, we call an $H^1$-valued process $u$ a \textit{mild solution} to \eqref{eqn:skdvgeneral} with initial condition $u(0,\cdot)=u_0$ if the mild formula
\begin{align}
    u(t)=S(t)u_0-\int_0^t S(t-s) \partial_x\bigl(u^2(s)\bigr) \ \d s+\sigma\int_0^tS(t-s) M\bigl(u(s)\bigr) [\d W^Q_s]\label{eqn:mild}
\end{align}
holds for all $t>0$. The integral equation \eqref{eqn:mild} generalizes the variation of constants formula, or Duhamel formula, to the SPDE setting and allows for solutions of lower regularity. As is the case with PDEs, a mild solution is also classical if it takes values in the domain of the linear operator.

Alternatively, one can also consider weak solutions: an $H^1$-valued process $u$ is a \textit{weak solution} to \eqref{eqn:skdvgeneral} with initial condition $u(0,\cdot)=u_0$ if the identity 
\begin{align}
    \langle u(t), \zeta\rangle_{L^2}=&\langle u_0,\zeta \rangle_{L^2}-\int_0^t \langle u_x(s),\partial_x^2 \zeta \rangle_{L^2}-2\langle u(s)u_x(s),\zeta \rangle_{L^2}\ \d s \nonumber\\
    &+\sigma \int_0^t \langle  M(u(s))[\d W_s^Q],\zeta \rangle_{L^2} \label{eqn:weak}
\end{align}
holds for all $\zeta\in H^2(\R)$ and $t>0$. See for instance \cite{weaksol}, which asserts the existence of weak solutions for \eqref{eqn:skdvgeneral} posed on a bounded domain with scalar noise. This solution type is, however, less common in the stochastic KdV literature. We refer to \cite[Appendix G]{liurockner} for an overview of various solution concepts for SPDEs and the relationships between them.

Let us illustrate the broad applicability of this abstract setting with three examples, which we use throughout the paper to showcase our results.

\subsubsection{Example I: Scalar multiplicative noise}
\label{subsubsec:scalar}
As a first example, consider the KdV equation perturbed by a single Brownian motion $\beta_t$, which we write as
\begin{align}\label{eqn:scalar}
    \d u=-(\partial_x^3 u +2u \partial_x u) \ \d t+\sigma u \ \d \beta_t.
\end{align}
Here, the Brownian motion takes values in the Hilbert space $\mathcal{H}=\R$. Since the perturbation is uniform in space, the covariance operator acts trivially as $Q=I_{\R}$. In this case, we have $\mathcal{H}_Q=\R$ and the multiplication between the noise and the function $u$ is simply given by $M_{I}(u):\R\to L^2(\R)$, which acts on $h\in \R$ as 
\[M_{I}(u)[h]=h u.\] 
By the defining identity \eqref{eqn:adjoint}, the formal adjoint $M_{I}^*(u)$ must satisfy
\[\langle hu,f\rangle_{L^2}=hM_{I}^*(u)[f], \] 
which implies that $M_{I}^*(u):L^2(\R)\to \R$ acts on $f\in L^2(\R)$ as
\[M_{I}^*(u)[f]=\langle u,f\rangle_{L^2(\R)}.\]
We furthermore compute
\begin{align*}
    T_{\alpha,\xi}\bigl[M_{I}(u)[h]\bigr]=&h T_{\alpha,\xi}u=M_{I}(T_{\alpha,\xi}u)[h],\\
    \partial_x\bigl[M_{I}(u)[h]\bigr]=&h u_x=M_{I}(u_x)[h],
\end{align*}
which shows that $M_{I}$ fits items 3 and 4 of Assumption~\ref{ass:multiplicative} with $\hat{T}_{\alpha,\xi}=I_{\R}$ and $\hat{\partial}_x$ acting as $\hat{\partial}_x h=0$ for all $h\in \R$.

As pointed out in \cite{cartwright}, we remark that \eqref{eqn:scalar} can be transformed into a KdV equation with a random coefficient. Indeed, by factoring out the geometric Brownian motion 
\[g(t)=e^{-\frac{\sigma^2}{2}t+\sigma\beta_t}\] 
as $v=g^{-1}u$, we find  
\[\d v=-(\partial_x^3 v+2 g(t)v\partial_x v) \ \d t.\]
This reduces matters of well-posedness to the well-posedness of a KdV equation where the nonlinearity is multiplied with a varying coefficient. We are, however, unaware of results that establish the well-posedness of such an equation with a non-smooth coefficient.

\subsubsection{Example I\joinR I: Translation-invariant colored noise}
\label{subsubsec:spacetimecolored}
Our setting also allows for space-time noise from the Hilbert space $\mathcal{H}=L^2(\R)$, with the usual inner product. 
We consider translation-invariant noise, with a spatial correlation structure given by an even function $q\in H^1(\R)\cap L^1(\R)$ that has a non-negative Fourier transform $\hat{q}$. We then introduce the covariance operator $Q$ that acts on $f\in L^2(\R)$ as the convolution
\begin{align}\label{eqn:convolution2}
    Qf(x)=\int_{\R}q(x-y)f(y)\ \d y
\end{align}
with respect to the kernel $q$. The integrability of $q$ ensures that $Q$ is a bounded operator on $L^2(\R)$, and the non-negativity of the Fourier transform $\hat{q}$ provides the remaining properties of $Q$ in Assumption~\ref{ass:Q}.

Well-posedness in this setting is asserted in \cite{debouardwellposedness}, where the authors construct mild solutions as per \eqref{eqn:mild}. We have the formal covariance identity
\begin{align*}
    \mathbb{E}\bigl[\d W^Q(x,t)\d W^Q(y,s)\bigr]=\delta(t-s)q(x-y),
\end{align*}
which quantifies how $q$ determines the spatial correlation of the noise, depending only on the distance between two points. 

Using the kernel $q$, it is possible to provide an explicit formulation of the operator $Q^{1/2}$. To this end, note that $Q$ acts as a Fourier multiplier with symbol $\hat{q}$. It follows from the assumption $q\in H^1(\R)\cap L^1(\R)$ and elementary properties of the Fourier transform that $\hat{q}$ lies in $L^1(\R)$ and is bounded. As a consequence, the function $\xi \mapsto \sqrt{\hat{q}(\xi)}$ is also bounded, and defines a Fourier multiplier which is bounded on $L^2(\R)$. Denoting the inverse Fourier transform of $\sqrt{\hat{q}}$ by $q_{1/2}$ allows us to characterize $Q^{1/2}$ as
\begin{align*}
    Q^{1/2}f(x)=\int_{\R}q_{1/2}(x-y)f(y)\ \d y.
\end{align*}
For future reference, it is convenient to introduce here a rescaled family $\{Q_\alpha\}_{\alpha>0}$ of the convolution operator $Q$, which rescales correlation lengths of the kernel $q$ as 
\begin{align}
    Q_\alpha f=\alpha q(\alpha \cdot)*f. \label{eqn:convolutionfamily}
\end{align}
Similarly, we set
\[ (Q^{1/2})_\alpha f=\alpha q_{1/2}(\alpha \cdot)*f.\]

In applications, one often encounters the Gaussian kernel 
\begin{align}\label{eqn:gaussian}
    q(x)=\frac{1}{2\zeta}e^{\frac{-\pi x^2}{4 \zeta^2}},
\end{align}
where the parameter $\zeta>0$ is a measure for the correlation length. In this case, $\hat{q}$ is given by 
\[\hat{q}(\xi)=\frac{1}{\sqrt{2\pi}\zeta^2}e^{\frac{-\zeta^2\xi^2}{\pi}}.\] We note that, with the kernel \eqref{eqn:gaussian}, $Q_\alpha$ acts as
\[Q_\alpha f= \frac{\alpha}{2\zeta}e^{\frac{-\pi \alpha^2 x^2}{4 \zeta^2}}*f,\]
so that the correlation length $\zeta$ of the Gaussian kernel is effectively rescaled to $\zeta/\alpha$.

In the setting of this example, the space-time noise enters the KdV equation via the operator $M_{\RomanII}(u):L^2(\R) \to L^2(\R)$ that acts on $h\in L^2(\R)$ as the point-wise multiplication 
\[\bigl(M_{\RomanII}(u)[h]\bigr)(x)=h(x) u(x).\]
In order to verify that $M_{\RomanII}(u)$ is in the class $HS(L_Q^2,L^2)$ for each $u\in L^2(\R)$, we compute
\begin{align*}
    \|M_{\RomanII}(u)\|^2_{HS(L_Q^2,L^2)}=&\sum_{k=0}^\infty \langle u Q^{1/2}e_k,u Q^{1/2}e_k\rangle_{L^2}\\
    =&\sum_{k=0}^\infty \int_{\R}u(x)^2\bigl\langle q_{1/2}(x-\cdot),e_k\bigr\rangle_{L^2}^2 \d x\\
    =&\int_{\R}u(x)^2\bigl\langle q_{1/2}(x-\cdot),q_{1/2}(x-\cdot)\bigr\rangle_{L^2}  \d x=\|q_{1/2}\|_{L^2}^2\|u\|_{L^2}^2.
\end{align*}
We then note that 
\[\|q_{1/2}\|_{L^2}^2=\|\sqrt{\hat{q}}\|_{L^2}^2=\|\hat{q}\|_{L^1},\]
via Parseval's theorem, and we conclude that the operator $M_{\RomanII}(u)$ is indeed Hilbert-Schmidt.

Applying the adjoint identity \eqref{eqn:adjoint} to the multiplication operator $M_{\RomanII}$ yields
\[\langle hu,f\rangle_{L^2}=\bigl\langle h,M_{\RomanII}^*(u)[f]\bigr\rangle_{L^2}, \] 
which reveals that $M_{\RomanII}(u)$ is self-adjoint:
\[M_{\RomanII}^*(u)[f]=uf=M_{\RomanII}(u)[f].\]
We can furthermore compute
\begin{align*}
    T_{\alpha,\xi}\bigl[M_{\RomanII}(u)[h]\bigr]=&\ T_{\alpha,\xi}uT_{\alpha,\xi}h=M_{\RomanII}(T_{\alpha,\xi}u)[T_{\alpha,\xi} h],\\
    \partial_x\bigl[M_{\RomanII}(u)[h]\bigr]=&\ u_x h+u h_x=M_{\RomanII}(u_x)[h]+M_{\RomanII}(u)[h_x],
\end{align*}
which shows that $M_{\RomanII}$ satisfies items 3 and 4 of Assumption~\ref{ass:multiplicative} with $\hat{T}_{\alpha,\xi}=T_{\alpha,\xi}$ and $\hat{\partial}_x=\partial_x$. Lemma~\ref{lem:scaledQ} then shows that we have the identities
\begin{align}\hat{T}_{\alpha,\xi}Q\hat{T}^*_{\alpha,\xi}=\alpha^{-1}Q_\alpha \quad \text{and} \quad
\hat{T}_{\alpha,\xi}\hat{\partial}_x Q\hat{T}^*_{\alpha,\xi}=\alpha^{-2}Q_\alpha \partial_x,\label{eqn:transidenties}
\end{align}
where we recall that the operator $Q_\alpha$ is defined in \eqref{eqn:convolutionfamily}.
Consequently, item 5 of Assumption~\ref{ass:multiplicative} is met. 

\subsubsection{Example I\joinR I\joinR I: Space-time white noise}
\label{subsubsec:spacetimewhite}
The setting described above formalizes translation-invariant space-time noise with arbitrary correlation length. It is inviting to consider the limiting case where the correlation length $\zeta$ in \eqref{eqn:gaussian} tends to zero. This leads to a space-time white noise $W_t$, which is completely uncorrelated in space and time as specified by the formal identity 
\begin{align*}
    \mathbb{E}[\d W(x,t)\d W(y,s)]=\delta(t-s)\delta(x-y).
\end{align*}
Upon doing so, the regularising effect of the convolution with respect to the kernel $q$ is lost. As a consequence, it is unclear how to interpret the noise term $M_{\RomanII}(u) \d W_t$, since $M_{\RomanII}(u)$ is not in the class $HS(L^2,L^2)$ and violates item 1 of Assumption~\ref{ass:multiplicative}. Indeed, to our knowledge, no well-posedness theory is currently available for the KdV equation with multiplicative space-time white noise. 

Despite these limitations, we can proceed formally by considering the covariance operator $Q=I_{L^2}$, since convolution with respect to the Dirac distribution acts as the identity operator $I_{L^2}$. We stress that this procedure only amounts to a formal computation, but as we shall see in the sequel the results are very insightful. In this case, 
the translational invariance identities simplify to
\[\hat{T}_{\alpha}Q\hat{T}^*_{\alpha}=\alpha^{-1} \quad \text{and} \quad \hat{T}_{\alpha}\hat{\partial}_x Q\hat{T}^*_{\alpha}=\alpha^{-2} \partial_x.\]

\subsection{Frozen-frame transformation}
\label{subsec:frozen}

In this section, we formulate an SPDE that governs the evolution of the remainder $v$ in the co-moving frame introduced in \eqref{eqn:remainder}.
We postulate that the rescaling and translation processes $\alpha$ and $\xi$ satisfy the system
\begin{align}
    \d \alpha &= \gamma^\sigma_d (v,\alpha) \ \d t +\gamma^\sigma_s(v,\alpha,\xi) \ \d W_t^Q, \label{eqn:alphaform}\\
    \d \xi &= \mu^\sigma_d(v,\alpha) \ \d t +\mu^\sigma_s(v,\alpha,\xi) \ \d W_t^Q,\label{eqn:xiform}
\end{align}
where $\gamma^\sigma_d, \mu^\sigma_d$ are real-valued and $\gamma^\sigma_s, \mu^\sigma_s$ are linear operators from $\mathcal{H}$ to $\R$. The drift components $\gamma^\sigma_d, \mu^\sigma_d$ and martingale components $\gamma^\sigma_s, \mu^\sigma_s$ will be explicitly provided in \S\ref{subsec:modulation}. For now we note that $\gamma^\sigma_s$ and $\mu^\sigma_s$ can be written as
\begin{align}\label{eqn:gammasrepresentation}
    \gamma^\sigma_s (v,\alpha,\xi)[h] = -\sigma \alpha \langle \hat{T}_{\alpha,\xi}h , \overline{\gamma}_s(v)\rangle_{\mathcal{H}},\\ 
    \mu^\sigma_s (v,\alpha,\xi)[h] = -\sigma \alpha \langle \hat{T}_{\alpha,\xi} h, \overline{\mu}_s(v) \rangle_{\mathcal{H}},\label{eqn:musrepresentation}
\end{align}
where $\overline{\gamma}_s$ and $\overline{\mu}_s$ are $\mathcal{H}$-valued, depending only on $v$. In addition, the drift components $\gamma^\sigma_d$ and $\mu^\sigma_d$ are of the form
\begin{align}\label{eqn:gammadrepresentation}
    \gamma^\sigma_d(v,\alpha)=&-\alpha^{-2}\overline{\gamma}_d^0(v)+\sigma^2 \overline{\gamma}_d(v,\alpha),\\
    \mu^\sigma_d(v,\alpha)=&\ c_*\alpha^{-2}-\alpha^{-2}\overline{\mu}_d^0(v)+\sigma^2 \overline{\mu}_d(v,\alpha)\label{eqn:mudrepresentation}
\end{align}
where $\overline{\gamma}_d^0,\overline{\mu}_d^0,\overline{\gamma}_d$ and $\overline{\mu}_d$ are real-valued.

In \ref{app:frozen} we apply It\^o's lemma\footnote{There are several It\^o-type formulas available in the literature, tailored to different solution types. In \ref{app:frozen}, we apply the regular It\^o formula \cite[Theorem 4.32]{daprato} to the weak formulation \eqref{eqn:weak}. Applying the mild It\^o formula \cite[Theorem 1]{mildito} to the mild formulation \eqref{eqn:mild} gives, after tedious computations, an equivalent result.} 
to show that the remainder $v$ defined in \eqref{eqn:remainder2} with $\xi$ and $\alpha$ as in \eqref{eqn:alphaform}-\eqref{eqn:xiform} satisfies the SPDE 
\begin{align}\label{eqn:perturbationgeneral}
    \d v =&  \alpha^{-3}\mathcal{L}_{c_*}v\ \d t + R^\sigma( v,\alpha)\ \d t+ \sigma S(v)[\hat{T}_{\alpha,\xi}\d W_t^Q],
\end{align}
where the drift term $R^\sigma(v,\alpha)$ is of the form 
\begin{align}R^\sigma(v,\alpha)=\alpha^{-3} [N(v)+R_0(v)]+\sigma^2\sum_{i=1}^6 R_i(v,\alpha).\label{eqn:generaldrift}
\end{align} 
Here $N(v)$ is the KdV nonlinearity 
\[N(v)=-\partial_x (v^2)=-2 v\partial_x v,\] 
while $R_0$ through $R_6$ are given by
\begin{align}
R_0(v)=&-\overline{\gamma}_d^0(v) (2+x\partial_x)[\phi_{c_*}+v]-\overline{\mu}_d^0 (v)\partial_x[\phi_{c_*}+v],\nonumber\\
    R_1(v,\alpha)=&\ \tfrac{1}{2}\bigl\| Q^{1/2}\hat{T}^*_{\alpha}\overline{\mu}_s(v)\bigr\|_{\mathcal{H}}^2  \partial_x^2[\phi_{c_*}+v],\nonumber\\
    R_2(v,\alpha)=&\ \bigl\| Q^{1/2}\hat{T}^*_{\alpha}\overline{\gamma}_s(v)\bigr\|_{\mathcal{H}}^2 ( \tfrac{1}{2}x^2\partial_x^2+2  x\partial_x+1)[\phi_{c_*}+v],\nonumber\\
    R_3(v,\alpha)=&\ \bigl\langle Q^{1/2}\hat{T}^*_{\alpha}\overline{\gamma}_s(v),Q^{1/2}\hat{T}^*_{\alpha}\overline{\mu}_s(v) \bigr\rangle_{\mathcal{H}} (x\partial_x^2+2\partial_x)[\phi_{c_*}+v],\nonumber\\
    R_4(v,\alpha)=&-2  M(\phi_{c_*}+v)\bigl[\hat{T}_{\alpha}Q\hat{T}_{\alpha}^*\overline{\gamma}_s(v)\bigr]-  x  M(\partial_x\phi_{c_*}+v_x)\bigl[\hat{T}_{\alpha}Q\hat{T}^*_{\alpha}\overline{\gamma}_s(v)\bigr]  \nonumber\\
&-   \alpha x  M(\phi_{c_*}+v)\bigl[\hat{T}_{\alpha}\hat{\partial}_x Q\hat{T}^*_{\alpha}\overline{\gamma}_s(v)\bigr],\nonumber\\
R_5(v,\alpha)=&- M(\partial_x\phi_{c_*}+v_x)\bigl[\hat{T}_{\alpha}Q\hat{T}^*_{\alpha} \overline{\mu}_s(v)\bigr]-   \alpha M(\phi_{c_*}+v)\bigl[\hat{T}_{\alpha}\hat{\partial}_x Q\hat{T}^*_{\alpha} \overline{\mu}_s(v)\bigr],\nonumber\\
R_6(v,\alpha)=&\ \alpha^{-1}\bigl(\overline{\gamma}_d(v,\alpha) (2+x\partial_x)[\phi_{c_*}+v]+\overline{\mu}_d(v,\alpha)\partial_x[\phi_{c_*}+v]\bigr).\label{eqn:Rterms}
\end{align}
In addition, the operator $S$ in \eqref{eqn:perturbationgeneral} acts on $h\in \mathcal{H}$ as
\begin{align}\label{eqn:martingalepart}
   S(v)[h]=&    M(\phi_{c_*}+v)[h] - (x \partial_x+2)[\phi_{c_*}+v]\langle h,\overline{\gamma_s}(v)\rangle_{\mathcal{H}}\nonumber\\
   &-  \partial_x[\phi_{c_*}+v] \langle h,\overline{\mu}_s(v)\rangle_{\mathcal{H}}.
\end{align}

Note that the noise term in \eqref{eqn:perturbationgeneral} has been transformed via $\hat{T}_{\alpha,\xi}$. The translation invariance of $\hat{T}_{\alpha,\xi}Q\hat{T}^*_{\alpha,\xi}$ (see item 5 of Assumption~\ref{ass:multiplicative}) implies that the transformed noise process $\hat{T}_{\alpha,\xi} W^Q_t$ does, in distribution, not depend on $\xi$. In what follows, we therefore omit the dependence of the noise transformation on $\xi$. We do stress that it should be taken into account during numerical simulations if a pathwise correspondence is desired. 

\subsection{Modulation equations}
\label{subsec:modulation}

The SPDE \eqref{eqn:perturbationgeneral} describes the evolution of the remainder \eqref{eqn:remainder} where the shift $\xi(t)$ and rescaling by $\alpha(t)$ are of the form \eqref{eqn:gammasrepresentation}-\eqref{eqn:mudrepresentation}. This leaves the freedom to make an informed choice for the drift components $\gamma^\sigma_d, \mu^\sigma_d$ and the martingale components $\overline{\gamma}_s, \overline{\mu}_s$ of $\alpha(t)$ and $\xi(t)$. The objective underlying this choice is to ensure that the remainder $v$ remains small in the weighted spaces $L^2_a$ defined in \eqref{eqn:weightedspace}, with $0<a<\sqrt{c_*}$. 

We note that \eqref{eqn:perturbationgeneral} is a non-autonomous semi-linear equation on account of the fact that the linear operator $\mathcal{L}_{c_*}$ carries a $(t,\omega)$-dependent factor $\alpha^{-3}(t,\omega)$. This comes at no surprise in view of the  scaling symmetry \eqref{eqn:symmetry}, where the time-variable receives a factor $\alpha^3$. By transforming the time-variable $t$ path-wise as $(t,\omega)\mapsto \int_0^t \alpha^{-3}(s,\omega)\ \d s$, we can scale out the $(t,\omega)$-dependent factor in \eqref{eqn:perturbationgeneral} as 
\begin{align}\label{eqn:perturbationtransformed}
    \d \tilde{v} =&\ \mathcal{L}_{c_*}\tilde{v}\ \d \tau +\tilde{\alpha}^3 R^\sigma( \tilde{v},\tilde{\alpha})\ \d \tau+ \tilde{\alpha}^{3/2}\sigma S(\tilde{v})[\hat{T}_{\tilde{\alpha}} \d W^Q_\tau],
\end{align}
where $\tilde{v},\tilde{\alpha}$ and $\tilde{\xi}$ are time-transformed versions of $v,\alpha$ and $\xi$. For details, we refer to \cite[Lemma 6.2]{hamsterscalar}, where the same argument is carried out in the setting of reaction-diffusion equations. After this transformation, \eqref{eqn:perturbationtransformed} is semi-linear, and we can recast it into the mild form
\begin{align}\label{eqn:perturbationmild}
    \tilde{v}(\tau) =& \int_0^\tau \tilde{\alpha}^3 e^{\mathcal{L}_{c_*}(\tau-\tau^\prime)} R^\sigma( \tilde{v},\tilde{\alpha})\ \d \tau^\prime+ \int_0^\tau \tilde{\alpha}^{3/2}e^{\mathcal{L}_{c_*}(\tau-\tau^\prime)}\sigma S(\tilde{v}) [\hat{T}_{\tilde{\alpha}} \d W^Q_{\tau^\prime}].
\end{align}
In order to control the drift and martingale components of $\tilde{v}$, and equivalently $v$, we demand that the drift component $R^{\sigma} ( v,\alpha)$
only takes values in the subspace of $L^2_a$ characterized by \eqref{eqn:orthogonality}, where the semigroup generated by $\mathcal{L}_{c_*}$ is contractive. Similarly, we demand that the stochastic integrand $S(v)$ defined in \eqref{eqn:martingalepart} maps $\mathcal{H}$ into this stable subspace. 

\paragraph{Martingale components}

Let us determine what form the martingale components of $\alpha$ and $\xi$ must have to ensure that $S$, as defined in \eqref{eqn:martingalepart}, maps into the stable subspace.  More precisely, we require that
\begin{align*}
\bigl\langle S(v)[h], \zeta_{c_*} \bigr\rangle_{L^2}= \bigl\langle S(v)[h], \phi_{c_*} \bigr\rangle_{L^2}=0,
\end{align*}
for each $v\in L^2(\R)$ and $h\in\mathcal{H} $. This is achieved if $\overline{\gamma}_s,\overline{\mu}_s$ are chosen in such a way that 
\begin{align*}
K(v)
\begin{bmatrix}
\bigl\langle h,\overline{\gamma}_s(v)\bigr\rangle_{\mathcal{H}}\\
\bigl\langle h,\overline{\mu}_s(v)\bigr\rangle_{\mathcal{H}}
\end{bmatrix}=&
\begin{bmatrix}
\bigl\langle
M(v+\phi_{c_*})[h], \phi_{c_*}\bigr\rangle_{L^2}\\
\bigl\langle
M(v+\phi_{c_*})[h], \zeta_{c_*}\bigr\rangle_{L^2}\end{bmatrix}
\end{align*}
holds for all $v\in L^2(\R)$ and $h\in\mathcal{H} $, where $K(v)$ is the matrix
\begin{align}
    K(v)=\begin{bmatrix} \label{eqn:kmatrix}
\bigl\langle (x\partial_x+2)[\phi_{c_*}+v],\phi_{c_*}\bigr\rangle_{L^2}& \bigl\langle  \partial_x v, \phi_{c_*}\bigr\rangle_{L^2}\\
\bigl\langle (x\partial_x+2)[\phi_{c_*}+v],\zeta_{c_*}\bigr\rangle_{L^2} &\bigl\langle \partial_x [\phi_{c_*}+v],\zeta_{c_*}\bigr\rangle_{L^2}
\end{bmatrix}.
\end{align}
The matrix $K(v)$ is invertible in case $\|v\|_{L^2_a}<\delta$ for some $\delta>0$, since $K(v)$ is invertible at $v=0$ and the mapping $v \mapsto \det K(v)$ is continuous from $L^2_a$ to $\R$. We then set
\begin{align}\label{eqn:martingalechoice}
\begin{bmatrix}
\overline{\gamma}_s(v)\\
\overline{\mu}_s(v)
\end{bmatrix}=
 K^{-1}(v)
\begin{bmatrix}
M^*(v+\phi_{c_*})[\phi_{c_*}]\\
 M^*(v+\phi_{c_*})[\zeta_{c_*}]
\end{bmatrix}.
\end{align}

\paragraph{Drift components}

Applying the orthogonality conditions \eqref{eqn:orthogonality} to the drift component $R^\sigma(v,\alpha)$ also gives rises to a system of two linear equations
\begin{align*}
\langle R^\sigma(v,\alpha), \zeta_{c_*} \rangle_{L^2}= \langle R^\sigma(v,\alpha), \phi_{c_*} \rangle_{L^2}=0,
\end{align*}
which is solved by setting
\begin{align}
    \begin{bmatrix}
\overline{\gamma}_d^0(v)\\
\overline{\mu}_d^0(v)
\end{bmatrix}
&=K(v)^{-1}
\begin{bmatrix}
\langle N(v),\phi_{c_*}\rangle_{L^2} \\
\langle N(v),\zeta_{c_*}\rangle_{L^2}
\end{bmatrix},\label{eqn:gammad0}
\end{align}
together with
\begin{align}\label{eqn:driftcomponentgeneral}
\begin{bmatrix}
\overline{\gamma}_d(v,\alpha)\\
\overline{\mu}_d(v,\alpha)
\end{bmatrix}=
-\alpha K^{-1}(v)\sum_{i=1}^5 \begin{bmatrix}
\langle R_i (v,\alpha),\phi_{c_*}\rangle_{L^2}\\
\langle R_i (v,\alpha),\zeta_{c_*}\rangle_{L^2}
\end{bmatrix}.
\end{align}

We collect that the position and scaling processes $\xi$ and $\alpha$ are governed by the modulation system
\begin{align}
    \d v =&\  \alpha^{-3} \mathcal{L}_{c_*}v\ \d t + R^\sigma( v,\alpha)\ \d t+ \sigma S(v) [\hat{T}_{\alpha}\d W^Q_t],\label{eqn:modulationv}\\
     \d \alpha =&\ \bigl[-\alpha^{-2}\overline{\gamma}_d^0(v)+\sigma^2 \overline{\gamma}_d(v,\alpha)\bigr]\ \d t -\sigma \alpha \bigl\langle \hat{T}_{\alpha} \d W_t^Q,\overline{\gamma}_s(v)\bigr\rangle_{\mathcal{H}},\label{eqn:modulationalpha}\\
    \d \xi =&\ \bigl[c_*\alpha^{-2}-\alpha^{-2}\overline{\mu}_d^0(v)+\sigma^2 \overline{\mu}_d(v,\alpha)\bigr] \ \d t -\sigma \alpha \bigl\langle \hat{T}_{\alpha} \d W_t^Q,\overline{\mu}_s(v)\bigr\rangle_{\mathcal{H}},\label{eqn:modulationxi}
\end{align}
and remark that the $v$-dependence on the right-hand side of \eqref{eqn:modulationv} can be summarised as
\[\d v = \alpha^{-3}\bigl[\mathcal{L}_{c_*} v + O(v^2)\bigr]\ \d t+\sigma^2 O(1)\ \d t+ \sigma O(1)\ \hat{T}_{\alpha}\d W_t^Q.\] 
We now return to the examples presented in Sections~\ref{subsubsec:scalar}-\ref{subsubsec:spacetimewhite}, which allows various terms in the modulation system \eqref{eqn:modulationv}-\eqref{eqn:modulationxi} to be simplified. 

\subsubsection{Example I: Modulation equations for scalar noise}\label{subsubsec:scalarnoisemod}
 In the setting of \S\ref{subsubsec:scalar}, the modulation system takes the form 
 \begin{align}\label{eqn:scalarmodulationv}
      \d v=&\ \alpha^{-3} \mathcal{L}_{c_*}v\ \d t +R^{\sigma}_{I}( v,\alpha)\ \d t+ \sigma S_{I}(v) \ \d \beta_t,\\
     \d \alpha =&\ \bigl[-\alpha^{-2} \overline{\gamma}_d^0(v)+\sigma^2 \alpha \overline{\gamma}_{d;I}(v)\bigr]\ \d t -\sigma \alpha \overline{\gamma}_{s;I}(v)\ \d \beta_t \label{eqn:scalarmodulationalpha},\\
    \d \xi =&\ \bigl[c_*\alpha^{-2}-\alpha^{-2} \overline{\mu}_d^0(v)+\sigma^2 \alpha \overline{\mu}_{d;I}(v)\bigr] \ \d t -\sigma \alpha \overline{\mu}_{s;I}(v)\ \d \beta_t. \label{eqn:scalarmodulationxi}
\end{align}
Here the drift component is given by
\[R_I^\sigma(v,\alpha)=\alpha^{-3} \bigl[N(v)+R_0(v)\bigr]+\sigma^2\sum_{i=1}^6 R_{i;I}(v),\] 
where
\begin{align}
    R_{1;I}(v)=&\ \tfrac{1}{2} \overline{\mu}_{s;I}(v)^2  \partial_x^2[\phi_{c_*}+v],\label{eqn:scalardrift1}\\
    R_{2;I}(v)=&\ \overline{\gamma}_{s;I}(v)^2 ( \tfrac{1}{2}x^2\partial_x^2+2  x\partial_x+1)[\phi_{c_*}+v],\label{eqn:scalardrift2}\\
    R_{3;I}(v)=&\ \overline{\gamma}_{s;I}(v)\overline{\mu}_{s;I}(v)  (x\partial_x^2+2\partial_x)[\phi_{c_*}+v],\label{eqn:scalardrift3}\\
    R_{4;I}(v)=&-\bigl(2  \overline{\gamma}_{s;I}(v)(\phi_{c_*}+v)+  x  \overline{\gamma}_{s;I}(v)(\partial_x\phi_{c_*}+v_x)\bigr),\label{eqn:scalardrift4}\\
R_{5;I}(v)=&-  \overline{\mu}_{s;I}(v)(\partial_x\phi_{c_*}+v_x),\label{eqn:scalardrift5}\\
R_{6;I}(v)=&\ \overline{\gamma}_{d;I}(v) (2+x\partial_x)[\phi_{c_*}+v]+\overline{\mu}_{d;I}(v) \partial_x[\phi_{c_*}+v],\label{eqn:scalardrift6}
\end{align}
while the martingale component is given by
\begin{align*}
    S_{I}(v)=&\ \phi_{c_*}+v- 2 \overline{\gamma}_{s;I}(v)[\phi_{c_*}+v]   -\overline{\gamma}_{s;I}(v)x\partial_x[\phi_{c_*}+v]  -   \overline{\mu}_{s;I}(v)\partial_x[\phi_{c_*}+v] .
\end{align*}
The martingale components $\overline{\gamma}_{s;I}$ and $\overline{\mu}_{s;I}$ are real-valued and take the form
\begin{align}\label{eqn:scalargammas}
\begin{bmatrix}
\overline{\gamma}_{s;I}(v)\\
\overline{\mu}_{s;I}(v)
\end{bmatrix}=
K^{-1}(v)
\begin{bmatrix}
\langle v+\phi_{c_*},\phi_{c_*}\rangle_{L^2}\\
\langle v+\phi_{c_*},\zeta_{c_*}\rangle_{L^2}
\end{bmatrix},
\end{align}
where we recall that the matrix $K(v)$ is defined in \eqref{eqn:kmatrix}. The drift components $\overline{\gamma}_{d;I}$ and $\overline{\mu}_{d;I}$ are given by 
\begin{align}
    \begin{bmatrix}
\overline{\gamma}_{d;I}(v)\\
\overline{\mu}_{d;I}(v)
\end{bmatrix}=&-
\overline{\mu}_{s;I}(v)^2
\begin{bmatrix}
\overline{\gamma}_d^1(v)\\
\overline{\mu}_d^1(v)
\end{bmatrix}
- \overline{\gamma}_{s;I}(v)^2
\begin{bmatrix}
\overline{\gamma}_d^2(v)\\
\overline{\mu}_d^2(v)
\end{bmatrix}- \overline{\gamma}_{s;I}(v)\overline{\mu}_{s;I}(v)
\begin{bmatrix}
\overline{\gamma}_d^3(v)\\
\overline{\mu}_d^3(v)
\end{bmatrix}\\
&+ \overline{\gamma}_{s;I}(v) K(v)^{-1}
\begin{bmatrix}
 \bigl\langle (x \partial_x+2)[\phi_{c_*}+v] ,\phi_{c_*}\bigr\rangle_{L^2} \\
 \bigl\langle (x \partial_x+2)[\phi_{c_*}+v] ,\zeta_{c_*}\bigr\rangle_{L^2}
\end{bmatrix}\nonumber\\
&+ \overline{\mu}_{s;I}(v) K(v)^{-1}
\begin{bmatrix}
 \bigl\langle \partial_x[\phi_{c_*}+v],\phi_{c_*}\bigr\rangle_{L^2} \\
 \bigl\langle \partial_x[\phi_{c_*}+v],\zeta_{c_*}\bigr\rangle_{L^2}
\end{bmatrix}  ,\label{eqn:driftcomponentscalar}
\end{align}
where the terms $\overline{\gamma}_d^1,\ldots, \overline{\gamma}_d^3$ and $\overline{\mu}_d^1,\ldots, \overline{\mu}_d^3$ are defined by the expressions
\begin{align}
\begin{bmatrix}
\overline{\gamma}_d^1(v)\\
\overline{\mu}_d^1(v)
\end{bmatrix}&=K(v)^{-1}
\begin{bmatrix}
\tfrac{1}{2}\bigl\langle \partial_x^2[\phi_{c_*}+v], \phi_{c_*}\bigr\rangle_{L^2}\\
\tfrac{1}{2}\bigl\langle \partial_x^2[\phi_{c_*}+v],\zeta_{c_*}\bigr\rangle_{L^2} 
\end{bmatrix},\label{eqn:gammad1}\\
\begin{bmatrix}
\gamma^2_d (v)\\
\overline{\mu}_d^2(v)
\end{bmatrix}
&=K(v)^{-1}
\begin{bmatrix}
\bigl\langle (\tfrac{1 }{2}  x^2 \partial_x^2+2x \partial_x+1)[\phi_{c_*}+v], \phi_{c_*}\bigr\rangle_{L^2}\\
  \bigl\langle (\tfrac{1 }{2}  x^2 \partial_x^2+2x \partial_x+1)[\phi_{c_*}+v], \zeta_{c_*}\bigr\rangle_{L^2}
\end{bmatrix},\label{eqn:gammad2}\\
\begin{bmatrix}
\overline{\gamma}_d^3(v)\\
\overline{\mu}_d^3(v)
\end{bmatrix}
&=K(v)^{-1}
\begin{bmatrix}
 \bigl\langle (x \partial_x^2+2\partial_x)[\phi_{c_*}+v],\phi_{c_*}\bigr\rangle_{L^2} \\
 \bigl\langle (x \partial_x^2+2\partial_x)[\phi_{c_*}+v],\zeta_{c_*}\bigr\rangle_{L^2} 
\end{bmatrix}\label{eqn:gammad3}
.
\end{align}
We remark that $\overline{\gamma}_{d;I}$ and $\overline{\mu}_{d;I}$ induce an $O(\sigma^2)$ drift on $\alpha$ and $\xi$.

\subsubsection{Example I\joinR I: Modulation equations for translation-invariant colored noise}
\label{subsubsec:spacetimecolored2}
In the setting of translation-invariant colored noise of \S\ref{subsubsec:spacetimecolored}, the modulation equations for $v, \alpha$ and $\xi$ take the form
\begin{align*}
    \d v =&\  \alpha^{-3}\mathcal{L}_{c_*}v\ \d t + R_{\RomanII}^\sigma( v,\alpha)\ \d t+ \sigma S_{\diamond}(v)[{T}_{\alpha}\d W_t^Q],\\
\d \alpha =&\ \bigl[-\alpha^{-2}\overline{\gamma}_d^0(v)+\sigma^2 \overline{\gamma}_{d;\RomanII}(v,\alpha)\bigr] \ \d t -\sigma \alpha \bigl\langle {T}_{\alpha}\d W_t^Q,\overline{\gamma}_{\diamond}(v)\bigr\rangle_{L^2},\\
    \d \xi =&\ \bigl[c_*\alpha^{-2}-\alpha^{-2}\overline{\mu}_d^0(v)+\sigma^2 \overline{\mu}_{d;\RomanII}(v,\alpha)\bigr] \ \d t -\sigma \alpha \bigl\langle {T}_{\alpha}\d W_t^Q,\overline{\mu}_{\diamond}(v)\bigr\rangle_{L^2}.
\end{align*}
Here, the martingale component $S_{\diamond}(v)$ acts on $h\in L^2(\R)$ as
\begin{align}\label{eqn:martingalespacetime}
   S_{\diamond}(v)[h]=&\    (\phi_{c_*}+v)h - (x \partial_x+2)[\phi_{c_*}+v]\bigl\langle h,\overline{\gamma}_{\diamond}(v)\bigr\rangle_{L^2}\nonumber\\
   &-  \partial_x[\phi_{c_*}+v] \bigl\langle h,\overline{\mu}_{\diamond}(v)\bigr\rangle_{L^2}.
\end{align}
The functions
$\overline{\gamma}_{\diamond}$ and $\overline{\mu}_{\diamond}$ are $L^2$-valued and given by
\begin{align}\label{eqn:FGL2}
\begin{bmatrix}
\overline{\gamma}_{\diamond}(v)\\
\overline{\mu}_{\diamond}(v)
\end{bmatrix}=K^{-1}(v)
\begin{bmatrix}
 (v+\phi_{c_*})\phi_{c_*}\\
 (v+\phi_{c_*})\zeta_{c_*}\end{bmatrix}.
\end{align}
The drift component $R_{\RomanII}^\sigma(v,\alpha)$ follows from the general formulation \eqref{eqn:generaldrift}, where one evaluates the terms $R_4$ and $R_5$ using \eqref{eqn:transidenties}. These identities also show that the inner products in $R_{1},R_{2}$ and $R_3$ can be computed for $f,g\in L^2(\R)$ as
\begin{align*}
\langle Q^{1/2}\hat{T}^*_{\alpha} f,Q^{1/2}\hat{T}^*_{\alpha}g \rangle_{L^2}=\langle \hat{T}_{\alpha}Q\hat{T}^*_{\alpha} f,g \rangle_{L^2}=
    \alpha^{-1}\langle Q_{\alpha}^{1/2} f,Q_{\alpha}^{1/2}g \rangle_{L^2},
\end{align*}
where $Q_\alpha$ is a rescaled version of the covariance operator $Q$ as introduced in \eqref{eqn:convolutionfamily}. 
The drift components $\overline{\gamma}_{d;\RomanII}$ and $\overline{\mu}_{d;\RomanII}$ are given by the expression
\begin{align}
    \begin{bmatrix}
\overline{\gamma}_{d;\RomanII}(v,\alpha)\\
\overline{\mu}_{d;\RomanII}(v,\alpha)
\end{bmatrix}=&-
\bigl\| Q_{\alpha}^{1/2} \overline{\mu}_{\diamond}(v)\bigr\|^2_{L^2}
\begin{bmatrix}
\overline{\gamma}_d^1(v)\\
\overline{\mu}_d^1(v)
\end{bmatrix}
- \bigl\| Q_{\alpha}^{1/2} \overline{\gamma}_{\diamond}(v)\bigr\|^2_{L^2}
\begin{bmatrix}
\overline{\gamma}_d^2(v)\\
\overline{\mu}_d^2(v)
\end{bmatrix}\nonumber\\
&-\bigl\langle Q_{\alpha}^{1/2} \overline{\gamma}_{\diamond}(v),Q_{\alpha}^{1/2}\overline{\mu}_{\diamond}(v) \bigr\rangle_{L^2}
\begin{bmatrix}
\overline{\gamma}_d^3(v)\\
\overline{\mu}_d^3(v)
\end{bmatrix}\nonumber\\
&+ K(v)^{-1}
\begin{bmatrix}
 \bigl\langle (x\partial_x+2)[\phi_{c_*}+v]Q_{\alpha}\overline{\gamma}_{\diamond}(v) ,\phi_{c_*}\bigr\rangle_{L^2} \\
 \bigl\langle (x\partial_x+2)[\phi_{c_*}+v]Q_{\alpha}\overline{\gamma}_{\diamond}(v) ,\zeta_{c_*}\bigr\rangle_{L^2}
\end{bmatrix}\nonumber\\
&+K(v)^{-1}
\begin{bmatrix}
 \bigl\langle \partial_x[\phi_{c_*}+v]Q_{\alpha}\overline{\mu}_{\diamond}(v),\phi_{c_*}\bigr\rangle_{L^2} \\
 \bigl\langle \partial_x[\phi_{c_*}+v]Q_{\alpha}\overline{\mu}_{\diamond}(v),\zeta_{c_*}\bigr\rangle_{L^2}
\end{bmatrix}\nonumber\\
&+K(v)^{-1} 
\begin{bmatrix}
 \bigl\langle (x Q_{\alpha} \partial_x\overline{\gamma}_{\diamond}(v)+Q_{\alpha} \partial_x \overline{\mu}_{\diamond}(v))(\phi_{c_*}+v),\phi_{c_*}\bigr\rangle_{L^2} \\
 \bigl\langle   (x Q_{\alpha} \partial_x \overline{\gamma}_{\diamond}(v)+Q_{\alpha} \partial_x \overline{\mu}_{\diamond}(v))(\phi_{c_*}+v),\zeta_{c_*}\bigr\rangle_{L^2}
\end{bmatrix},\label{eqn:driftcomponentcolored}
\end{align}
where we recall that the terms $\overline{\gamma}_d^0,\ldots, \overline{\gamma}_d^3$ and $\overline{\mu}_d^0,\ldots, \overline{\mu}_d^3$ are defined in \eqref{eqn:gammad0} and \eqref{eqn:gammad1}-\eqref{eqn:gammad3}.

\subsubsection{Example I\joinR I\joinR I: Modulation equations for space-time white noise}\label{subsubsec:whitenoisemod}
In the setting of space-time noise of \S\ref{subsubsec:spacetimewhite}, the modulation system takes a slightly simpler form, in the sense that the dependence on the rescaling process $\alpha$ is more straightforward:
\begin{align}
\d v =&\  \alpha^{-3} \mathcal{L}_{c_*}v\ \d t + R_{\RomanIII}^\sigma(v,\alpha)\ \d t+ \sigma S_{\diamond}(v) [{T}_{\alpha}\d W_t],\label{eqn:modulationspacetime}\\
\d \alpha =&\ \bigl[-\alpha^{-2}\overline{\gamma}_d^0(v)+\sigma^2 \overline{\gamma}_{d;\RomanIII}(v)\bigr]\ \d t -\sigma \alpha \bigl\langle {T}_{\alpha}\d W_t,\overline{\gamma}_{\diamond}(v)\bigr\rangle_{L^2}\label{eqn:modulationalphaspacetime},\\
    \d \xi =&\ \bigl[c_*\alpha^{-2}-\alpha^{-2}\overline{\mu}_d^0(v)+\sigma^2 \overline{\mu}_{d;\RomanIII}(v)\bigr]\ \d t -\sigma \alpha \bigl\langle {T}_{\alpha}\d W_t,\overline{\mu}_{\diamond}(v)\bigr\rangle_{L^2},\label{eqn:modulationxispacetime}
\end{align}
where
\[R_{\RomanIII}^\sigma(v,\alpha)=\alpha^{-3} \bigl[N(v)+R_0(v)\bigr]+\sigma^2\alpha^{-1}\sum_{i=1}^6 R_{i;\RomanIII}(v),\] 
and
\begin{align*}
    R_{1;\RomanIII}(v)=&\ \tfrac{1}{2}\bigl\| \overline{\mu}_{\diamond}(v)\bigr\|_{\mathcal{H}}^2  \partial_x^2[\phi_{c_*}+v]\\
    R_{2;\RomanIII}(v)=&\ \bigl\|\overline{\gamma}_{\diamond}(v)\bigr\|_{\mathcal{H}}^2 ( \tfrac{1}{2}x^2\partial_x^2+2  x\partial_x+1)[\phi_{c_*}+v]\\
    R_{3;\RomanIII}(v)=&\ \bigl\langle \overline{\gamma}_{\diamond}(v),\overline{\mu}_{\diamond}(v) \bigr\rangle_{\mathcal{H}} (x\partial_x^2+2\partial_x)[\phi_{c_*}+v]\\
    R_{4;\RomanIII}(v)=&-2  (\phi_{c_*}+v)\overline{\gamma}_{\diamond}(v)-  x  (\partial_x\phi_{c_*}+v_x)\overline{\gamma}_{\diamond}(v) \\
    &-   x  (\phi_{c_*}+v){\partial}_x \overline{\gamma}_{\diamond}(v)\\
R_{5;\RomanIII}(v)=&- (\partial_x\phi_{c_*}+v_x)\overline{\mu}_{\diamond}(v)-   (\phi_{c_*}+v){\partial}_x \overline{\mu}_{\diamond}(v)\\
R_{6;\RomanIII}(v)=&\ \overline{\gamma}_{d;\RomanIII}(v) (2+x\partial_x)[\phi_{c_*}+v]+\overline{\mu}_{d;\RomanIII}(v)\partial_x[\phi_{c_*}+v].
\end{align*}
The martingale components $S_{\diamond},\overline{\gamma}_{\diamond}$ and $\overline{\mu}_{\diamond}$ are as in \eqref{eqn:FGL2} and the drift components of $\alpha$ and $\xi$ take the form
\begin{align}
    \begin{bmatrix}
\overline{\gamma}_{d;\RomanIII}(v)\\
\overline{\mu}_{d;\RomanIII}(v)
\end{bmatrix}=&-
\bigl\| \overline{\mu}_{\diamond}(v)\bigr\|^2_{L^2}
\begin{bmatrix}
\overline{\gamma}_d^1(v)\\
\overline{\mu}_d^1(v)
\end{bmatrix}- \bigl\|  \overline{\gamma}_{\diamond}(v)\bigr\|^2_{L^2}
\begin{bmatrix}
\overline{\gamma}_d^2(v)\\
\overline{\mu}_d^2(v)
\end{bmatrix}\nonumber\\&-\bigl\langle \overline{\gamma}_{\diamond}(v),\overline{\mu}_{\diamond}(v) \bigr\rangle_{L^2}
\begin{bmatrix}
\overline{\gamma}_d^3(v)\\
\overline{\mu}_d^3(v)
\end{bmatrix}\nonumber\\
&+ K(v)^{-1}
\begin{bmatrix}
 \bigl\langle (x\partial_x+2)[\phi_{c_*}+v]\overline{\gamma}_{\diamond}(v) ,\phi_{c_*}\bigr\rangle_{L^2} \\
 \bigl\langle (x\partial_x+2)[\phi_{c_*}+v]\overline{\gamma}_{\diamond}(v) ,\zeta_{c_*}\bigr\rangle_{L^2}
\end{bmatrix}\nonumber\\
&+K(v)^{-1}
\begin{bmatrix}
 \bigl\langle \partial_x[\phi_{c_*}+v]\overline{\mu}_{\diamond}(v),\phi_{c_*}\bigr\rangle_{L^2} \\
 \bigl\langle \partial_x[\phi_{c_*}+v]\overline{\mu}_{\diamond}(v),\zeta_{c_*}\bigr\rangle_{L^2}
\end{bmatrix}\nonumber\\
&+K(v)^{-1} 
\begin{bmatrix}
 \bigl\langle (x \partial_x \overline{\gamma}_{\diamond}(v)+\partial_x \overline{\mu}_{\diamond}(v))(\phi_{c_*}+v),\phi_{c_*}\bigr\rangle_{L^2} \\
 \bigl\langle   (x \partial_x \overline{\gamma}_{\diamond}(v)+\partial_x \overline{\mu}_{\diamond}(v))(\phi_{c_*}+v),\zeta_{c_*}\bigr\rangle_{L^2} \label{eqn:gammadIII}
\end{bmatrix}.
\end{align}
We recall again that the terms $\overline{\gamma}_d^0,\ldots, \overline{\gamma}_d^3$ and $\overline{\mu}_d^0,\ldots, \overline{\mu}_d^3$ are defined in \eqref{eqn:gammad0} and \eqref{eqn:gammad1}-\eqref{eqn:gammad3}.
\subsection{Numerical simulations}
\label{subsec:numerics}
Having fully laid out the modulation systems in the stochastic co-moving frame for our three example setups, we now explore the dynamics that the systems produce via numerical simulations. We restrict ourselves to simulations of scalar noise (Example I) and space-time white noise (Example I\joinR I\joinR I), as the modulation systems for these examples are the most tractable. See \ref{app:schemes} for the numerical schemes that were employed to simulate the stochastic KdV equation \eqref{eqn:SKDV}  and the modulation systems  \eqref{eqn:scalarmodulationv}-\eqref{eqn:scalarmodulationxi} and \eqref{eqn:modulationalphaspacetime}-\eqref{eqn:modulationxispacetime}. 

\paragraph{Pathwise simulation}

It is widely acknowledged that accurate numerical simulations of SPDEs are challenging to obtain \cite{lord}. In order to validate our results, we present a numerical comparison between our constructed modulation parameters and soliton parameters derived from a direct simulation. Figure~\ref{fig:pathwisesoliton} shows one realization obtained from a simulation of the KdV equation \eqref{eqn:scalar} with multiplicative scalar noise in both the original frame and co-moving frame. 

\begin{figure*}[h!]
    \centering
    \begin{subfigure}[t]{0.5\textwidth}
        \centering
        \includegraphics[height=2.3in]{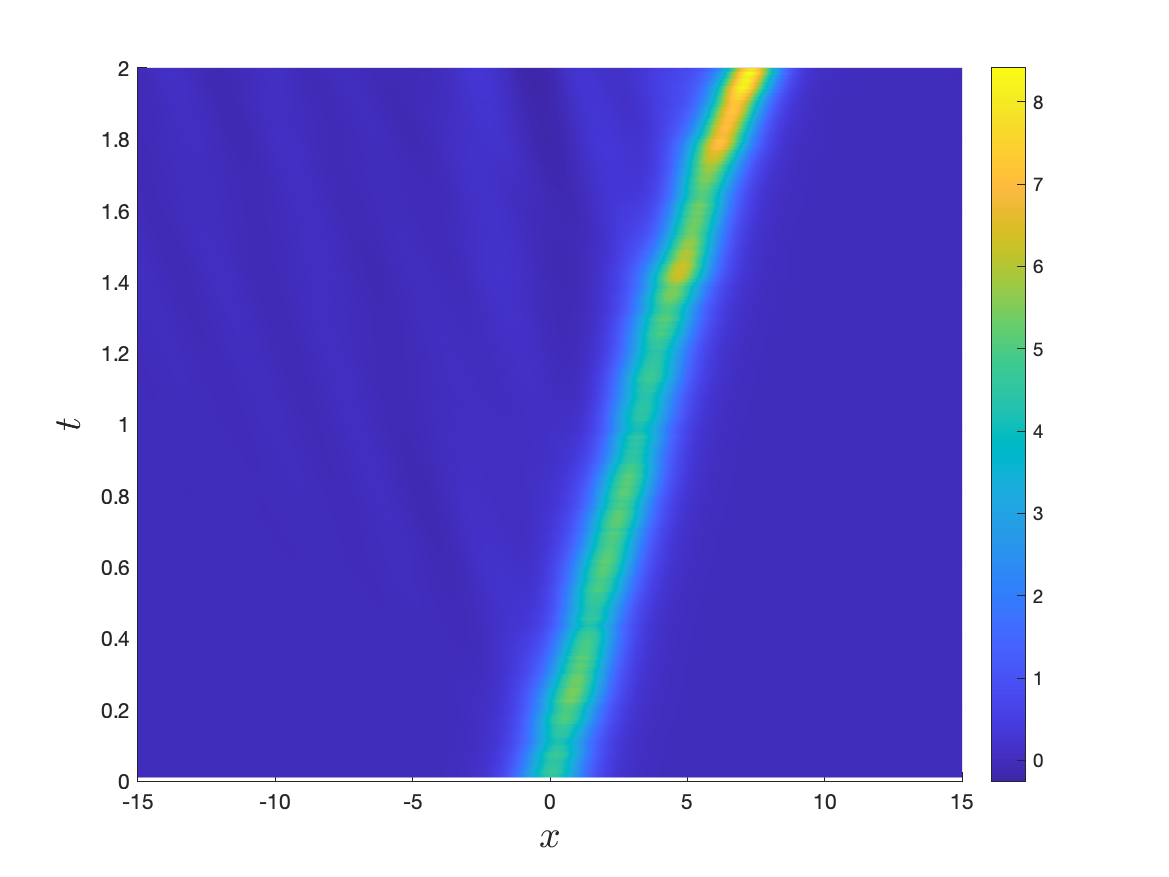}
        \caption{Original frame.}\label{subfig:pathwisesolitonoriginal}
    \end{subfigure}%
    ~ 
    \begin{subfigure}[t]{0.5\textwidth}
        \centering
        \includegraphics[height=2.3in]{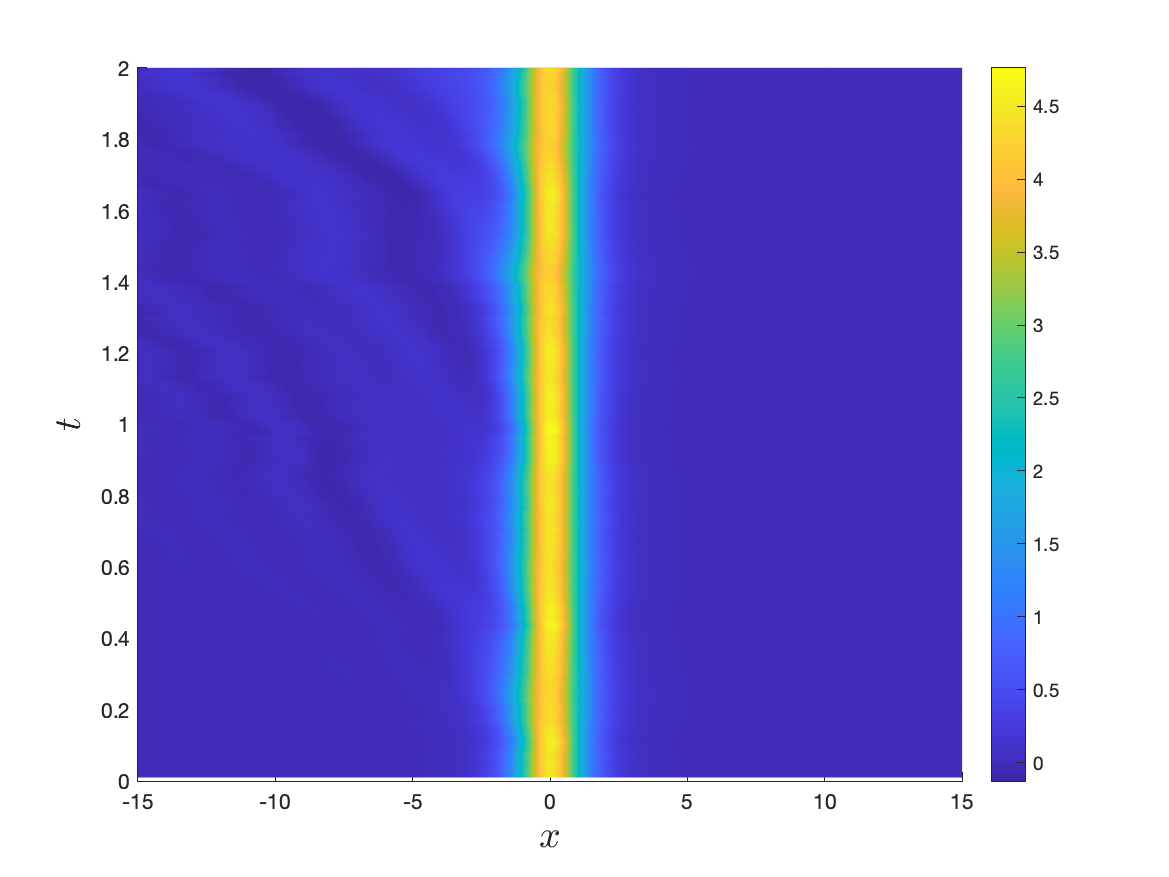}
        \caption{Stochastic co-moving frame.}\label{subfig:pathwisesolitonfrozen}
    \end{subfigure}
        \begin{subfigure}[t]{0.5\textwidth}
        \centering
        \includegraphics[height=2.3in]{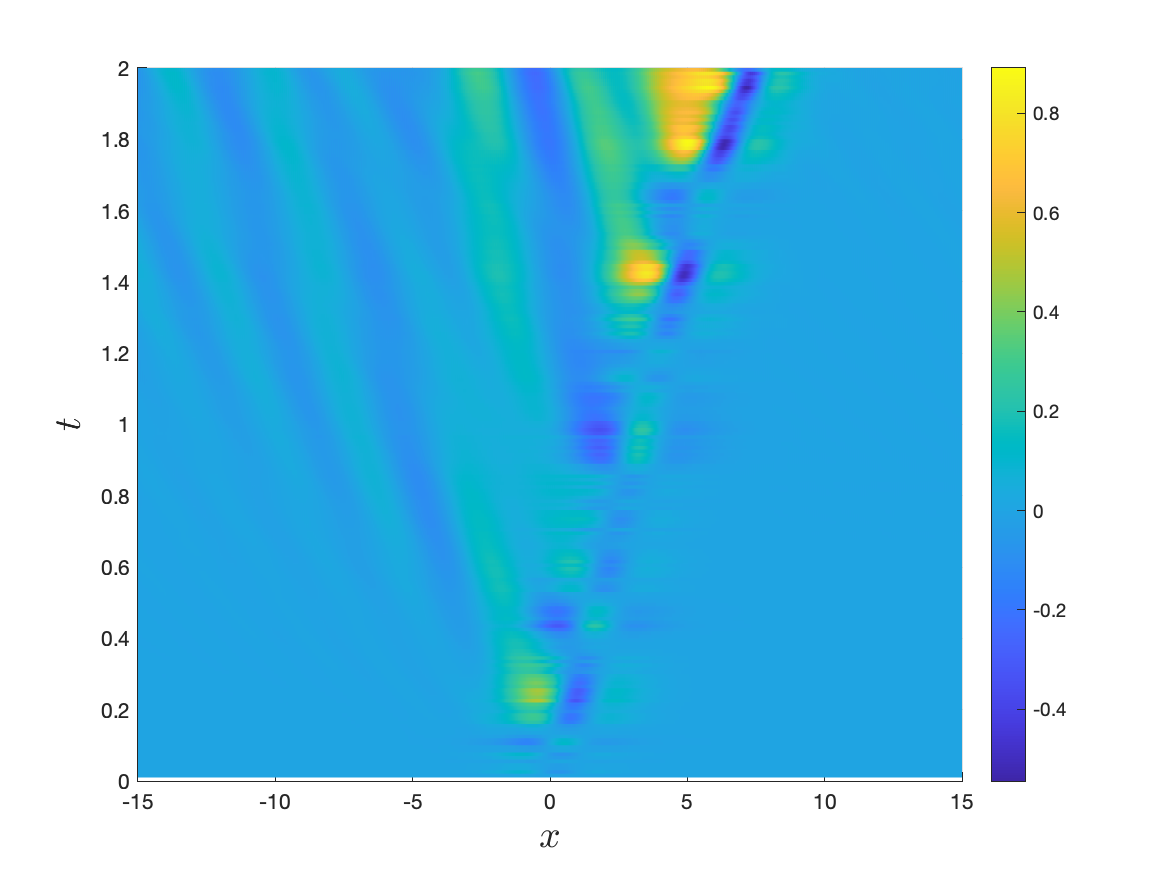}
        \caption{Original frame, soliton removed. }\label{subfig:originalp}
    \end{subfigure}%
    ~ 
    \begin{subfigure}[t]{0.5\textwidth}
        \centering
        \includegraphics[height=2.3in]{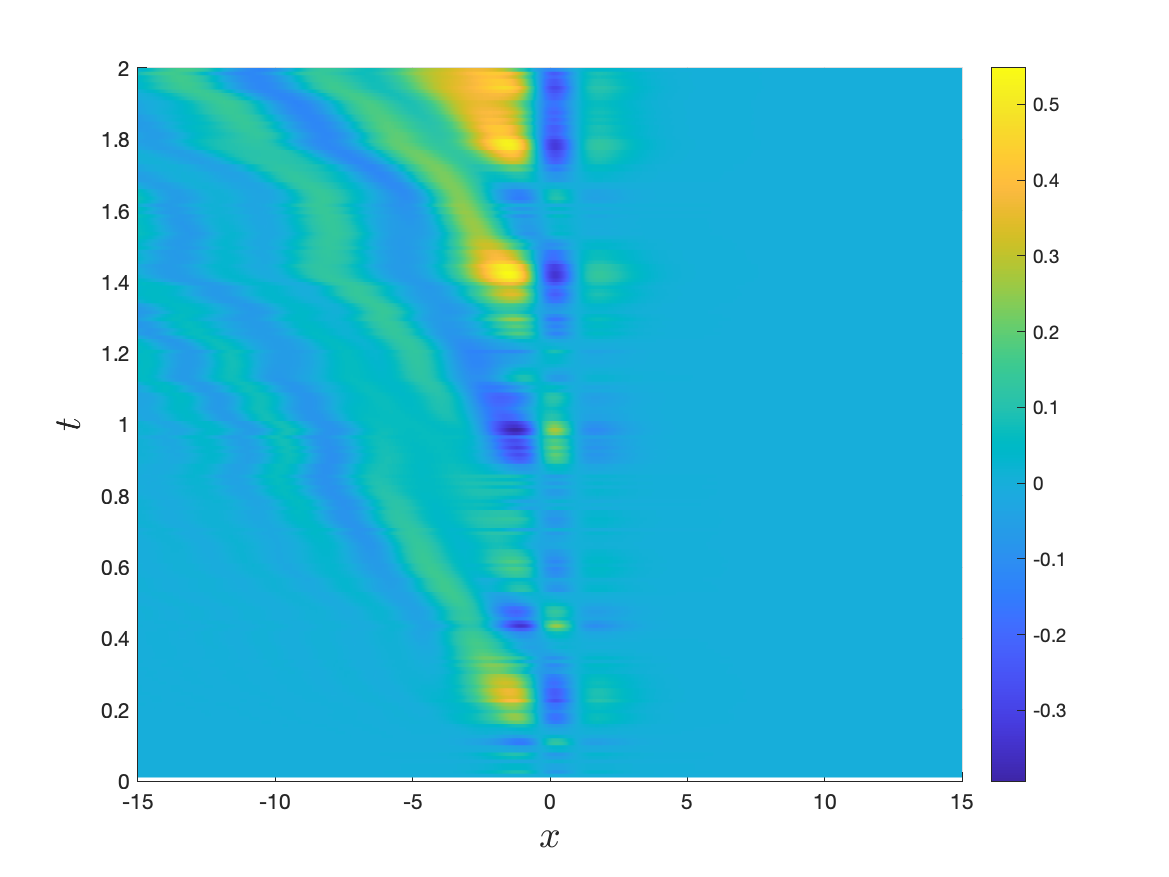}
        \caption{Stochastic co-moving frame, soliton removed.}\label{subfig:frozenp}
    \end{subfigure}
    \caption{Simulation of the KdV equation with scalar noise of strength $\sigma=0.25$. Panel  (a) shows the original frame realization $u(t,x)$, from a simulation of \eqref{eqn:scalar}. Panel (b) shows $ \phi_{c_*}(x)+v(t,x)$, from simulation in the frozen frame of \eqref{eqn:scalarmodulationv}-\eqref{eqn:scalarmodulationxi} with the same realization of the noise. Panels (c) and (d) show the perturbation with respect to the soliton, that is $u(t,x)-\phi_{c(t)}(x-\xi(t))$ with the phase-definitions \eqref{eqn:phasedefinitions} in panel (c), and $v(t,x)$ in panel (d).}\label{fig:pathwisesoliton}
\end{figure*}

In Figure~\ref{subfig:pathwisesolitonoriginal}, we observe that the soliton propagates approximately at a constant velocity, and at times slightly speeds up or slows down when it increases or decreases in amplitude, respectively. The transformation from the stochastic KdV equation \eqref{eqn:scalar} to the modulation system \eqref{eqn:scalarmodulationv}-\eqref{eqn:scalarmodulationxi} allows us to `freeze' the stochastic soliton. In Figure~\ref{subfig:pathwisesolitonfrozen}, the soliton remains centered and roughly has constant amplitude. To the left of the soliton we observe slight perturbations due to the noise. These can be observed more clearly upon removing the soliton in Figure~\ref{subfig:originalp}, which reveals the `wake' of the stochastic soliton. The stochastic perturbations encountered by the soliton result in a radiation field to the left of the soliton. In Figure~\ref{subfig:frozenp}, we furthermore observe that the radiation field has undergone a rescaling in the $x$-direction, as is evident from the distortion of the radiation waves. The effect of the stochastic frozen-frame transformation can also be visualised in the case of space-time white noise, see Figure~\ref{fig:pathwiseperturbation} in \ref{app:supplementary}.

For comparison purposes, we define `fitted' versions of the position $\xi_{\text{fit}}(t)$ and amplitude $c_{\text{fit}}(t)$ of a solution $u(t,x)$
to \eqref{eqn:skdvgeneral} implicitly via the identities 
\begin{align}
    \bigl\langle u\bigl(t,\cdot+\xi_{\text{fit}}(t)\bigr)-\phi_{c_{\text{fit}}(t)},\zeta_{c_{\text{fit}}(t)}\bigr\rangle_{L^2}&=0, \label{eqn:phasedefinitions}\\
    \bigl\langle u\bigl(t,\cdot+\xi_{\text{fit}}(t)\bigr)-\phi_{c_{\text{fit}}(t)},\phi_{c_{\text{fit}}(t)}\bigr\rangle_{L^2}&=0,\nonumber
\end{align}
which we solve numerically.
This allows us to compare the evolution of soliton parameters obtained by direct simulation of \eqref{eqn:skdvgeneral} and the modulation system \eqref{eqn:modulationv}-\eqref{eqn:modulationxi}. Recall therefore that the amplitude process $c(t)$ can be recovered from the rescaling process $\alpha(t)$ as $c(t)=c_*\alpha^{-2}(t)$. We also introduce the phase shift processes
\begin{align}
    \Omega(t)=\xi(t)-\int_0^t c(s) \ \d s \quad \text{and} \quad \Omega_{\text{fit}}(t)=\xi_{\text{fit}}(t)-\int_0^t c_{\text{fit}}(s) \ \d s,\label{eqn:phaseshift}
\end{align}
which track the deviation of the soliton position from the integrated stochastic velocities $c(t)$ and $c_{\text{fit}}(t)$ and isolate the noise-induced effects.

\begin{figure*}[h!]
    \centering
    \begin{subfigure}[t]{0.5\textwidth}
        \centering
        \includegraphics[height=2.1in]{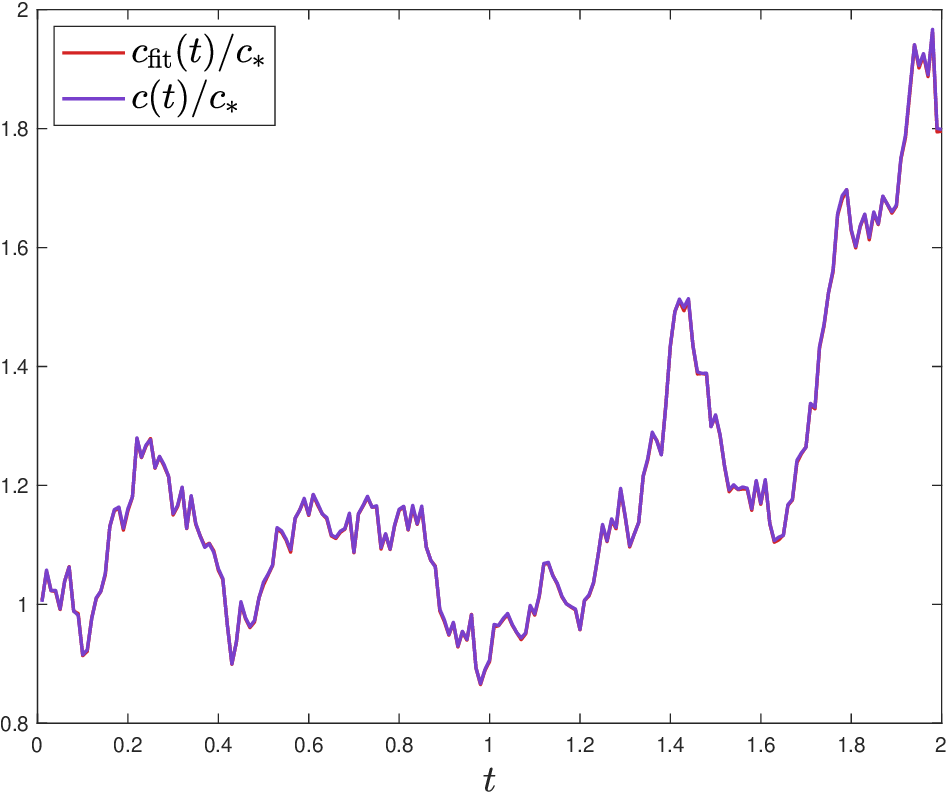}
    \end{subfigure}%
    ~ 
    \begin{subfigure}[t]{0.5\textwidth}
        \centering
        \includegraphics[height=2.1in]{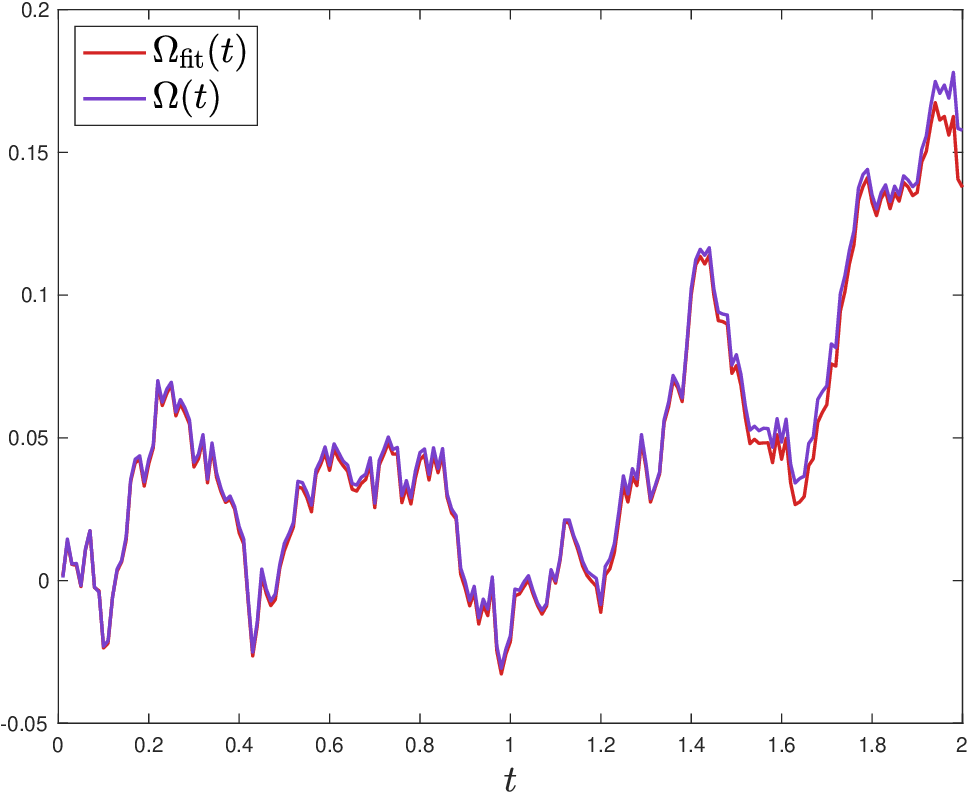}
    \end{subfigure}
    \caption{Path-wise comparison of soliton amplitudes $c(t)$ to $c_{\text{fit}}(t)$ (left) and phase shifts $\Omega(t)$ to $\Omega_{\text{fit}}(t)$ (right) at noise strength 
 $\sigma=0.25$ and initial amplitude $c_*=3$. The parameters $c_{\text{fit}}(t)$ and $\Omega_{\text{fit}}(t)$, defined in \eqref{eqn:phasedefinitions} and \eqref{eqn:phaseshift}, are obtained from direct simulation in the original frame of \eqref{eqn:scalar}. The soliton amplitude $c(t)$ and phase shift $\Omega(t)$ are obtained from simulation of the frozen frame system \eqref{eqn:scalarmodulationv}-\eqref{eqn:scalarmodulationxi}.}\label{fig:pathwiseparameters}
\end{figure*}

Figure~\ref{fig:pathwiseparameters} shows the correspondence between the evolution of the soliton amplitude and phase shift in both frames. Note that the soliton amplitude in this realization attains almost twice its original value at $t=2$. The phase shifts $\Omega$ and $\Omega_{\text{fit}}$ develop a small discrepancy over time, which we attribute mainly to
truncation effects and the fact that
errors in $c_{\text{fit}}$ are compounded through the integral
in \eqref{eqn:phaseshift}.

\begin{figure*}[h!]
\centering
    \begin{subfigure}[t]{0.5\textwidth}
    \centering
        \includegraphics[height=2in]{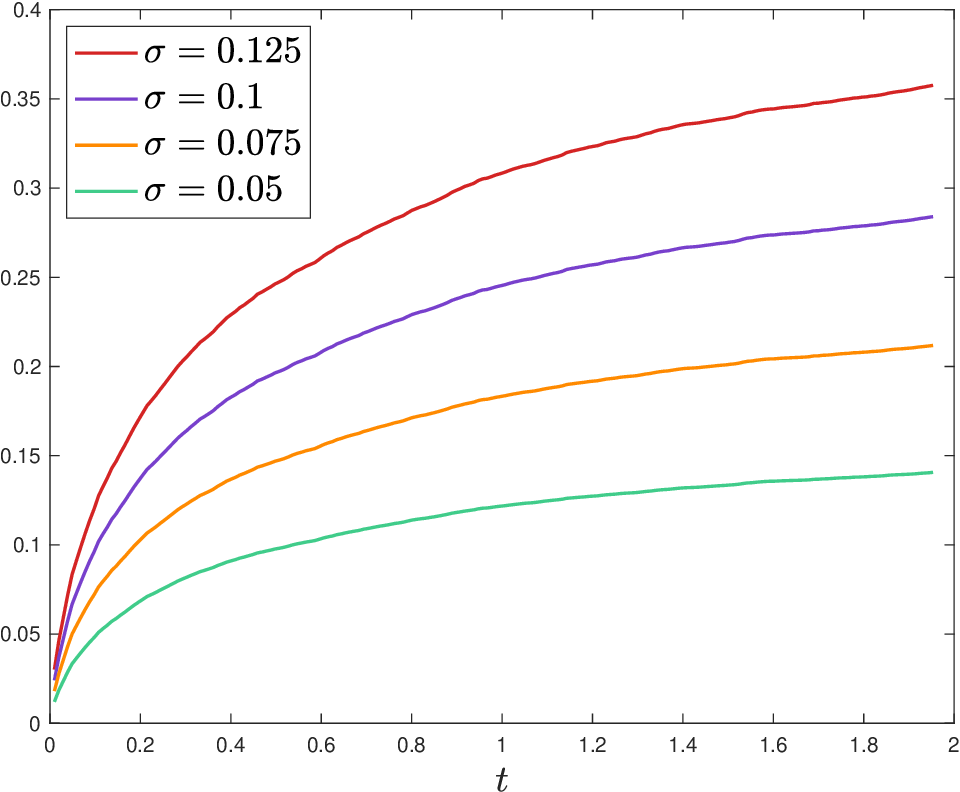}
        \caption{Perturbation size over time.}
    \end{subfigure}%
    ~
    \begin{subfigure}[t]{0.5\textwidth}
    \centering
        \includegraphics[height=2in]{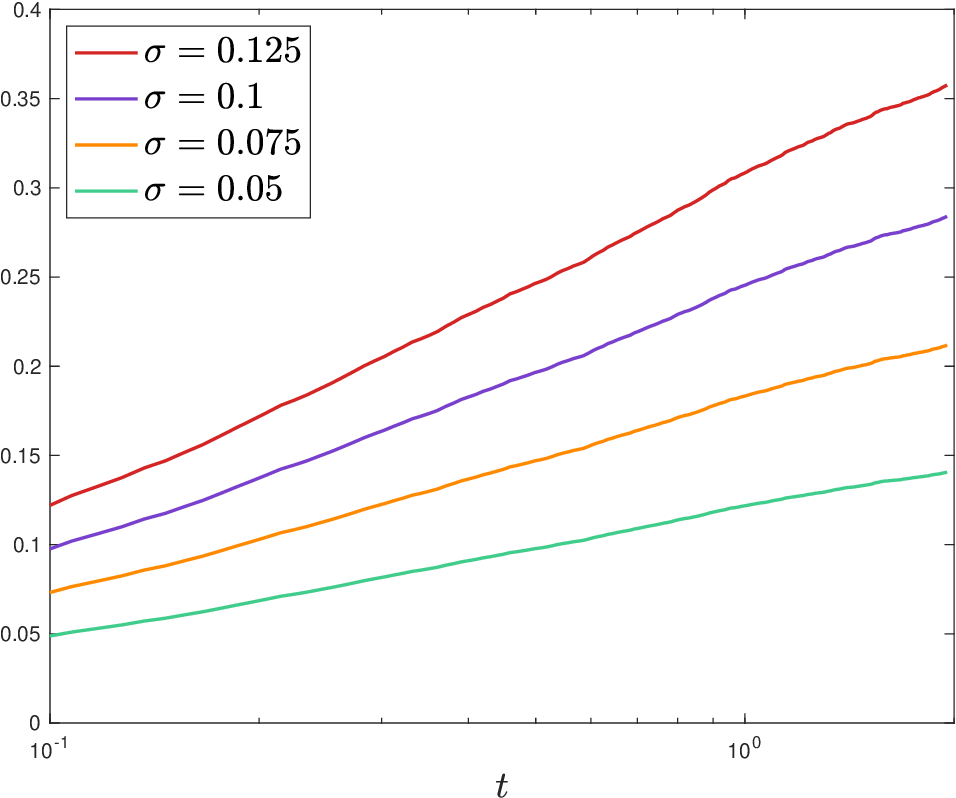}
        \caption{Perturbation size over time, log-scale.}\label{subfig:perturbationlog}
        \end{subfigure}
        \caption{Sample mean of the process $\sup_{s\leq t}\|v(s)\|_{L_a^2([-50,20])}$ for scalar noise, see \eqref{eqn:scalarmodulationv},
        computed over $500$ realisations for $\sigma\in\{0.05,0.075,0.1,0.125\}$ and $c_*=3$.  The exponential weight $e^{ax}$ in the $L^2_a$-norm strongly amplifies numerical effects entering from the right boundary of the computational domain $[-50,50]$. We take care to avoid these by computing the $L^2_a$-norm on $[-50,20]$, with $a=0.5$. For the initial soliton-parameter used in this simulation, the relevant dynamics occur well within $[-50,20]$ (see Figure~\ref{fig:pathwisesoliton}).}\label{fig:perturbationscalar}
\end{figure*}

\paragraph{Stability}
The construction of the modulation system in \S\ref{subsec:modulation} should ensure that the perturbation $v$ remains small in the exponentially weighted spaces $L^2_a$ defined in \eqref{eqn:weightedspace}. Figure~\ref{fig:perturbationscalar} shows the average growth of the in $L_a^2$-norm of the perturbation $v$ with respect to the soliton. The spatial norm of the perturbation appears to grows logarithmically, as indicated by Figure~\ref{subfig:perturbationlog}, where we observe a linear growth of the perturbation size on logarithmic scale. For the case of  space-time white noise we refer to Figure~\ref{fig:perturbationwhite}. 

The logarithmic growth strongly suggests that our soliton-tracking method is valid over exponentially long timescales. A logarithmic growth of the remaining perturbation is also observed in \cite[Figure 3.8]{hamster}, where the spirit of our approach is applied to traveling waves in reaction-diffusion equations. Using this fact, Hamster and the second author rigorously prove in \cite{hamsterstability} that the exit-time from the soliton family is exponentially long with respect to the parameter $1/\sigma$.

\section{Soliton dynamics}
\label{sec:solitondynamics}

In this section, we set out to derive explicit, tractable expansions to uncover 
the effects of the multiplicative noise on the soliton amplitude $c$ and position $\xi$. In \S\ref{sec:tracking}, we have seen that the dynamics of the soliton parameters $c$ and $\xi$ are governed by  the rescaling process $\alpha$ and the infinite dimensional perturbation $v$, which follow the coupled equations
\begin{align}
    v(t) =&\ I^\sigma_v(v,\alpha,t),\nonumber\\
    \alpha(t)=&\ I^\sigma_\alpha(v,\alpha,t).\label{eqn:integralformalpha}
\end{align}
For $v$, we choose to work with the mild formulation
\begin{align}\label{eqn:vintegral}
I^\sigma_v(v,\alpha,t)=& \int_0^t  e^{\int_s^t \alpha^{-3}(t^\prime)\d t^\prime \mathcal{L}_{c_*}} R^\sigma (v,\alpha)\ \d s\\
&+ \sigma \int_0^t e^{\int_s^t \alpha^{-3}(t^\prime)\d t^\prime \mathcal{L}_{c_*}} S(v)[\hat{T}_{\alpha}\d W_s^Q],\nonumber
\end{align}
which follows from \eqref{eqn:perturbationmild} by undoing the time transformation. This has the advantage of being suitable for constructing explicit approximations. For $I_\alpha^\sigma$ we use the strong form
\begin{align}
I^\sigma_\alpha(v,\alpha,t)=&\ 1+\int_0^t \bigl[-\alpha^{-2}\overline{\gamma}_d^0(v)+\sigma^2 \overline{\gamma}_d(v,\alpha)\bigr] \ \d s \nonumber\\
&-\sigma \int_0^t \alpha \langle \hat{T}_{\alpha}\d W_s^Q,\overline{\gamma}_s(v)\rangle_{\mathcal{H}},\label{eqn:strongalpha}
\end{align}
which corresponds to \eqref{eqn:scalarmodulationalpha}.

To develop our approximation procedure, it is relevant to note that the process $\alpha$ exhibits significant fluctuations, while the perturbation $v$ remains relatively small due to the damping of the semigroup. Indeed, we observe that $\alpha$ grows as $O(\sigma \sqrt{t})$, which can be anticipated by noting that $I^\sigma_\alpha$ contains no damping terms. On the other hand, Figure~\ref{fig:perturbationscalar} indicates that $v$ grows at a slower rate, namely $O(\sigma \ln t)$. As such, a carefully tailored approximation procedure is required to accurately capture the interplay between $\alpha$ and $v$. To construct approximations of $\alpha$ that account for the influence of the perturbation $v$, we introduce SDEs based on an expansion of the $v$-dependent coupling terms. In broad terms, we will expand the $\alpha$ dynamics in terms of $v$, while expanding the $v$ dynamics in terms of $\sigma$, treating $\alpha$ as an external input.

Below, in \S\ref{subsec:expansion} - \S\ref{subsec:expper}, we describe the expansions of the coupled system \eqref{eqn:integralformalpha} in more detail. 
Combining these, we obtain a sequence of increasingly refined approximations for the modulation parameters, which we present in \S\ref{subsec:combined}.  
We proceed by evaluating the first few approximations for the cases of Example I (scalar noise) and Example I\joinR I\joinR I (space-time white noise) in \S\ref{subsec:scalarnoiseapprox} and \S\ref{subsec:whitenoiseapprox}, respectively. In the first setting, the soliton amplitude roughly behaves as a geometric Brownian motion. The explicit approximations are used to compute leading-order statistical properties, which we compare with sample statistics of the numerical observations $c_{\text{fit}}$ and $\Omega_{\text{fit}}$ defined in \eqref{eqn:phasedefinitions} and \eqref{eqn:phaseshift}.

\subsection{Expansion of the rescaling process}
\label{subsec:expansion}

In order to unravel how the rescaling process $\alpha$ is influenced by the perturbation $v$, we introduce an expansion of $\alpha$ in terms of $v$. As a first step, consider the situation where the perturbation $v$ is set to zero in \eqref{eqn:integralformalpha}. That is, we introduce a process $A_0$ which satisfies
\begin{align}\label{eqn:alphastar}
A_0(t) = I^\sigma_\alpha(0,A_0,t)
\end{align}
and we supply this SDE with the initial condition $A_0(0)=1$. The process $A_0$ then constitutes a relatively crude first approximation to $\alpha$, which we  subsequently refine by increasing the order of the perturbation $v$ that we take into account. 

In particular, let us now include terms in the It\^o form $I_\alpha^\sigma(v,\alpha,t)$ that depend \textit{linearly} on $v$. We assume that there is an approximation of $v$ available, for which we introduce the variable $\tilde{v}_1$. This variable should be thought of as approximating $v$ with an error of $O(\sigma^2)$. The next approximation $A_1$, given the process $\tilde{v}_1$, is defined through the SDE
\[A_1(\tilde{v}_1,t)=I^\sigma_\alpha(0 ,A_1,t)+[I^\sigma_\alpha]^{(1)}(\tilde{v}_1;A_1,t).\]
Here, $[I^\sigma_\alpha]^{(1)}$ is defined as
\begin{align*}
[I^\sigma_\alpha]^{(1)}(v;\alpha,t)&=\sigma^2\int_0^t  [\overline{\gamma}_d]^{(1)}(v;\alpha) \ \d s-\sigma \int_0^t \alpha \bigl\langle \hat{T}_{\alpha}\d W_s^Q,[\overline{\gamma}_s]^{(1)}(v)\bigr\rangle_{\mathcal{H}},
\end{align*}
with $[\overline{\gamma}_d]^{(1)}$ and $[\overline{\gamma}_s]^{(1)}$ denoting the linear parts of the mappings $v\mapsto \overline{\gamma}_d(v,\alpha)$ and $v\mapsto \overline{\gamma}_s(v)$. We remark that the functional $\overline{\gamma}_d^0(v)$ in \eqref{eqn:strongalpha} contains no linear part, and is therefore not included in the definition of $[I^\sigma_\alpha]^{(1)}$.

In general, if $f$ is a map from a Banach space $X$ into a Banach space $Y$ that is $N+1$ times differentiable at $v=0$, we write 
\[f(v)=\sum_{k=0}^N[f]^{(k)}(v)+O\bigl(\|v\|^{N+1}_{X}\bigr)\] 
where $[f]^{(k)}(v)\sim v^k$ is the symmetric $k$-linear map that collects the order $k$ powers of $v$ in $f(v)$. Alternatively, one can say that $[f]^{(k)}(v)$ denotes the order $k$ term in the Taylor expansion of $v\mapsto f(v)$ around zero.

This expansion procedure extends naturally to higher orders. For the next approximation, we also include quadratic terms. Furthermore, we base this approximation on an additional variable $\tilde{v}_2$, which should be thought of as an approximation to $v$ with error $O(\sigma^3)$. Given two  processes $\tilde{v}_1,\tilde{v}_2$, we define $A_2$ as the solution to
\begin{align*}
    A_2(\tilde{v}_1,\tilde{v}_2,t)=I^\sigma_\alpha(0,A_2,t)+[I^\sigma_\alpha]^{(1)}(\tilde{v}_2;A_2,t)+[I^\sigma_\alpha]^{(2)}(\tilde{v}_1;A_2,t).
\end{align*}

In general, $A_k$ is defined implicitly in terms of the processes $\tilde{v}_1,\ldots,\tilde{v}_k$ via the SDE
\begin{align}\label{eqn:alphahigher}
A_k(\tilde{v}_1,\ldots,\tilde{v}_k,t)
=&I^\sigma_\alpha(0,A_k,t)+\sum_{i=1}^k[I^\sigma_\alpha]^{(i)}(\tilde{v}_{\lfloor \frac{k}{i} \rfloor};A_k,t)
\end{align}
with
\begin{align*}
[I^\sigma_\alpha]^{(k)}(v;\alpha,t)=&\int_0^t \bigl[-\alpha^{-2}[\overline{\gamma}_d^0]^{(k)}(v)+\sigma^2 [\overline{\gamma}]^{(k)}_d(v;\alpha)\bigr] \ \d s\\
&-\sigma \int_0^t \alpha \bigl\langle \hat{T}_{\alpha}\d W_s^Q,[\overline{\gamma}_s]^{(k)}(v)\bigr\rangle_{\mathcal{H}}.
\end{align*}

\subsection{Expansion of the phase shift}
Following the same procedure, we expand the modulation equation for the position $\xi$ in terms of $v$. In \S\ref{sec:tracking} we have seen that $\xi$ can be recovered from the identity
\begin{align*}
    \xi(t)&=\int_0^t c(s) \ \d s+I^\sigma_\Omega(v,\alpha,t),
\end{align*}
where
\begin{align*}
    I^\sigma_\Omega(v,\alpha,t) =& -\int_0^t \alpha^{-2}\overline{\mu}_d^0(v)\ \d s+\sigma^2 \int_0^t\overline{\mu}_d(v,\alpha)\ \d s -\sigma \int_0^t\alpha \bigl\langle \hat{T}_{\alpha} \d W_s^Q,\overline{\mu}_s(v)\bigr\rangle_{\mathcal{H}}.
\end{align*}
Note that the position primarily follows the velocity $c(t)$, with additional noise-induced corrections resulting in the phase shift
\begin{align}
    \Omega(t)=\xi(t)-\int_0^t c(s) \ \d s=I^\sigma_\Omega(v,\alpha,t).
\end{align}
Analogously to \eqref{eqn:alphastar} and \eqref{eqn:alphahigher}, we define approximations to the phase shift $\Omega(t)$ as
\begin{align*}
 \bar{\Omega}_0(\alpha,t)=&\ I_\Omega^\sigma(0,\alpha,t),
\end{align*}
and for $k\geq1$
\begin{align*}
 \bar{\Omega}_k(\tilde{v}_1,\ldots,\tilde{v}_k,\alpha,t)=&\ I^\sigma_\Omega(0,\alpha,t)+\sum_{i=1}^k[I^\sigma_\Omega]^{(i)}\bigl(\tilde{v}_{\lfloor \frac{k}{i} \rfloor};\alpha,t\bigr).
\end{align*}

\subsection{Expansion of the perturbation}
\label{subsec:expper}

We now turn to the complementary problem and examine how the perturbation $v$ depends on the rescaling process $\alpha$. Treating $\alpha$ as an input, we expand the perturbation $v$ in terms of the small parameter $\sigma$ as
\begin{align}\label{eqn:formalv}
V_k(\alpha,t)=\sigma V^{(1)}(\alpha,t)+\ldots +\sigma^k V^{(k)}(\alpha,t) 
\end{align}
based on the integral form \eqref{eqn:integralformalpha}. Here, $V^{(k)}$ collects all terms of $O(\sigma^k)$ in the It\^o form \eqref{eqn:integralformalpha}. Collecting all $O(\sigma)$ terms in $I_v^\sigma$, gives
\[V^{(1)}(\alpha,t)=\sigma^{-1}I_v^\sigma(0,\alpha,t)=\int_0^t e^{\int_s^t \alpha^{-3}(t^\prime)\d t^\prime \mathcal{L}_{c_*}} S(0)[\hat{T}_{\alpha}\d W_s^Q].\]

In order to find the subsequent term $V^{(2)}$ in the expansion \eqref{eqn:formalv}, we note first that the drift component in the It\^o form \eqref{eqn:vintegral} satisfies $R^\sigma(v,\alpha)=O(v^2+\sigma^2)$. Consequently, we can explicitly define $V^{(2)}$ in \eqref{eqn:formalv} using $V^{(1)}$. Indeed, collecting the $O(\sigma^2)$ terms in $I_v^\sigma$, we arrive at 
\begin{align*}  V^{(2)}(\alpha,t)=&\int_0^t  \alpha^{-3} e^{\int_s^t \alpha^{-3}(t^\prime)\d t^\prime \mathcal{L}_{c_*}}  \Bigl[N\bigl( V^{(1)}(\alpha,s)\bigr)+[R_0]^{(2)}\bigl( V^{(1)}(\alpha,s)\bigr)\Bigr]\ \d s\\
&+\int_0^t e^{\int_s^t \alpha^{-3}(t^\prime)\d t^\prime \mathcal{L}_{c_*}}\sum_{i=1}^6 R_i(0,\alpha)\ \d s\\
&+  \int_0^t e^{\int_s^t \alpha^{-3}(t^\prime)\d t^\prime \mathcal{L}_{c_*}} [S]^{(1)} \bigl(V^{(1)}(\alpha,s)\bigr)[\hat{T}_{\alpha}\d W_s^Q].
\end{align*}

Any subsequent term $V^{(k)}$ in \eqref{eqn:formalv} can now be found by continuing systematically. In general, we have
\begin{align*}  V^{(k)}(\alpha,t)=&\int_0^t  \alpha^{-3} e^{\int_s^t \alpha^{-3}(t^\prime)\d t^\prime \mathcal{L}_{c_*}}  \Bigl[\chi(k)N\bigl( V^{(\frac{k}{2})}(\alpha,s)\bigr)+\sum_{i\mid k}[R_0]^{(i)}\bigl( V^{(\frac{k}{i})}(\alpha,s)\bigr)\Bigr]\ \d s\\
&+\int_0^t e^{\int_s^t \alpha^{-3}(t^\prime)\d t^\prime \mathcal{L}_{c_*}}\sum_{i=1}^6 \sum_{j\mid(k-2)}[R_i]^{(j)}\bigl(V^{(\frac{k-2}{j})}(\alpha,s);\alpha\bigr)\ \d s\\
&+  \int_0^t e^{\int_s^t \alpha^{-3}(t^\prime)\d t^\prime \mathcal{L}_{c_*}} \sum_{i\mid(k-1)}[S]^{(i)} \bigl(V^{(\frac{k-1}{i})}(\alpha,s)\bigr)[\hat{T}_{\alpha}\d W_s^Q],
\end{align*}
where $\chi(k)=1$ if $k$ is even and $\chi(k)=0$ otherwise.

\subsection{The combined system approximation}\label{subsec:combined}
We now combine the expansion of $\alpha$ in $v$ and the expansion of $v$ in $\sigma$ to construct our full approximations to the coupled system  \eqref{eqn:integralformalpha}. We define for $k\geq 0$ approximations $\alpha_k$ to $\alpha$ as
\[\alpha_k(t)=A_k\bigl(V_1(\alpha_{k-1},\cdot),\dots,V_k(\alpha_{k-1},\cdot) ,t\bigr),\]
and approximations $\Omega_k$ to $\Omega$ as
\[\Omega_k(t)=\bar{\Omega}_k\bigl(V_1(\alpha_{k},\cdot),\dots,V_k(\alpha_{k},\cdot),\alpha_k ,t\bigr).\]
For $k\geq 1$ we furthermore introduce the approximations
    \[v_k(t)={V}_k(\alpha_{k-1},t),\]
to $v$.

The soliton amplitude directly follows from the rescaling process $\alpha$ via the relation $c(t)=c_*\alpha^{-2}(t)$. We can therefore define approximations $c_0,c_1,c_2,\ldots$ to $c$ by directly writing
\[c_k=c_*\alpha_k^{-2} \quad \text{for}\quad  k\geq 0.\]

We now examine what these approximation constructions produce for the examples discussed in \S\ref{sec:tracking}. 

\subsection{Example I: Soliton dynamics for scalar noise}\label{subsec:scalarnoiseapprox}
We first turn to the setting of Example I outlined in \S\ref{subsubsec:scalar}. Here, $v, \alpha$ and $\xi$ follow the modulation system \eqref{eqn:scalarmodulationv}-\eqref{eqn:scalarmodulationxi}. We observe that the fully decoupled approximation in this setting satisfies the SDE
\begin{align*}
    \d \alpha_0 =&\ \overline{\gamma}_{d;I}(0)\sigma^2 \alpha_0\  \d t-\overline{\gamma}_{s;I}(0)\sigma \alpha_0\  \d \beta_t \\
   =&\ (\tfrac{74}{135}+\tfrac{4\pi^2}{405})\sigma^2 \alpha_0\  \d t-\tfrac{2}{3}\sigma \alpha_0 \ \d \beta_t
\end{align*}
with $\alpha_0(0)=1$. Here we have used Table~\ref{tab:constants} to evaluate the constants. This SDE admits the explicit solution
\[\alpha_0(t)=e^{(\overline{\gamma}_{d;I}(0)-\frac{1}{2}\overline{\gamma}^2_{s;I}(0))\sigma^2t-\overline{\gamma}_{s;I}(0)\sigma \beta_t}=e^{(\frac{44}{135}+\frac{4 \pi^2}{405})\sigma^2 t-\frac{2}{3}\sigma \beta_t},\]
which is a geometric Brownian motion.

At first order, the perturbation is given by 
\begin{align*}
    v_1(t) =&\ \sigma \int_0^t e^{\int_s^t \alpha_0^{-3}(t^\prime)\d t^\prime \mathcal{L}_{c_*}} S_I(0)\ \d \beta_s
\end{align*}
where we remark that
\[S_{I}(0)=- \tfrac{1}{3}\phi_{c_*}   -\tfrac{2}{3}x\partial_x \phi_{c_*}  +\tfrac{2}{3} c_*^{-1/2} \partial_x \phi_{c_*} .\]

The subsequent approximation $\alpha_1$ to the rescaling process is the geometric Brownian motion
\[\d \alpha_{1} =\sigma^2K^{1,1}_{I}(t) \alpha_1 \ \d t-\sigma K^{2,1}_{I}(t)   \alpha_1 \ \d\beta_t,\]
with random coefficients $K^{1,1}_{I}(t),K^{2,1}_{I}(t)$ that are given explicitly by
\begin{align*}
    K^{1,1}_{I}(t)=&\overline{\gamma}_{d;I}(0)+[\overline{\gamma}_{d;I}]^{(1)}\bigl(v_1(t)\bigr), \\
    K^{2,1}_{I}(t)=&\overline{\gamma}_{s;I}(0)+[\overline{\gamma}_{s;I}]^{(1)}\bigl(v_1(t)\bigr).
\end{align*}

The second order approximation for the perturbation is given explicitly by 
\[v_2(t)=\sigma V^{(1)}(\alpha_1,t)+\sigma^2 V^{(2)}(\alpha_1,t),\]
using 
\[V^{(1)}(\alpha,t)=\int_0^t e^{\int_s^t \alpha^{-3}(t^\prime)\d t^\prime \mathcal{L}_{c_*}} S_I(0)\ \d \beta_s\]
and
\begin{align*}
    V^{(2)}(\alpha,t)=&\int_0^t  \alpha^{-3}(s) e^{\int_s^t \alpha^{-3}(t^\prime)\d t^\prime \mathcal{L}_{c_*}}  \Bigl[N\bigl( V^{(1)}(\alpha,s)\bigr)+[R_0]^{(2)}\bigl( V^{(1)}(\alpha,s)\bigr)\Bigr]\ \d s\\
&+\int_0^t e^{\int_s^t \alpha^{-3}(t^\prime)\d t^\prime \mathcal{L}_{c_*}}\sum_{i=1}^6 R_{i;I}(0)\ \d s\\
&+ \int_0^t e^{\int_s^t \alpha^{-3}(t^\prime)\d t^\prime \mathcal{L}_{c_*}} [S_I]^{(1)} \bigl(V^{(1)}(\alpha,s)\bigr)\ \d \beta_s.
\end{align*}
This approximation collects the leading-order drift effects in the perturbation $v$. The sample mean of $v$ and $v_2$ are displayed in Figure~\ref{fig:meanv}, which shows the development of an average radiation field induced by the noise. These sample means are computed by first subtracting the process $v_1$, which has mean zero. This eliminates the leading-order fluctuations and speeds up convergence to the mean. In contrast to wave profiles in the stochastic (FitzHugh)-Nagumo equations analysed in \cite{hamster}, the average perturbation does not localise around the wave profile and does not seem to converge in time. Rather, the noise leads to an average pattern of radiation waves that expands far behind the soliton.

 \begin{figure*}[h!]
 \centering
    \begin{subfigure}[t]{0.5\textwidth}
        \centering
        \includegraphics[height=2in]{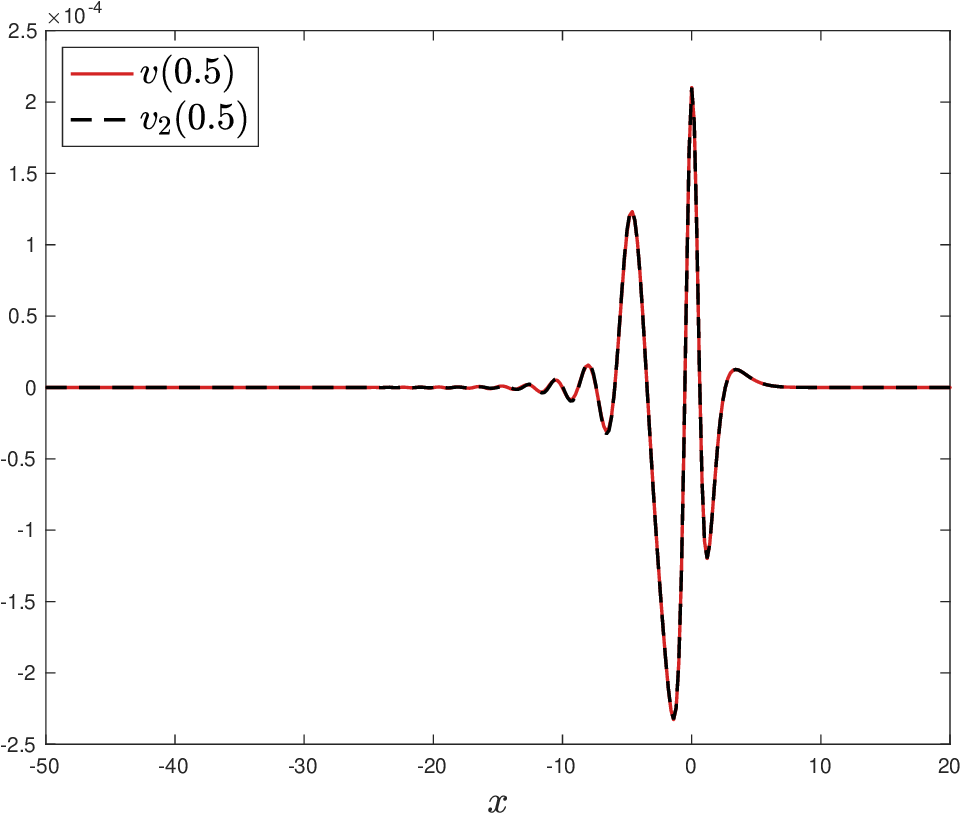}
    \end{subfigure}%
    ~
    \begin{subfigure}[t]{0.5\textwidth}
        \centering
        \includegraphics[height=2in]{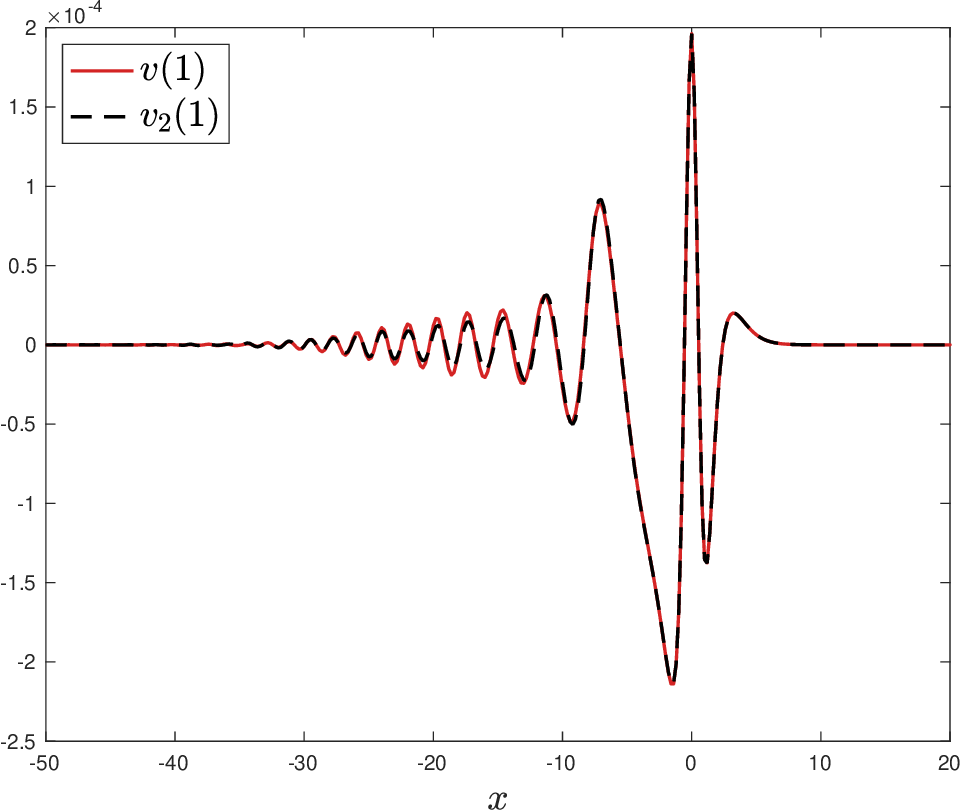}
    \end{subfigure}%
    \\
    \begin{subfigure}[t]{0.5\textwidth}
        \centering
        \includegraphics[height=2in]{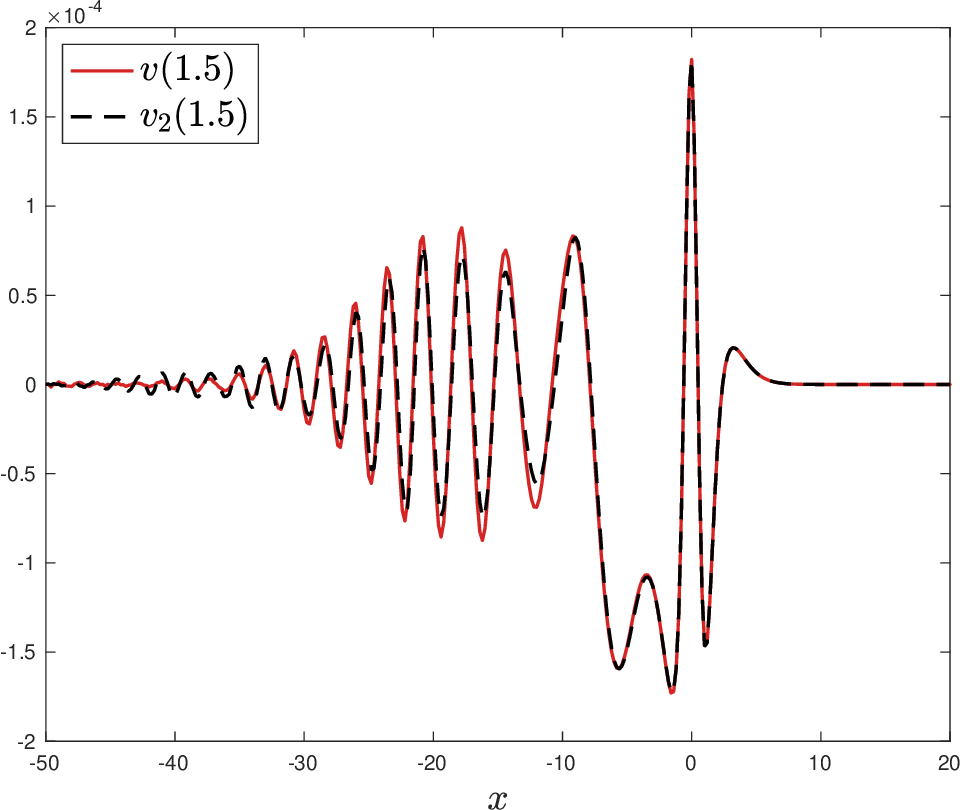}
    \end{subfigure}%
    ~
    \begin{subfigure}[t]{0.5\textwidth}
        \centering
        \includegraphics[height=2in]{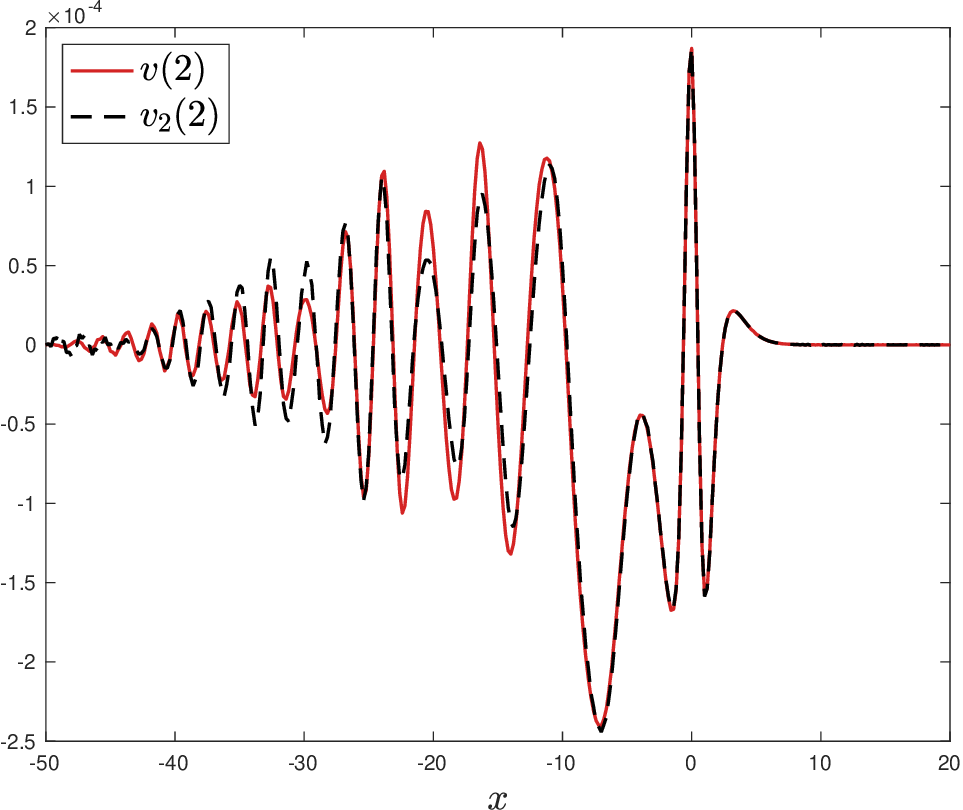}
    \end{subfigure}%
    \caption{Sample mean of $v$ (dashed) and the approximation $v_2$ (solid) as the perturbation develops between $t=0.5$ and $t=2$. Computed over 3000 realisations for $\sigma=0.03$.} \label{fig:meanv}
\end{figure*}

Using $v_2$, the approximation $\alpha_2$ is defined through the scalar SDE
\begin{align*}\d \alpha_2=&\ [- K^{0,2}_{I}(t)\alpha_2^{-2}+\sigma^2K^{1,2}_{I}(t)  \alpha_2] \ \d t-\sigma K^{2,2}_{I}(t)   \alpha_2 \ \d \beta_t 
\end{align*}
with random coefficients 
\begin{align*}
    K^{0,2}_{I}(t)=&\ [\overline{\gamma}_d^0]^{(2)}\bigl(\sigma V^{(1)}(\alpha_1,t)\bigr), \\
    K^{1,2}_{I}(t)=&\ \overline{\gamma}_{d;I}(0)+[\overline{\gamma}_{d;I}]^{(1)}\bigl(v_2(t)\bigr) +[\overline{\gamma}_{d;I}]^{(2)}\bigl(\sigma V^{(1)}(\alpha_1,t)\bigr), \\
    K^{2,2}_{I}(t)=&\ \overline{\gamma}_{s;I}(0)+[\overline{\gamma}_{s;I}]^{(1)}\bigl(v_2(t)\bigr)+[\overline{\gamma}_{s;I}]^{(2)}\bigl(\sigma V^{(1)}(\alpha_1,t)\bigr).
\end{align*}

\paragraph{Amplitude}

The first approximation for the soliton amplitude $c(t)$ is the geometric Brownian motion
\begin{align*}
    c_0(t)=c_* \alpha_0^{-2}(t) =c_*e^{-(\frac{88}{135}+\frac{8 \pi^2}{405})\sigma^2 t+\frac{4}{3}\sigma \beta_t}.
\end{align*}
We remark that for small noise strengths $\sigma$, the dynamics of $c_0(t)$ are largely determined by the factor $e^{\frac{4}{3}\sigma \beta_t}$. Note that this factor is also present in \eqref{eqn:heuristic}, which heuristically explains how the leading-order stochastic dynamics of $c(t)$ arise.
Using the exact expression
\begin{align}\label{eqn:scalarvariance}
\operatorname{Var}[c_0(t)]=c_*^2e^{(\frac{64}{135}-\frac{16\pi^2}{405})\sigma^2 t}(e^{\frac{16}{9}\sigma^2 t}-1),
\end{align}
we compare the variance of $c_0(t)$ to the sample variance of $c_{\text{fit}}(t)$ in Figure~\ref{fig:scalaramplitudevar}. Although the approximation $c_0$ is fully decoupled from the perturbation $v$, we see that its variance already agrees quite well with that of $c(t)$.

\begin{figure*}[h!]
        \centering
 \begin{subfigure}[t]{0.5\textwidth}
     \centering
        \includegraphics[height=2in]{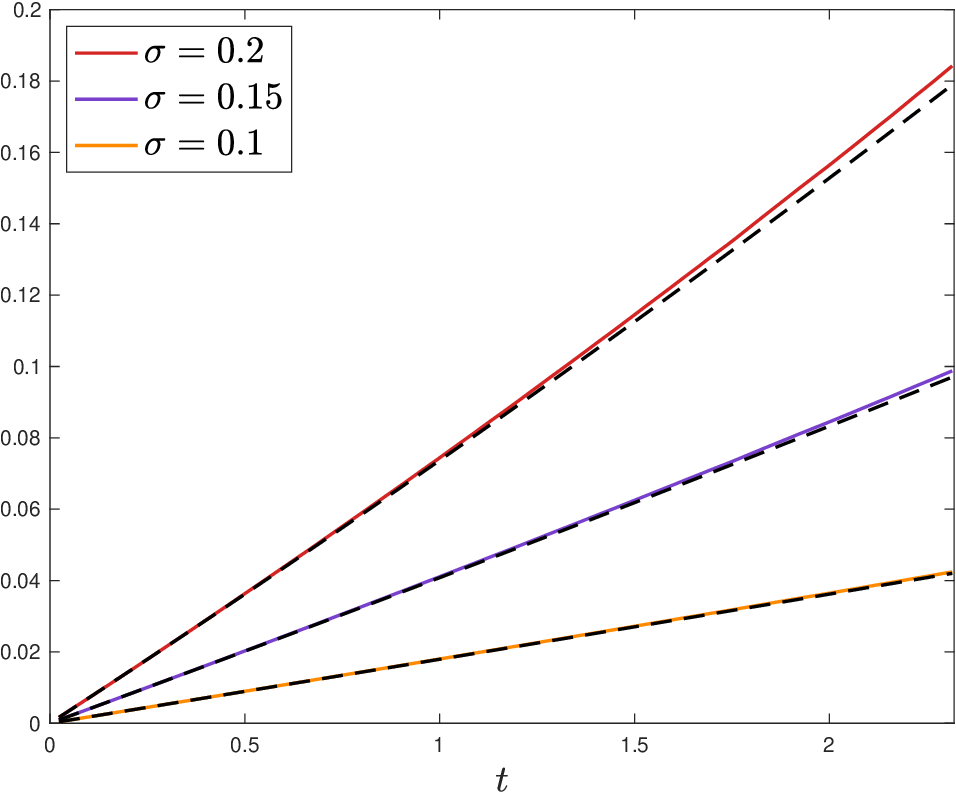}
    \caption{Scalar noise.}\label{fig:scalaramplitudevar}
\end{subfigure}%
~
\begin{subfigure}[t]{0.5\textwidth}
    \centering
    \includegraphics[height=2in]{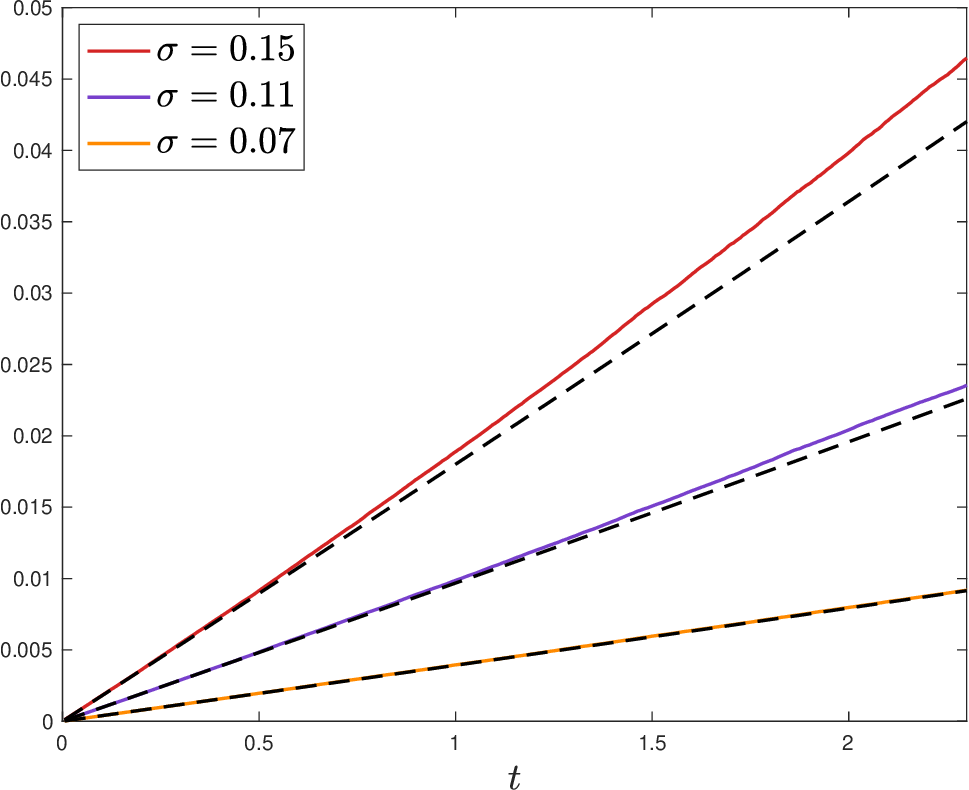}
    \caption{Space-time white noise.}\label{fig:whitenoiseamplitudevariance}
\end{subfigure}
\caption{Sample variance of the process $c_{\text{fit}}(t)/c_*$ for scalar noise and space-time white noise at various noise strengths $\sigma$. Solid lines indicate the sample variance, dashed lines indicate the variance of $c_0(t)$ as in \eqref{eqn:scalarvariance} and \eqref{eqn:whitenoisevariance}, respectively. The sample variance is computed over $204\cdot10^4$ and  $8\cdot10^4$ realizations, respectively.}
\end{figure*}

 Figure~\ref{fig:scalaramplitudemean} compares the mean of $c_{\text{fit}}(t)$ with that of the increasingly refined approximations $c_0(t), c_1(t)$ and $c_2(t)$, using
 \begin{align}
     \mathbb{E}[c_0(t)]=c_*e^{(\frac{32}{135}-\frac{8\pi^2}{405})\sigma^2 t}\label{eqn:meanc0}
 \end{align}
 for $c_0(t)$. The means of $c_1(t)$ and $c_2(t)$ are not as easily computed analytically due to the dependence on $\alpha$ and $v$ in the random coefficients $K^{0,2}_{I}(t),K^{1,2}_{I}(t)$ and $K^{2,2}_{I}(t)$. We therefore consider their sample means. We see in Figure~\ref{fig:scalaramplitudemean} that the mean of $c_{\text{fit}}(t)$ is not well-approximated by that of $c_0(t)$ or $ c_1(t)$. The sample mean of $c_2(t)$, however, agrees well with that of  $c_{\text{fit}}(t)$, indicating that the quadratic terms of $v$ contribute significantly to the evolution of the mean amplitude. 


 \begin{figure*}[h!]
 \centering

    \begin{subfigure}[t]{0.5\textwidth}
        \centering
        \includegraphics[height=2in]{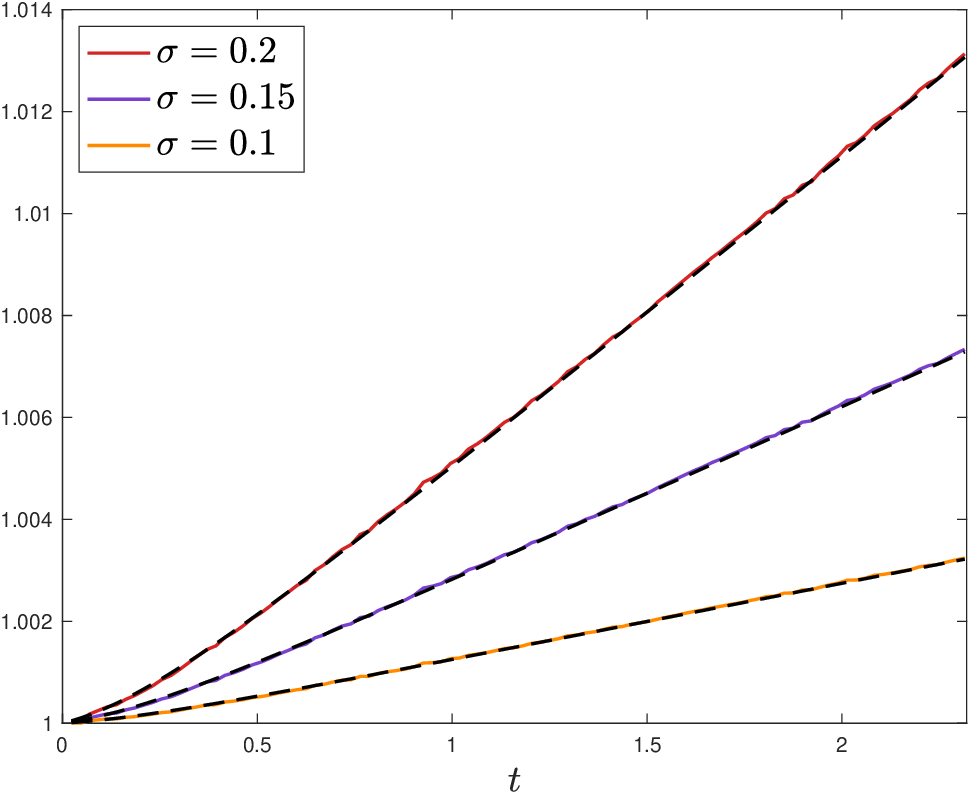}
        \caption{Mean of $c_{\text{fit}}(t)/c_*$ at various noise strengths.}
    \end{subfigure}%
    ~
    \begin{subfigure}[t]{0.5\textwidth}
        \centering
        \includegraphics[height=2in]{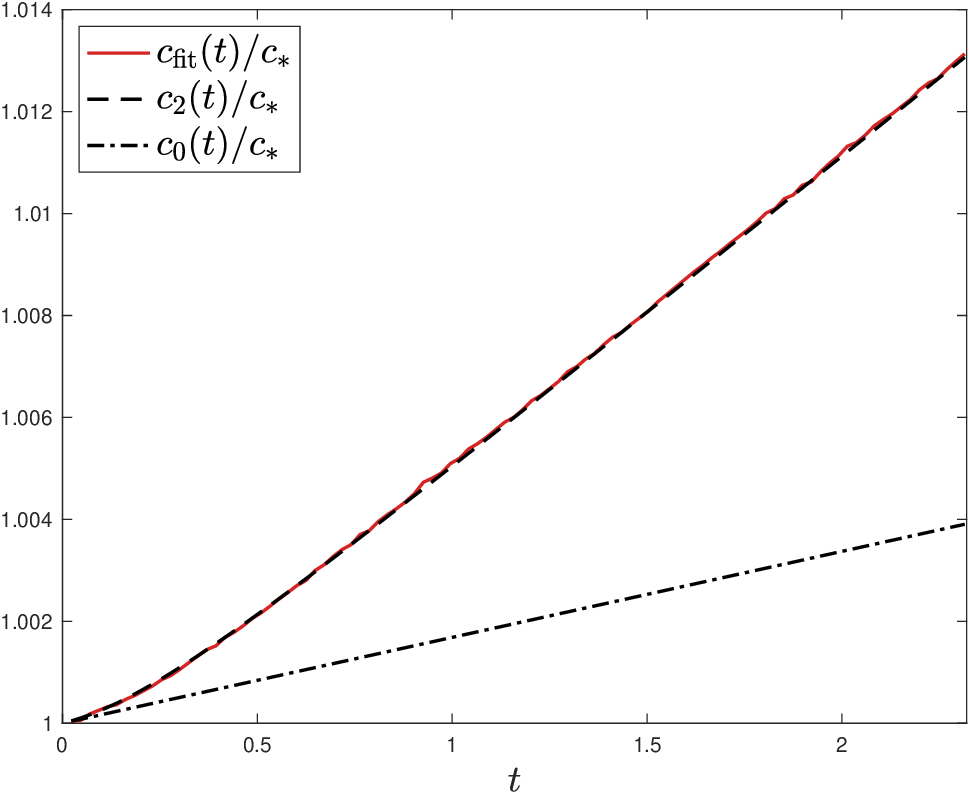}
        \caption{Means of $c_{\text{fit}}(t)/c_*$ and approximations $c_0(t)$ and $c_2(t)$.}
    \end{subfigure}%
    \caption{Sample mean of the process $c_{\text{fit}}(t)/c_*$ computed over $2\cdot10^4$ realisations for scalar noise. The mean of $c_2(t)$ (dashed) agrees well with that of $c_{\text{fit}}(t)$ at the simulated noise strength values $\sigma\in\{0.1,0.15,0.2\}$, whereas the mean of $c_0(t)$ (dash dot) as in \eqref{eqn:meanc0} fails to capture the correct amplitude drift. The mean of $c_1(t)$ provides no improvement, it can not be distinguished from that of $c_0(t)$ at these simulation values.}\label{fig:scalaramplitudemean}
\end{figure*}

\paragraph{Phase shift}
The first approximation to the phase shift process $\Omega(t)$ defined in \eqref{eqn:phaseshift} is given by
\begin{align*}
    \Omega_0(t) =&\ \overline{\mu}_{d;I}(0)\sigma^2 \int_0^t \alpha_0(s) \ \d s +\tfrac{2}{3}c_*^{-1/2}\sigma \int_0^t \alpha_0(s) \ \d \beta_s.
\end{align*}
We remark that numerical computations of $\Omega_{\text{fit}}(t)$ are unsuitable for ensemble simulations, due to a large truncation error (see Figure~\ref{fig:pathwiseparameters}). We therefore consider ensemble simulations of the more robust process $\Omega(t)$. 
Figure~\ref{fig:xivariancescalar} shows that the soliton position, on average, develops a phase lag from the velocity $c(t)$ which appears to grow almost linearly in time. The mean of the lowest approximation $\Omega_0(t)$ is evaluated as
\begin{align}
    \mathbb{E}[\Omega_0(t)]=&\ \overline{\mu}_{d;I}(0)\sigma^2 \int_0^t \mathbb{E}\bigl[\alpha_0(s)\bigr]\ \d s \nonumber\\
    =&\ \frac{(\frac{16\pi^2}{405}-\frac{34}{45})c_0^{-1/2}}{\frac{74}{135}+\frac{4\pi^2}{405}}\bigl(e^{(\frac{74}{135}+\frac{4\pi^2}{405})\sigma^2 t}-1\bigr).\label{eqn:ximeanscalar}
\end{align}
For the variance, we write
\begin{align*}
    \operatorname{Var}[\Omega_0(t)]=&\ \tfrac{4}{9}c_*^{-1}\sigma^2 \operatorname{Var}\bigl[ \int_0^t \alpha_0(s) \ \d \beta_s\bigr]\\
    &+\tfrac{4}{3}c_*^{-1/2}\overline{\mu}_{d;I}(0)\sigma^3 \operatorname{Cov}\bigl(\int_0^t \alpha_0(s) \ \d s,\int_0^t \alpha_0(s)\  \d \beta_s\bigr)\\
    &+\overline{\mu}_{d;I}(0)^2 \sigma^4\operatorname{Var}\bigl[\int_0^t \alpha_0(s)\ \d s\bigr],
\end{align*}
and explicitly compute the leading-order term
\begin{align}
    \sigma^2 \operatorname{Var}\bigl[ \int_0^t \alpha_0(s)  \ \d \beta_s\bigr]=&\ \tfrac{4}{9}c_*^{-1}\sigma^2\mathbb{E}\left(\int_0^t \alpha_0(s)\ \d \beta_s\right)^2\nonumber\\
    =&\ \tfrac{4}{9}c_*^{-1}\sigma^2\int_0^t \mathbb{E} \bigl[\alpha_0^2(s)\bigr]\ \d s\nonumber\\
    =&\ \tfrac{45}{2c_*(78+\pi^2)} \bigl(e^{(\frac{208}{135}+\frac{8\pi^2}{405})\sigma^2 t}-1\bigr)\label{eqn:xivariancescalar}.
\end{align}
In Figure~\ref{fig:xivariancescalar} we compare these statistics of $\Omega_0(t)$ to sample statistics of $\Omega(t)$. The sample variance agrees well with the prediction \eqref{eqn:xivariancescalar}. The sample mean, however, differs slightly from the prediction \eqref{eqn:ximeanscalar}. We observe that the approximation is significantly improved upon considering the sample mean of $\Omega_2(t)$.

\begin{figure*}[h!]
    \centering
    \begin{subfigure}[t]{0.5\textwidth}
        \centering
        \includegraphics[height=2in]{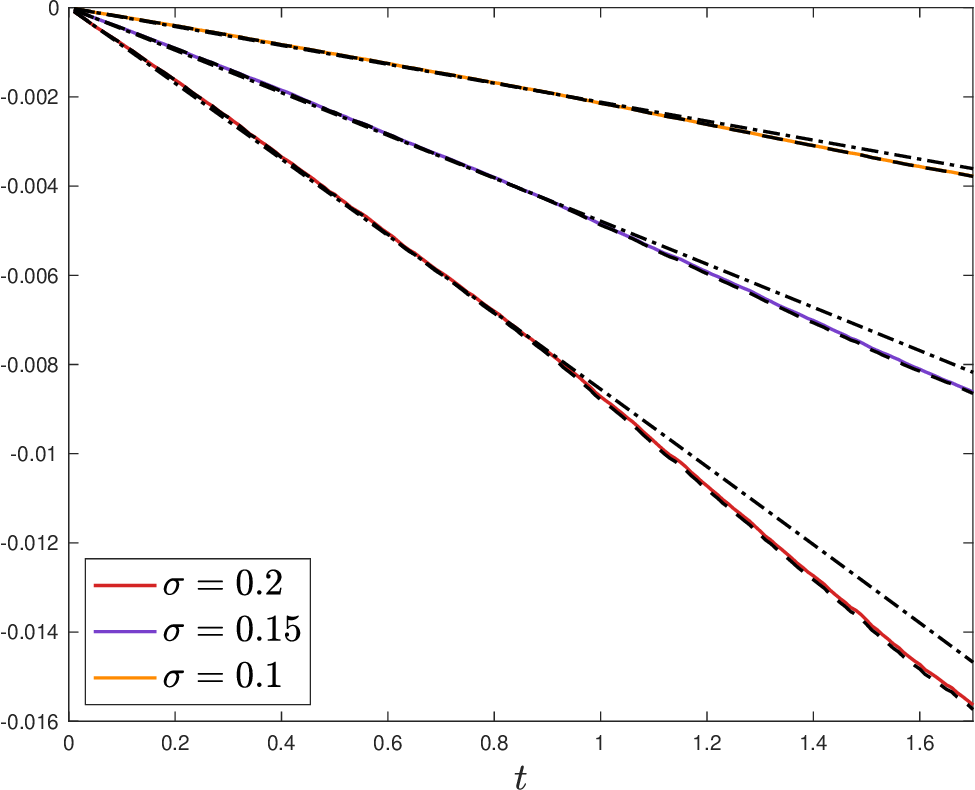}
        \caption{Sample mean of $ \Omega(t)$. }
    \end{subfigure}%
    ~
\begin{subfigure}[t]{0.5\textwidth}
        \centering
        \includegraphics[height=2in]{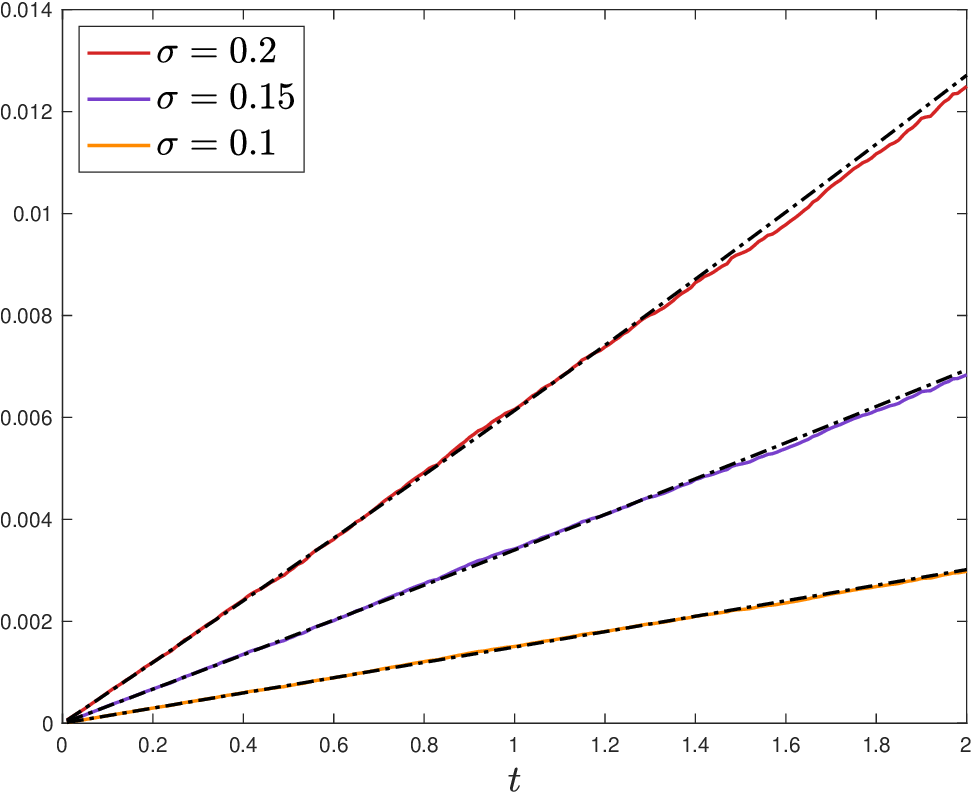}
        \caption{Sample variance of $ \Omega(t)$.}
    \end{subfigure}
    \caption{Sample statistics of the process $\Omega(t)$ for scalar noise. Computed over $3\cdot10^3$ realizations for $\sigma\in\{0.1,0.15,0.2\}.$ Dash-dotted lines indicate the theoretical mean and (leading-order) variance of $\Omega_0(t)$ as in \eqref{eqn:ximeanscalar} and \eqref{eqn:xivariancescalar}, respectively. Panel (a) also shows the sample mean of $\Omega_2(t)$ (dotted), which gives significant improvement over the mean of $\Omega_0(t)$. }\label{fig:xivariancescalar}
\end{figure*}


\paragraph{Remainders}
In order to confirm that our approximation procedure only neglects higher-order noise effects, we numerically investigate the resulting error. Figure~\ref{fig:errors} shows the growth of the remainders $\|v(t)- v_1(t)\|_{L_a^2}$ and $|c(t)-c_2(t)|$. Indeed, the size of the remainders $\|v(t)- v_1(t)\|_{L_a^2}$ and $|c(t)-c_2(t)|$ decreases significantly with decreasing values of $\sigma$. An estimation of the order $\beta$ at which these remainders depend on $\sigma$ (see Figure~\ref{fig:errororder} in \ref{app:supplementary}) reveals that the remainders $\|v(t)- v_1(t)\|_{L_a^2}$ and $|c-c_2|$ scale with a power of $\sigma$ higher than 2 and 3, respectively. This indicates that $c_{2}$ indeed captures all effects of $O(\sigma^2)$.

\begin{figure*}[h!]
    \centering
    \begin{subfigure}[t]{0.5\textwidth}
        \centering
        \includegraphics[height=2in]{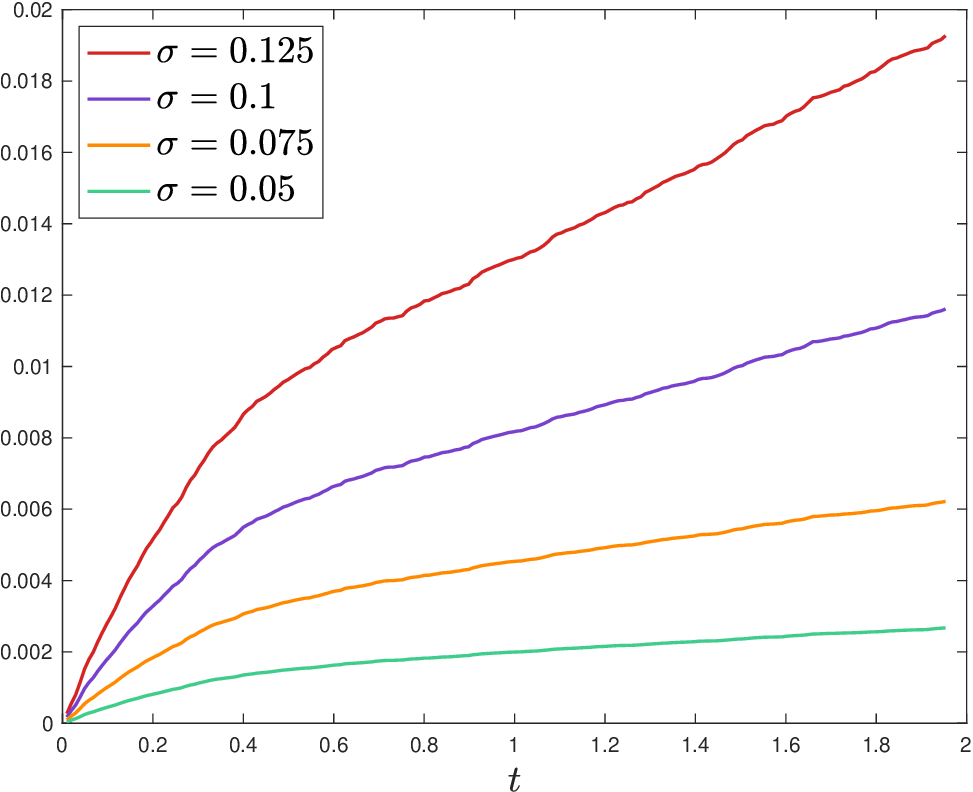}
        \caption{Sample mean of $\sup_{s\leq t}\|v(s)- v_1(s)\|_{L_a^2}$. }
    \end{subfigure}%
    ~
    \begin{subfigure}[t]{0.5\textwidth}
        \centering
        \includegraphics[height=2in]{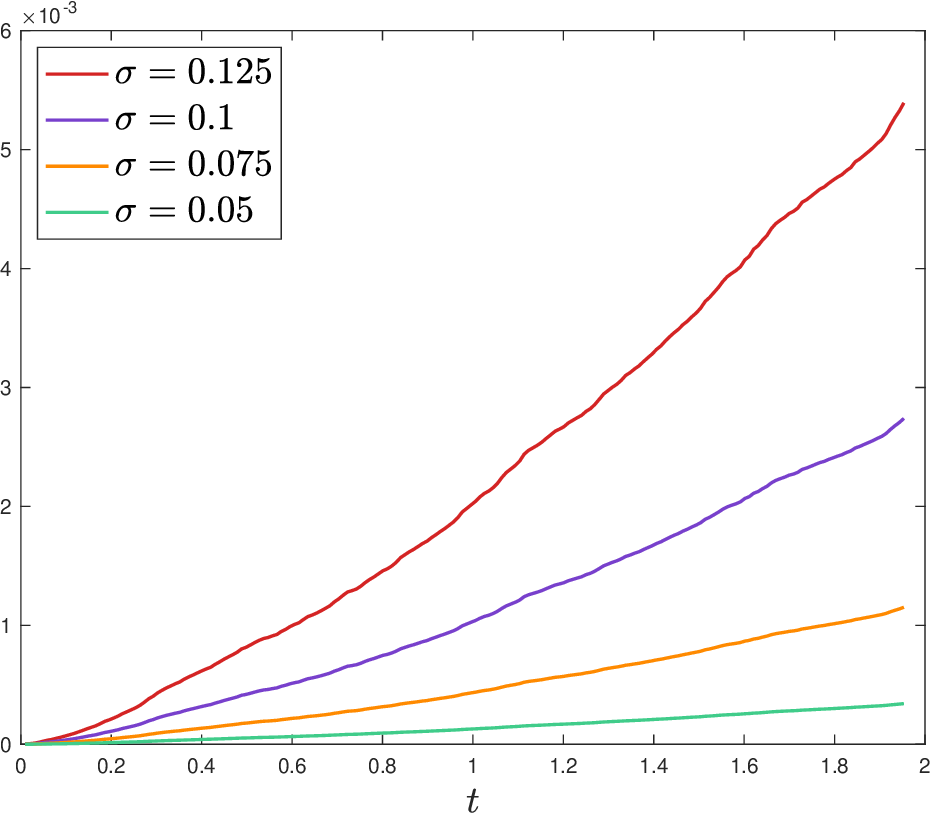}
        \caption{Sample mean of $\sup_{s\leq t}|c(s)-c_2(s)|$. }
    \end{subfigure}%
    \caption{The error made by approximating $v$ with the first order expansion $v_1$ and $c$ with the second order expansion $ c_{2}$, for scalar noise. Computed over 200 simulations, for $\sigma\in\{0.05,0.075,0.1,0.125\}$.}\label{fig:errors}
\end{figure*}

\subsection{Example I\joinR I\joinR I: Soliton dynamics for space-time white noise}
\label{subsec:whitenoiseapprox}
We now turn to the setting of space-time white noise, introduced in \S\ref{subsubsec:spacetimewhite}, where $\alpha$ and $\xi$ follow the modulation equations \eqref{eqn:modulationalphaspacetime}-\eqref{eqn:modulationxispacetime}.  A key difference in this example is that the noise is space-dependent. The modulation equations are formulated in the stochastic co-moving frame, where the noise undergoes a spatial rescaling by $T_{\alpha}$.  We discuss this noise transformation in \ref{app:noiserescaling}, where we show that the process $\tilde{W}_t:=\alpha^{1/2}T_{\alpha}W_t$ generates the same statistics as the white noise $W_t$. In what follows, we formulate the modulation equations \eqref{eqn:modulationalphaspacetime}-\eqref{eqn:modulationxispacetime} using the space-time white noise $\tilde{W}_t$.

The first approximation to the rescaling process $\alpha$ is then defined as the process $\alpha_0$  which satisfies the SDE
\begin{align}\label{eqn:alphafirstorderspacetime}
    \d \alpha_0=\overline{\gamma}_{d;\RomanIII}(0)\sigma^2 \ \d t-\sigma \alpha_0^{1/2}\bigl\langle \d \tilde{W}_t,\overline{\gamma}_{\diamond}(0)\bigr\rangle_{L^2}, 
\end{align}
with $\alpha_0(0)=1$. Note that, in distribution, the Brownian motion driving this SDE equals
\[-\sigma \bigl\langle \tilde{W}_t,\overline{\gamma}_{\diamond}(0)\bigr\rangle_{L^2}\stackrel{d}{=}\bigl\|\overline{\gamma}_{\diamond}(0)\bigr\|_{L^2}\sigma \beta_t=\sqrt{\tfrac{4}{35}}c_*^{1/4}\sigma \beta_t,\] 
so that \eqref{eqn:alphafirstorderspacetime} is of the form
\begin{align}\d X(t)=\delta \ \d t+s  X^{1/2}(t)\ \d \beta_t.\label{eqn:Bessel}
\end{align}
The solution to \eqref{eqn:Bessel} with $\delta\geq0$ is known as a squared Bessel process \cite{bessel}. The squared Bessel process $X(t)$ remains strictly positive for $\delta\geq \tfrac{1}{2}s^2$, and therefore the drift component $\overline{\gamma}_{d;\RomanIII}(0)\approx  0.093c_*^{1/2} \geq \tfrac{2}{35}c_*^{1/2}$ in \eqref{eqn:alphafirstorderspacetime} is large enough to ensure that $\alpha_0(t)$ remains strictly positive. The mean of the approximation $\alpha_0$ is easily computed as\footnote{Higher order and negative moments can in principle be computed explicitly by using that a Bessel process $X(t)$ as defined by \eqref{eqn:Bessel} has the noncentral Chi-square distribution $\frac{4}{s^2 t}X(t)\sim \chi^2_{4\delta/s^2}\bigl(\frac{4}{s^2 t}X(0)\bigr)$.}
\[\mathbb{E}\bigl[\alpha_0(t)\bigr]=1+\overline{\gamma}_{d;\RomanIII}(0)\sigma^2 t\approx 1+  0.093c_*^{1/2}\sigma^2 t. \]

At first order, the perturbation is given by 
\begin{align*}
    v_1(t) =&\ \sigma \int_0^t \alpha^{-1/2}(s)e^{\int_s^t \alpha_0^{-3}(t^\prime)\d t^\prime \mathcal{L}_{c_*}} S_{\diamond}(0)[\d W_s] ,
\end{align*}
where
\begin{align}
S_{\diamond}(0)[h]=&\ \phi_{c_*} h - \tfrac{1}{9}c_*^{-3/2}(x \partial_x+2)\phi_{c_*}\langle h,\phi_{c_*}^2\rangle_{L^2}\nonumber\\
   &-  \tfrac{2}{9}\partial_x\phi_{c_*} \langle h,c_*^{-2}\phi_{c_*}^2-c_*^{-1/2}\phi_{c_*}\zeta_{c_*}\rangle_{L^2}.
   \end{align}
The second order approximation for the perturbation is given explicitly as
\[v_2(t)=\sigma V^{(1)}(\alpha_1,t)+\sigma^2 V^{(2)}(\alpha_1,t),\]
using 
\begin{align*}
    V^{(1)}(\alpha,t) =&  \int_0^t \alpha^{-1/2}(s)e^{\int_s^t \alpha^{-3}(t^\prime)\d t^\prime \mathcal{L}_{c_*}} S_{\diamond}(0)[\d W_s]
\end{align*}
and
\begin{align*}
    V^{(2)}(\alpha,t)=&\int_0^t  \alpha^{-3}(s) e^{\int_s^t \alpha^{-3}(t^\prime)\d t^\prime \mathcal{L}_{c_*}}  \Bigl[N\bigl( V^{(1)}(\alpha,s)\bigr)+[R_0]^{(2)}\bigl( V^{(1)}(\alpha,s)\bigr)\Bigr]\ \d s\\
&+\int_0^t \alpha^{-1}(s)e^{\int_s^t \alpha^{-3}(t^\prime)\d t^\prime \mathcal{L}_{c_*}}\sum_{i=1}^6 R_{i;\RomanIII}(0)\ \d s\\
&+ \int_0^t \alpha^{-1/2}(s)e^{\int_s^t \alpha^{-3}(t^\prime)\d t^\prime \mathcal{L}_{c_*}} [S_{\diamond}]^{(1)} \bigl(V^{(1)}(\alpha,s)\bigr)[\d W_s].
\end{align*}
The subsequent approximation $\alpha_1$ to the rescaling process satisfies the square root SDE
\[\d \alpha_{1} =\sigma^2K^{1,1}_{\RomanIII}(t)  \ \d t- \sigma  \alpha_1^{1/2} \bigl\langle \d \tilde{W}_t,K^{2,1}_{\RomanIII}(t)\bigr\rangle_{L^2},\]
with random coefficient
\begin{align*}
    K^{1,1}_{\RomanIII}(t)=&\ \overline{\gamma}_{d;\RomanIII}(0)+[\overline{\gamma}_{d;\RomanIII}]^{(1)}\bigl(v_1(t)\bigr), 
\end{align*}
where $K^{2,1}_{\RomanIII}$ is the process
\begin{align*}
    K^{2,1}_{\RomanIII}(t)=&\ \overline{\gamma}_{\diamond}(0)+[\overline{\gamma}_{\diamond}]^{(1)}\bigl(v_1(t)\bigr).
\end{align*}
Subsequently, the approximation $\alpha_2$ is given by the SDE
\begin{align*}\d \alpha_2=&\ [- K^{0,2}_{\RomanIII}(t)\alpha_2^{-2}+ \sigma^2 K^{1,2}_{\RomanIII}(t)] \ \d t-\sigma   \alpha_2^{1/2} \bigl\langle \d \tilde{W}_t,K^{2,2}_{\RomanIII}(t)\bigr\rangle_{L^2}
\end{align*}
with random coefficients 
\begin{align*}
    K^{0,2}_{\RomanIII}(t)=&\ [\overline{\gamma}_d^0]^{(2)}\bigl(\sigma V^{(1)}(\alpha_1,t)\bigr),\\
    K^{1,2}_{\RomanIII}(t)=&\ \overline{\gamma}_{d;\RomanIII}(0)+[\overline{\gamma}_{d;\RomanIII}]^{(1)}\bigl(v_2(t)\bigr) +[\overline{\gamma}_{d;\RomanIII}]^{(2)}\bigl(\sigma V^{(1)}(\alpha_1,t)\bigr),
\end{align*}
where $K^{2,2}_{\RomanIII}(t)$ is the process
\begin{align*}
    K^{2,2}_{\RomanIII}(t)=&\ \overline{\gamma}_{\diamond}(0)+[\overline{\gamma}_{\diamond}]^{(1)}\bigl(v_2(t)\bigr)+[\overline{\gamma}_{\diamond}]^{(2)}\bigl(\sigma V^{(1)}(\alpha_1,t)\bigr).
\end{align*}
\paragraph{Amplitude}
Using It\^{o}'s lemma and \eqref{eqn:alphafirstorderspacetime}, we find that the first approximation $c_0(t)=c_* \alpha_0^{-2}(t) $ for the amplitude process $c(t)$ has the It\^o form
\begin{align*}
c_0(t) =&\ c_*+2c_*\sigma \int_0^t \alpha_0^{-5/2} \bigl\langle \d  \tilde{W}_s,\overline{\gamma}_{\diamond}(0) \bigr\rangle_{L^2}\\
 &+c_*\sigma^2\bigl[-2 \overline{\gamma}_{d;\RomanIII}(0)+3  \| \overline{\gamma}_{\diamond}(0)\|^2_{L^2}\bigr]\int_0^t \alpha_0^{-3}\ \d s.
\end{align*}
We compute the leading-order variance
\begin{align}\label{eqn:whitenoisevariance}
    \operatorname{Var}[ c_0(t)]=&\ 4  \sigma^2 c_*^2 \mathbb{E}\biggl[\Bigl(\int_0^t \alpha_0^{-5/2}\bigl\langle \d W_s,\overline{\gamma}_{\diamond}(0) \bigr\rangle_{L^2} \Big)^2\biggr] +O(\sigma^3)\nonumber\\
    =&\ 4\sigma^2 c_*^2 \bigl\|\overline{\gamma}_{\diamond}(0)\bigr\|^2_{L^2}\int_0^t \mathbb{E}\bigl[\alpha_0^{-5}\bigr] \ \d s+O(\sigma^3)\nonumber\\
    =&\ \tfrac{16}{35}\sigma^2 c_*^{5/2}\int_0^t \mathbb{E}\bigl[\alpha_0^{-5}\bigr] \ \d s+O(\sigma^3)
\end{align}
and compare this expression to the sample variance of $c_{\text{fit}}(t)$ in Figure~\ref{fig:whitenoiseamplitudevariance}.

Figure~\ref{fig:whitenoiseamplitudemean} compares the mean of $c_{\text{fit}}(t)$ with that of the increasingly refined approximations $c_0(t), c_1(t)$ and $c_2(t)$. Here we use 
\begin{align}
    \mathbb{E}[c_0(t)] =& \, c_*+c_*\sigma^2\Bigl[-2 \overline{\gamma}_{d;\RomanIII}(0)+3  \bigl\| \overline{\gamma}_{\diamond}(0)\bigr\|^2_{L^2}\Bigr]\int_0^t \mathbb{E}\bigl[\alpha_0^{-3}(s)\bigr]\ \d s\\
    \approx & \, c_*+0.16 c_*^{3/2}\sigma^2\int_0^t \mathbb{E}\bigl[\alpha_0^{-3}(s)\bigr]\ \d s,
\end{align}
and numerically compute the negative moment $\mathbb{E}[\alpha_0^{-3}(t)]$\footnote{We remark that, alternatively, this moment can be evaluated numerically by repeated integration of the moment generating function of $\alpha_0(t)$, which has a noncentral Chi-square distribution.}. For the means of $c_1(t)$ and $c_2(t)$ we, once more, resort to the sample mean. As is the case for scalar noise, the quadratic terms of $v$ contribute significantly to the mean amplitude, and the mean of $c_{\text{fit}}(t)$ is not well-approximated by that of $c_0(t)$. We see in Figure~\ref{fig:whitenoiseamplitudemean} that the amplitude drift is well-approximated by the mean of $c_2(t)$.


\begin{figure*}[h!]
        \centering
 \begin{subfigure}[t]{0.5\textwidth}
        \centering
        \includegraphics[height=2in]{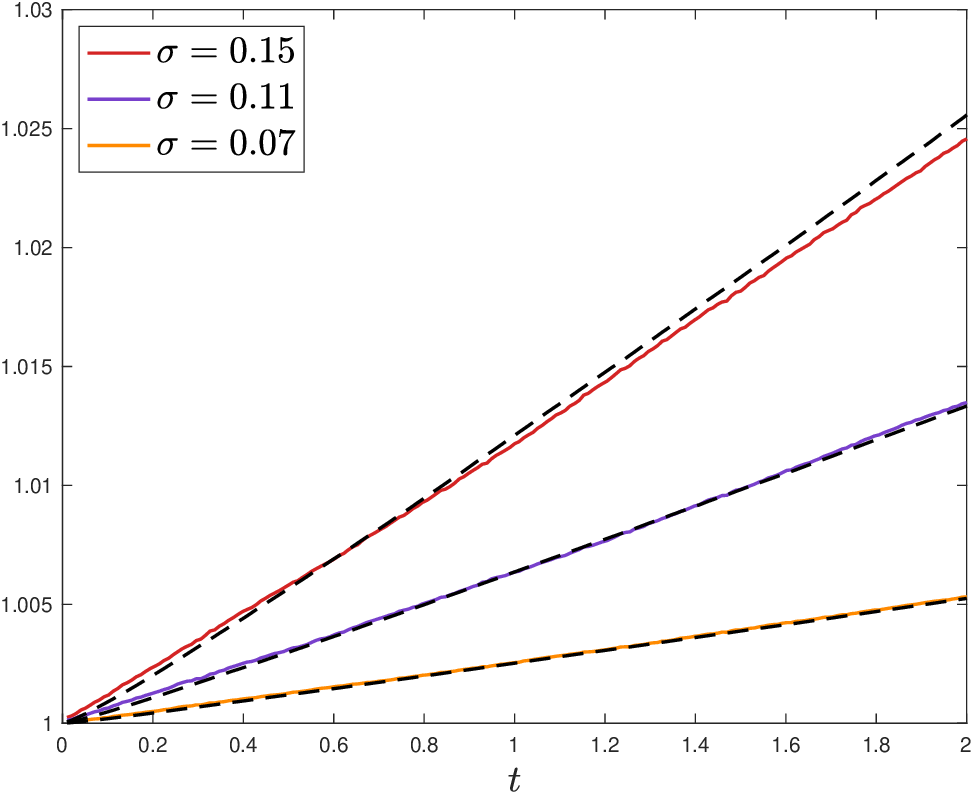}
        \caption{Mean of $c_{\text{fit}}(t)/c_*$ at various noise strengths.}
    \end{subfigure}%
    ~
    \begin{subfigure}[t]{0.5\textwidth}
        \centering
        \includegraphics[height=2in]{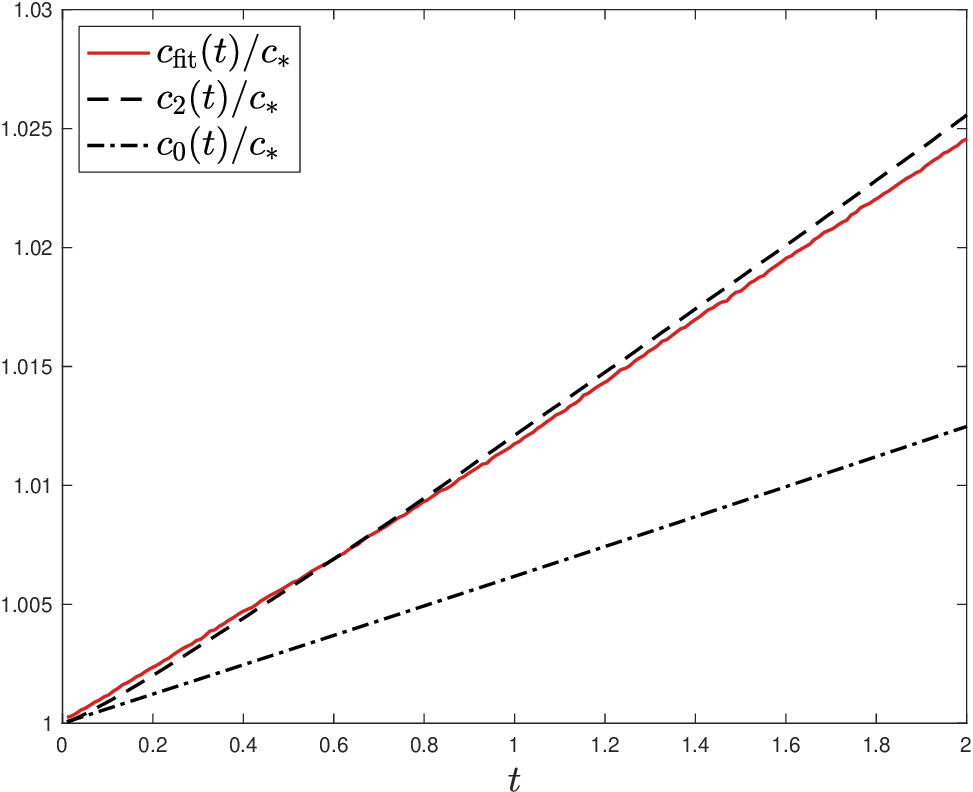}
        \caption{Means of $c_{\text{fit}}(t)/c_*$ and approximations $c_0(t)$ and $c_2(t)$.}
    \end{subfigure}%

    \caption{Sample mean of the process $c_{\text{fit}}(t)/c_*$ computed over $36\cdot10^3$ realisations for space-time white noise. The mean of $c_2(t)$ (dashed) agrees well with that of $c_{\text{fit}}(t)$ at the simulated noise strength values $\sigma\in\{0.07,0.11,0.15\}$, whereas the mean of $c_0(t)$ (dash dot) as in \eqref{eqn:meanc0} fails to capture the correct amplitude drift. The mean of $c_1(t)$ provides no improvement, it can not be distinguished from that of $c_0(t)$ at these simulation values.}\label{fig:whitenoiseamplitudemean}
\end{figure*}

\paragraph{Phase shift}
The first approximation $\Omega_0(t)$ to the phase shift process $\Omega(t)$ defined in \eqref{eqn:phaseshift} is given by
\begin{align}
\Omega_0(t)=&\ \sigma^2 \overline{\mu}_{d;\RomanIII}(0) t -\sigma \int_0^t\alpha_0^{1/2} \bigl\langle \d \tilde{W}_s,\overline{\mu}_{\diamond}(0)\bigr\rangle_{L^2}.\label{eqn:omega0}
\end{align}
We can explicitly compute the variance of the phase shift approximation $\Omega_0(t)$:
\begin{align}\label{eqn:xivariancewhite}
    \operatorname{Var}[ \Omega_0(t)]
    =&\   \sigma^2 \mathbb{E}\left( \int_0^t \alpha_0^{1/2}\langle \d \tilde{W}_s,\overline{\mu}_{\diamond}(0)\rangle_{L^2} \right)^2\nonumber\\
    =&\ \sigma^2\bigl\|\langle \cdot,\overline{\mu}_{\diamond}(0)\rangle_{L^2}\bigr\|^2_{HS(L^2,\R)}\int_0^t \mathbb{E}\bigl[\alpha_0 (s)\bigr]\ \d s\nonumber\\
    =&\ \sigma^2\bigl\|\overline{\mu}_{\diamond}(0)\bigr\|_{L^2}^2 (t+\tfrac{1}{2}\overline{\gamma}_{d;\RomanIII}(0)\sigma^2 t^2) \nonumber\\
    \approx &\ 0.435 c_*^{-1/2} \sigma^2t+0.02\sigma^4 t^2.
\end{align}
This expression agrees well with the sample variance of $\Omega(t)$, as shown in Figure~\ref{fig:xivariancewhite}.

\begin{figure*}[h!]
\centering
 \begin{subfigure}[t]{0.33\textwidth}
        \centering
        \includegraphics[height=1.45in]{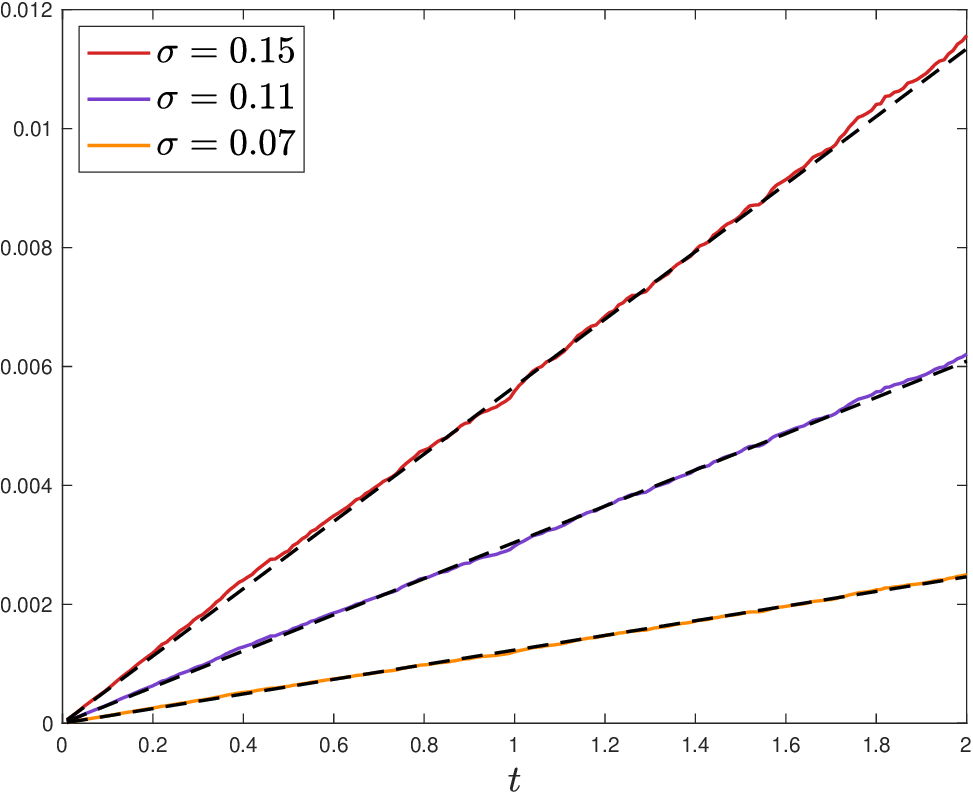}
    \caption{Sample variance of $\Omega(t)$.}\label{fig:xivariancewhite}
    \end{subfigure}%
    ~
    \begin{subfigure}[t]{0.33\textwidth}
        \centering
        \includegraphics[height=1.5in]{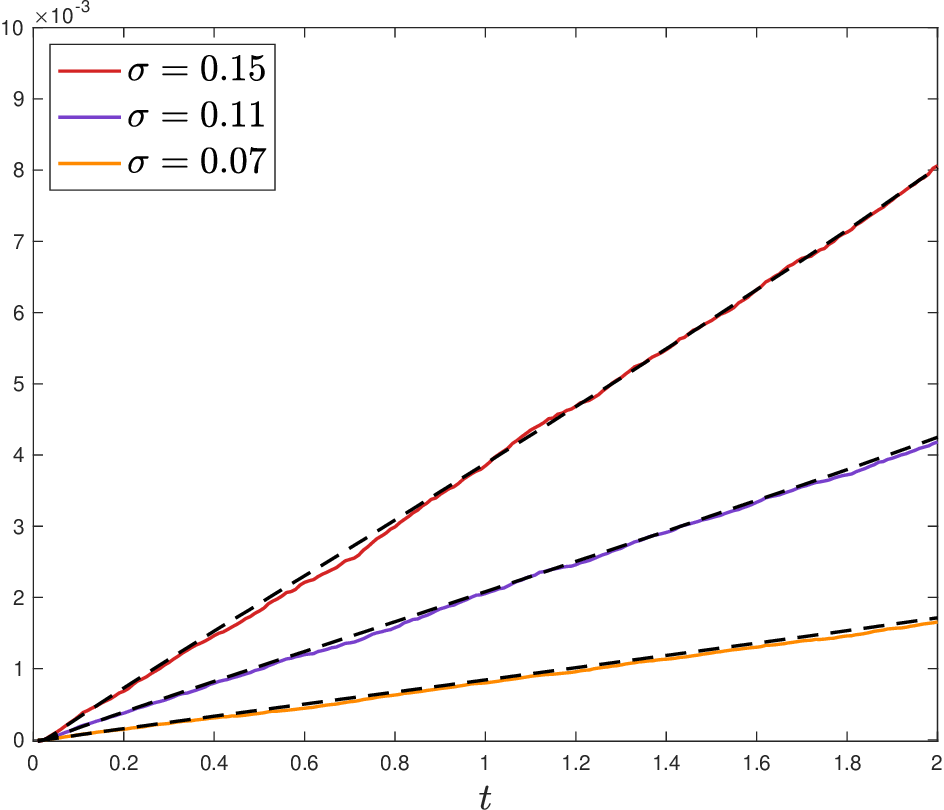}
    \caption{Sample mean of $\Omega(t)$.}\label{fig:ximeanwhite}
    \end{subfigure}%
    ~
    \begin{subfigure}[t]{0.33\textwidth}
        \centering
        \includegraphics[height=1.5in]{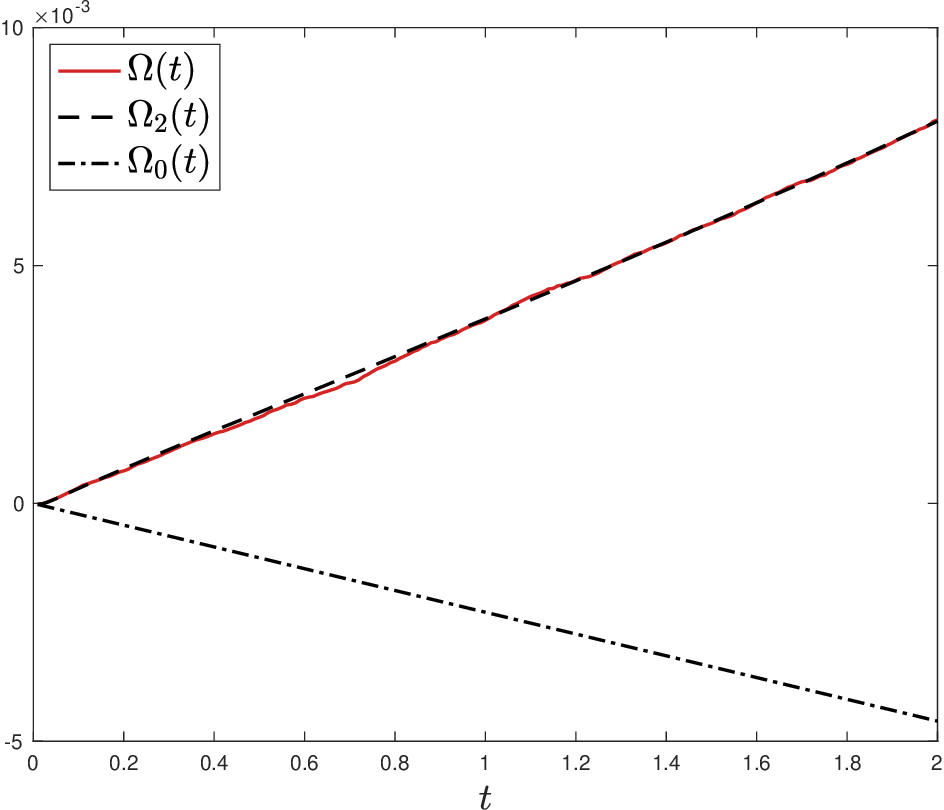}
    \caption{Means of $\Omega_{\text{fit}}(t)$ and approximations $\Omega_0(t)$ and $\Omega_2(t)$.}\label{fig:ximeanwhite0}
    \end{subfigure}%
    \caption{Sample statistics (solid) of the process $\Omega(t)$ for space-time white noise, at noise strengths $\sigma\in\{0.07,0.11,0.15\}$. Dashed lines indicate the theoretical variance of $\Omega_0(t)$ as in \eqref{eqn:xivariancewhite} and the sample mean of $\Omega_2(t)$, respectively. The dash-dotted line in panel (c) shows the mean of $\Omega_0(t)$, see \eqref{eqn:omega0}, which fails to capture the correct phase drift. The sample variance is computed over $2500$ realizations and the sample mean over $10^4$ realizations.}
\end{figure*}

The mean phase shift from the primary velocity $c(t)$ is influenced by quadratic terms of the perturbation $v$. We thus consider the process 
\begin{align*}
    \Omega_2(t) =&-\int_0^t M^{0,2}_{\RomanIII}(s)\alpha_2^{-2}(s)\ \d s+\sigma^2 \int_0^t M^{1,2}_{\RomanIII}(s) \ \d s \\
    &-\sigma \int_0^t\alpha_2^{1/2} \bigl\langle  \d \tilde{W}_s,M^{2,2}_{\RomanIII}(s)\bigr\rangle_{L^2}
\end{align*}
with random coefficients 
\begin{align*}
    M^{0,2}_{\RomanIII}(t)=&\ [\overline{\mu}_d^0]^{(2)}\bigl(V_1(\alpha_2,t)\bigr),\\ 
    M^{1,2}_{\RomanIII}(t)=&\ \overline{\mu}_{d;\RomanIII}(0)+[\overline{\mu}_{d;\RomanIII}]^{(1)}\bigl(V_2(\alpha_2,t)\bigr) +[\overline{\mu}_{d;\RomanIII}]^{(2)}\bigl(V_1(\alpha_2,t)\bigr)
    \end{align*}
and where $M^{2,2}_{\RomanIII}$ is the process
\begin{align*}
    M^{2,2}_{\RomanIII}(t)=&\ \overline{\mu}_{\diamond}(0)+[\overline{\mu}_{\diamond}]^{(1)}\bigl(V_2(\alpha_2,t)\bigr)+[\overline{\mu}_{\diamond}]^{(2)}\bigl(V_1(\alpha_2,t)\bigr).
\end{align*}
We compare a numerical evaluation of the mean
\begin{align*}
    \mathbb{E}[\Omega_2(t)] =&-\int_0^t \mathbb{E}\bigl[ M^{0,2}_{\RomanIII}(s)\alpha_2^{-2}(s)\bigr]\ \d s+\sigma^2 \int_0^t  \mathbb{E}\bigl[M^{1,2}_{\RomanIII}(s) \bigr]\ \d s
\end{align*}
to the sample mean of $\Omega(t)$ in Figure~\ref{fig:ximeanwhite} and observe that they agree quite well.

\appendix

\renewcommand{\thesection}{\appendixname\Alph{section}}
\section{Frozen-frame transformation}
\label{app:frozen}

Our goal here is to derive the SPDE \eqref{eqn:perturbationgeneral} for the remainder $v=\alpha^2 T_{\alpha,\xi}u-\phi_{c_*}$, where we recall that the processes $\alpha$ and $\xi$ satisfy \eqref{eqn:alphaform} and \eqref{eqn:xiform}, respectively.
We define  $\Phi_{\alpha,\xi} u=\alpha^2 T_{\alpha,\xi} u$ so that $v=\Phi_{\alpha,\xi}u-\phi_{c_*}$. For an arbitrary test-function $\zeta$, we now characterize the evolution of the real-valued process
\[\langle v, \zeta \rangle_{L^2}=\langle \Phi_{\alpha,\xi}u-\phi_{c_*},\zeta \rangle_{L^2}.\]
To this end, we collect the first and second-order (Fr\'{e}chet) derivatives of the mapping $(u,\alpha,\xi)\mapsto \langle \Phi_{\alpha,\xi} u,\zeta\rangle_{L^2} $. For the derivatives with respect to $u,\alpha$ and $\xi$ we can write
\begin{align*}
\partial_{u}\langle \Phi_{\alpha,\xi} u,\zeta\rangle_{L^2} [v]=&\ \alpha^2 \langle T_{\alpha,\xi}v,\zeta\rangle_{L^2} \\
\partial_\alpha \langle \Phi_{\alpha,\xi}u,\zeta \rangle_{L^2} 
=&\ 2\alpha \langle T_{\alpha,\xi}u,\zeta\rangle_{L^2}+\alpha^2 \langle x T_{\alpha,\xi}u_x, \zeta\rangle_{L^2} \\
\partial_\xi \langle \Phi_{\alpha,\xi}u,\zeta \rangle_{L^2}=&\ \alpha^2 \langle T_{\alpha,\xi}u_x,\zeta \rangle_{L^2},
\end{align*}
and for the second derivatives we find
\begin{align*}
\partial^2_{u}\langle \Phi_{\alpha,\xi}u,\zeta\rangle_{L^2}[v,w]=&\ 0,\\
\partial^2_{\alpha } \langle \Phi_{\alpha,\xi}u,\zeta\rangle_{L^2} =&\ 2\langle T_{\alpha,\xi}u,\zeta\rangle_{L^2}+4\alpha \langle xT_{\alpha,\xi}u_x,\zeta\rangle_{L^2}+\alpha^2 \langle x^2T_{\alpha,\xi}u_{xx},\zeta\rangle_{L^2},\\
\partial^2_{\xi} \langle \Phi_{\alpha,\xi}u,\zeta\rangle_{L^2}=&\ \alpha^2 \langle T_{\alpha,\xi} u_{xx},\zeta\rangle_{L^2},\\
\partial_{\alpha \xi} \langle \Phi_{\alpha,\xi} u,\zeta\rangle_{L^2}=&\  2\alpha \langle T_{\alpha,\xi}u_x,\zeta\rangle_{L^2}+\alpha^2 \langle xT_{\alpha,\xi}u_{xx},\zeta \rangle_{L^2},\\
\partial_{u \alpha} \langle \Phi_{\alpha,\xi}u,\zeta\rangle_{L^2}[v]=&\ 2\alpha \langle T_{\alpha,\xi} v,\zeta \rangle_{L^2}+\alpha^2 \langle x T_{\alpha,\xi}v_x,\zeta \rangle_{L^2},\\
\partial_{u \xi} \langle \Phi_{\alpha,\xi}u,\zeta \rangle_{L^2}[v]=&\ \alpha^2 \langle T_{\alpha,\xi}v_x,\zeta \rangle_{L^2}.
\end{align*}
Upon choosing an orthonormal basis $\{e_k\}_{k=0}^\infty$ of $\mathcal{H}$, It\^o's lemma \cite[Theorem 4.32]{daprato} now gives
\begin{align*}
\d \langle v ,\zeta \rangle_{L^2} =& \sum_{i=0}^6 \overline{R}^\sigma_i(u,\alpha,\xi,\zeta)\ \d t+ \overline{S}^\sigma(u,\alpha,\xi,\zeta)[\d W_t^Q],
\end{align*}
where
\begin{align*} \overline{R}^\sigma_0(u,\alpha,\xi,\zeta)&=\partial_{u}\langle \Phi_{\alpha,\xi} u,\zeta\rangle_{L^2} [-u_{xxx}-2u u_x],\\
\overline{R}^\sigma_1(u,\alpha,\xi,\zeta)&= \tfrac{1}{2}\sum_{k=0}^\infty\gamma^{\sigma}_s[Q^{1/2} e_k]^2 \partial^2_{ \alpha} \langle \Phi_{\alpha,\xi}u,\zeta\rangle_{L^2},\\
\overline{R}^\sigma_2(u,\alpha,\xi,\zeta)&=\tfrac{1}{2} \sum_{k=0}^\infty \mu^{\sigma}_s [Q^{1/2} e_k]^2  \partial^2_{ \xi} \langle \Phi_{\alpha,\xi}u,\zeta\rangle_{L^2}, \\
\overline{R}^\sigma_3(u,\alpha,\xi,\zeta)&=\sum_{k=0}^\infty\gamma^{\sigma}_s [Q^{1/2} e_k]\mu^{\sigma}_s [Q^{1/2} e_k] \partial_{\alpha \xi} \langle \Phi_{\alpha,\xi} u,\zeta\rangle_{L^2},\\
\overline{R}^\sigma_4(u,\alpha,\xi,\zeta)&= \sigma \sum_{k=0}^\infty\gamma^{\sigma}_s [Q^{1/2}e_k] \partial_{u \alpha} \langle \Phi_{\alpha,\xi}u,\zeta\rangle_{L^2}\bigl[  M(u) [Q^{1/2}e_k]\bigr],\\
\overline{R}^\sigma_5(u,\alpha,\xi,\zeta)&= \sigma \sum_{k=0}^\infty\mu^{\sigma}_s [Q^{1/2}e_k] \partial_{u \xi} \langle \Phi_{\alpha,\xi}u,\zeta \rangle_{L^2}\bigl[ M(u) [Q^{1/2}e_k]\bigr],\\
\overline{R}^\sigma_6(u,\alpha,\xi,\zeta)&=\gamma^\sigma_d \partial_\alpha \langle \Phi_{\alpha,\xi}u,\zeta \rangle_{L^2}  +\mu^\sigma_d \partial_\xi \langle \Phi_{\alpha,\xi}u,\zeta \rangle_{L^2},
\end{align*}
and
\begin{align*}
\overline{S}^\sigma(u,\alpha,\xi,\zeta)[h]=&\ \sigma \partial_{u}\langle \Phi_{\alpha,\xi} u,\zeta\rangle_{L^2} \bigl[ M(u) [h] \bigr]+ \partial_\alpha \langle \Phi_{\alpha,\xi}u,\zeta \rangle_{L^2} \gamma^{\sigma}_s[h]\\&+ \partial_\xi \langle \Phi_{\alpha,\xi}u,\zeta \rangle_{L^2}\mu^{\sigma}_s[h]. \end{align*}
In the term $\overline{R}^\sigma_3(u,\alpha,\xi,\zeta)$, we simplify the summation as
\begin{align*}
   \sum_{k=0}^\infty \gamma^{\sigma}_s [Q^{1/2} e_k]\mu^{\sigma}_s [Q^{1/2} e_k]=&\ \sigma^2\alpha^2\sum_{k=0}^\infty \langle \hat{T}_{\alpha,\xi}Q^{1/2}e_k,\overline{\gamma}_s \rangle_{\mathcal{H}}\langle \hat{T}_{\alpha,\xi} Q^{1/2}e_k,\overline{\mu}_s \rangle_{\mathcal{H}}\\
   =&\ \sigma^2\alpha^2\langle Q^{1/2}\hat{T}^*_{\alpha,\xi}\overline{\gamma}_s,Q^{1/2}\hat{T}^*_{\alpha,\xi}\overline{\mu}_s \rangle_{\mathcal{H}},
\end{align*}
where we have used that $\gamma_s^\sigma$ is of the form \eqref{eqn:gammasrepresentation}. In the mixed derivative
\begin{align*}
    \partial_{\alpha \xi} \langle \Phi_{\alpha,\xi} u,\zeta\rangle_{L^2}=2\alpha \langle T_{\alpha,\xi}u_x,\zeta\rangle_{L^2}+\alpha^2 \langle xT_{\alpha,\xi}u_{xx},\zeta \rangle_{L^2}
\end{align*}
we substitute 
\begin{align}\label{eqn:subst}
   T_{\alpha,\xi}[\partial_x^j u]= \alpha^{-(j+2)}  \partial_x^j [\phi_{c_*}+v] ,
\end{align}
for $j=1,2$ and find
\begin{align*}
\overline{R}^\sigma_3(u,\alpha,\xi,\zeta)=&\ 2\sigma^2\langle Q^{1/2}\hat{T}^*_{\alpha,\xi}\overline{\gamma}_s,Q^{1/2}\hat{T}^*_{\alpha,\xi}\overline{\mu}_s \rangle_{\mathcal{H}}
\langle \partial_x[\phi_{c_*}+v],\zeta\rangle_{L^2}\\
&+\sigma^2\langle Q^{1/2}\hat{T}^*_{\alpha,\xi}\overline{\gamma}_s,Q^{1/2}\hat{T}^*_{\alpha,\xi}\overline{\mu}_s \rangle_{\mathcal{H}}\langle x \partial_x^2[\phi_{c_*}+v],\zeta \rangle_{L^2}.
\end{align*}
An analogous computation shows that
\begin{align*}
\overline{R}^\sigma_1(u,\alpha,\xi,\zeta)=&\ \sigma^2\| Q^{1/2}\hat{T}^*_{\alpha,\xi}\overline{\gamma}_s\|_{\mathcal{H}}^2\bigl(2\langle \phi_{c_*}+v,\zeta\rangle_{L^2}\\&+4 \langle x\partial_x[\phi_{c_*}+v],\zeta\rangle_{L^2}+ \langle x^2\partial_x^2[\phi_{c_*}+v],\zeta\rangle_{L^2}\bigr)
\end{align*}
and 
\begin{align*}
    \overline{R}^\sigma_2(u,\alpha,\xi,\zeta)=\sigma^2\| Q^{1/2}\hat{T}^*_{\alpha,\xi}\overline{\mu}_s\|_{\mathcal{H}}^2 \langle \partial_x[\phi_{c_*}+v],\zeta\rangle_{L^2}.
\end{align*}
In order to simplify the term $\overline{R}^\sigma_5(u,\alpha,\xi,\zeta)$, we rewrite the summands as
\begin{align*}
    \mu^{\sigma}_s [Q^{1/2}e_k] \partial^2_{u \xi} \langle \Phi_{\alpha,\xi}u,\zeta\rangle_{L^2}&\bigl[ M(u) [Q^{1/2}e_k]\bigr]\\
    &=-\sigma\alpha^3 \langle \hat{T}_{\alpha,\xi}Q^{1/2}e_k,\overline{\mu}_s \rangle_{\mathcal{H}}  \langle T_{\alpha,\xi}\partial_x M(u) [Q^{1/2}e_k],\zeta \rangle_{L^2}\\
    &=\sigma \alpha^3 \langle e_k,Q^{1/2}\hat{T}^*_{\alpha,\xi}\overline{\mu}_s \rangle_{\mathcal{H}}  \bigl\langle e_k, Q^{1/2}M^*(u)[\partial_x T^*_{\alpha,\xi}\zeta] \bigr\rangle_{\mathcal{H}},
\end{align*}
so that
\begin{align*}
    \overline{R}^\sigma_5(u,\alpha,\xi,\zeta)
    =&\ \sigma^2 \alpha^3 \bigl\langle Q^{1/2}\hat{T}^*_{\alpha,\xi} \overline{\mu}_s, Q^{1/2}M^*(u)[\partial_x T^*_{\alpha,\xi}\zeta]\bigr\rangle_{\mathcal{H}} \\
    =&-\sigma^2 \alpha^3 \langle T_{\alpha,\xi}\partial_x M(u)[Q\hat{T}^*_{\alpha,\xi} \overline{\mu}_s],  \zeta\rangle_{L^2} \\
    =&-\sigma^2 \alpha^3 \langle  M(T_{\alpha,\xi}u_x)[\hat{T}_{\alpha,\xi}Q\hat{T}^*_{\alpha,\xi} \overline{\mu}_s],  \zeta\rangle_{L^2}\\&-\sigma^2 \alpha^3 \langle  M(T_{\alpha,\xi}u)[\hat{T}_{\alpha,\xi}\hat{\partial}_x Q\hat{T}^*_{\alpha,\xi} \overline{\mu}_s],  \zeta\rangle_{L^2}. 
\end{align*}
After substituting \eqref{eqn:subst} for $j=0,1$ we find
\begin{align*}
    \overline{R}^\sigma_5(u,\alpha,\xi,\zeta)
    =&-\sigma^2  \bigl\langle  M\bigl(\partial_x[\phi_{c_*}+v]\bigr)[\hat{T}_{\alpha,\xi}Q\hat{T}^*_{\alpha,\xi} \overline{\mu}_s],  \zeta\bigr\rangle_{L^2}\\&-\sigma^2 \alpha \langle  M(\phi_{c_*}+v)[\hat{T}_{\alpha,\xi}\hat{\partial}_x Q\hat{T}^*_{\alpha,\xi} \overline{\mu}_s],  \zeta\rangle_{L^2}. 
\end{align*}
An analogous computation shows that
\begin{align*}
    \overline{R}^\sigma_4(u,\alpha,\xi,\zeta)
    =&-2\sigma^2 \langle M(\phi_{c_*}+v)[\hat{T}_{\alpha,\xi}Q\hat{T}_{\alpha,\xi}^*\overline{\gamma}_s], \zeta \rangle_{L^2} \\&-\sigma^2 \bigl\langle x  M\bigl(\partial_x[\phi_{c_*}+v]\bigr)[\hat{T}_{\alpha,\xi}Q\hat{T}^*_{\alpha,\xi}\overline{\gamma}_s] , \zeta\bigr\rangle_{L^2}\\
    &-\sigma^2\alpha \langle x  M(\phi_{c_*}+v)[\hat{T}_{\alpha,\xi}\hat{\partial}_x Q\hat{T}^*_{\alpha,\xi}\overline{\gamma}_s] , \zeta\rangle_{L^2}.
\end{align*}
In the term $\overline{R}^\sigma_0(u,\alpha,\xi,\zeta)$
we substitute \eqref{eqn:subst} for $j=0,1,3$ and 
\[-\partial_{x}^3[\phi_{c_*}+v]-2[\phi_{c_*}+v]\partial_{x} [\phi_{c_*}+v]=\mathcal{L}_{c_*}v+N(v)-c_*\partial_x [\phi_{c_*}+v]\]
to obtain
\[\overline{R}^\sigma_0(u,\alpha,\xi,\zeta)=  \alpha^{-3}  \langle \mathcal{L}_{c_*}v+N(v)-c_* \partial_x [\phi_{c_*}+v],\zeta \rangle_{L^2}.\]
In the martingale component $\overline{S}^\sigma$, we substitute the derivatives and \eqref{eqn:gammasrepresentation}-\eqref{eqn:musrepresentation}, which gives
\begin{align*}\overline{S}^\sigma(u,\alpha,\xi,\zeta)[h]=&\ \sigma \alpha^2 \langle M(T_{\alpha,\xi}u)[\hat{T}_{\alpha,\xi}h],\zeta\rangle_{L^2} \\
&-\sigma  (2\alpha^2 \langle T_{\alpha,\xi}u,\zeta\rangle_{L^2}+\alpha^3 \langle x T_{\alpha,\xi}u_x, \zeta\rangle_{L^2}) \langle \hat{T}_{\alpha,\xi}h,\overline{\gamma}_s\rangle_{\mathcal{H}}\\
&-\sigma \alpha^3 \langle T_{\alpha,\xi}u_x,\zeta \rangle_{L^2}\langle \hat{T}_{\alpha,\xi}h,\overline{\mu}_s\rangle_{\mathcal{H}}.
\end{align*}
Substituting \eqref{eqn:subst} then yields
\[\overline{S}^\sigma(u,\alpha,\xi,\zeta)[h]=\langle \sigma S(v)[\hat{T}_{\alpha,\xi}h],\zeta \rangle_{L^2},\]
where $S$ is defined in \eqref{eqn:martingalepart}.
We collect also that 
\[\overline{R}^\sigma_0(u,\alpha,\xi,\zeta)=\langle\alpha^{-3} \mathcal{L}_{c_*} v,\zeta \rangle_{L^2}+\langle R^\sigma_0(v,\alpha,\xi),\zeta \rangle_{L^2}\]
and
\[\overline{R}^\sigma_i(u,\alpha,\xi,\zeta)=\langle R^\sigma_i(v,\alpha,\xi),\zeta \rangle_{L^2},\]
for $i=1,\ldots,6$, where $R^\sigma_0,\ldots,R^\sigma_6$ are defined in \eqref{eqn:Rterms}. Since the test function $\zeta$ was arbitrary, we conclude that \eqref{eqn:perturbationgeneral} follows.

\section{Noise rescaling}
\label{app:noiserescaling}
In this appendix we collect several useful properties of the family $\{Q_\alpha\}_{\alpha>0}$ defined in \eqref{eqn:convolutionfamily}  in relation to the transformation operators $T_{\alpha,\xi}$ defined in \eqref{eqn:mappingT}.

\renewcommand{\thesection}{\Alph{section}}
\begin{Lemma}
\label{lem:scaledQ}
Let $\alpha,\beta>0$ and $\xi\in \R$. Denote by $T_{\alpha,\xi}^*$ the $L^2(\R)$-adjoint of $T_{\alpha,\xi}$. Then we have  the identities
\begin{enumerate}
    \item $(Q_{\alpha})^{1/2} =(Q^{1/2})_\alpha$;
    \item $T_{\alpha,\xi}Q_{\beta} = Q_{\beta\alpha} T_{\alpha,\xi}$;
    \item $T_{\alpha,\xi}^*=\alpha^{-1} T_{\alpha^{-1},-\alpha^{-1}\xi} .$
\end{enumerate}
\end{Lemma}
\renewcommand{\thesection}{\appendixname\Alph{section}}
\begin{proof}

Recall that $\hat{q}$ denotes the Fourier transform of the convolution kernel $q$, and $q_{1/2}:=\mathcal{F}^{-1}\{\sqrt{\hat{q}}\}$, where $\mathcal{F}^{-1}$ denotes the inverse Fourier transform. Let $f\in L^2(\R)$.
\begin{enumerate}
\item  We compute
\begin{align*}
   \alpha^2 q_{1/2}(\alpha \cdot) * q_{1/2}(\alpha \cdot) * f=&\ \mathcal{F}^{-1}\bigl\{\sqrt{\hat{q}}(\alpha^{-1}\cdot) \sqrt{\hat{q}}(\alpha^{-1} \cdot) \hat{f}\bigr\}\\
   =&\ \mathcal{F}^{-1}\bigl\{\hat{q}(\alpha^{-1}\cdot)  \hat{f}\bigr\}= \alpha q(\alpha \cdot)*f=Q_\alpha f.
\end{align*}
    \item Substituting $y=\alpha z+\xi$, we have
    \begin{align*}
        T_{\alpha,\xi}[Q_{\beta} f] =&\ \beta\int_{\R}q(\beta(\alpha x+\xi-y))f(y)\ \d y\\
    =&\ \beta\alpha \int_{\R}q(\beta\alpha (x-z))f(\alpha z+\xi)\ \d z\\
    =&\ \beta\alpha q(\beta\alpha \cdot)*T_{\alpha,\xi}[f]= Q_{\beta\alpha} T_{\alpha,\xi}[f].
\end{align*}
\item Substituting $y=\alpha x+\xi$, we have
\begin{align*}
    \langle T_{\alpha,\xi}[f],g\rangle_{L^2} =& \int_{\R}f(\alpha x+\xi)g(x)\ \d x\\
    =&\ \alpha^{-1} \int_{\R} f(y)g(\alpha^{-1}(y-\xi))\ \d y\\
    =&\ \langle f,\alpha^{-1} T_{\alpha^{-1},-\alpha^{-1}\xi}[g]\rangle_{L^2}.
\end{align*}
\end{enumerate}
\end{proof}

We now discuss the effect of the transformation $T_\alpha$ on the space-time white noise $W_t$. A defining property of the space-time white noise $W_t$ is the isometry
\begin{align}\label{eqn:whitenoiseisometry}
    \mathbb{E}\bigl[\int_0^{t_1}\langle \d W_s,w_1\rangle_{L^2}\int_0^{t_2}\langle \d W_s,w_2\rangle_{L^2}\bigr]=(t_1 \wedge t_2)\langle w_1,w_2\rangle_{L^2},
\end{align}
which holds for $t_1,t_2>0$ and $w_1,w_2\in L^2(\R)$. From the isometry \eqref{eqn:whitenoiseisometry} one obtains, upon differentiating with respect to $t_1$ and $t_2$ and choosing $w_1=\delta_x$ and $w_2=\delta_y$, the formal covariance identity
\begin{align*}
    \mathbb{E}\bigl[\frac{\d W_{t_1}(x)}{\d t}\frac{\d W_{t_2}(y)}{\d t}\bigr]=\delta(t_1-t_2)\delta(x-y).
\end{align*}
Let us consider the covariance structure that results from rescaling the space-time white noise $W_t$ via $T_{\alpha}$. Picking an orthonormal basis $\{e_k\}_{k=0}^\infty$ of $L^2(\R)$, we compute
\begin{align*}
    \langle T_{\alpha}e_k,T_{\alpha}e_j\rangle_{L^2}=&\int_{\R}e_k(\alpha x)e_j(\alpha x)\ \d x=\alpha^{-1}\int_{\R}e_k(y)e_j(y)\ \d y,
\end{align*}
which shows that $\{\alpha^{1/2}T_{\alpha}e_k\}_{k=0}^\infty$ is also an orthonormal basis of $L^2(\R)$. We then find via It\^o's isometry that,
\begin{align*}
    &\mathbb{E}\bigl[\int_0^{t_1}\langle  \alpha^{1/2}{T}_{\alpha} [\d W_s],w_1\rangle_{L^2}\int_0^{t_2}\langle  \alpha^{1/2}{T}_{\alpha} [\d W_s],w_2\rangle_{L^2}\bigr]\\
    &=\mathbb{E}\int_0^{t_1\wedge t_2} \sum_{k=0}^\infty  \langle \alpha^{1/2}{T}_{\alpha} e_k, w_1 \rangle_{L^2} \langle \alpha^{1/2}{T}_{\alpha} e_k, w_2\rangle_{L^2}  \ \d s,
\end{align*}
which via \eqref{eqn:whitenoiseisometry} gives the covariance identity 
\begin{align}
\label{eqn:scaledwhitenoiseisometry}
    \mathbb{E}\bigl[\int_0^{t_1}\langle  \alpha^{1/2}{T}_{\alpha} \d W_s,w_1\rangle_{L^2}\int_0^{t_2}\langle  \alpha^{1/2}{T}_{\alpha} \d W_s,w_2\rangle_{L^2}\bigr]=(t_1 \wedge t_2)\langle w_1,w_2\rangle_{L^2}.
\end{align}
Here we have used that $\{\alpha^{1/2}T_{\alpha}e_k\}_{k=0}^\infty$ is an orthonormal basis of $L^2(\R)$. We thus observe that the process $\tilde{W}_t:=\alpha^{1/2}T_{\alpha}W_t$ generates the same statistics as the white noise $W_t$.

\section{Expansions}
Table~\ref{tab:constants} collects evaluations the various constants that appear in the expansion of the modulation system in \S\ref{sec:solitondynamics}.
The evaluations  denoted with the symbol `$\approx$' have been computed numerically. For the exact evaluations, we have used that the defining equations consist of inner products between $\phi_{c_*}, \zeta_{c_*}$ and derivatives thereof. These can be written as integrals over hyperbolic functions, for which exact evaluations are available.

\begin{table}
\begin{center}
\begin{tabular}{||c | c | l||} 
 \hline
Constant/Function & Defining equation & Value \\ [0.5ex] 
 \hline\hline
 \hline
  $\overline{\gamma}_{d;I}(0)$ & \eqref{eqn:driftcomponentscalar} & $\tfrac{74}{135}+\tfrac{4\pi^2}{405}$ \\ 
  \hline
  $\overline{\mu}_{d;I}(0)$ & \eqref{eqn:driftcomponentscalar} & $(\tfrac{16\pi^2}{405}-\tfrac{34}{45})c_*^{-1/2}$ \\ 
  \hline
  $\overline{\gamma}_{d;\RomanIII}(0)$ & \eqref{eqn:gammadIII} & $\approx 0.093c_*^{1/2}$ \\ 
  \hline
  $\overline{\mu}_{d;\RomanIII}(0)$ & \eqref{eqn:gammadIII} & $\approx -0.1$ \\ 
 \hline
 $\overline{\gamma}_{s;I}(0)$& \eqref{eqn:scalargammas}&$\tfrac{2}{3}$\\ 
 \hline
 $\overline{\mu}_{s;I}(0)$& \eqref{eqn:scalargammas}&$-\tfrac{2}{3}c_*^{-1/2}$\\ 
 \hline
 $\overline{\gamma}_{\diamond}(0)$& \eqref{eqn:FGL2}&$\tfrac{1}{9}c_*^{-3/2}\phi_{c_*}^2$\\ 
 \hline
$\|\overline{\gamma}_{\diamond}(0)\|_{L^2}^2$ & \eqref{eqn:FGL2} & $\tfrac{4}{35}c_*^{1/2}$ \\
 \hline
$\overline{\mu}_{\diamond}(0)$ & \eqref{eqn:FGL2} & $\tfrac{2}{9}c_*^{-2}\phi_{c_*}^2-\tfrac{2}{9}c_*^{-1/2}\phi_{c_*}\zeta_{c_*}$ \\
 \hline
 $\|\overline{\mu}_{\diamond}(0)\|_{L^2}^2$ & \eqref{eqn:FGL2}&  $\approx 0.435c_*^{-1/2}$ \\
 \hline
\end{tabular}
\caption{\label{tab:constants}Values of various constants that appear in \S\ref{sec:solitondynamics}.}
\end{center}
\end{table}

We proceed by outlining how the expansions of the functionals $\overline{\gamma}_d^0, \overline{\mu}_d^0,\overline{\gamma}_d,\overline{\mu}_d$ and mappings $\overline{\gamma}_s,\overline{\mu}_s, R_0$ and $S$ that appear in \S\ref{sec:solitondynamics} can be computed.
Since these mappings are defined as products with the matrix $K^{-1}(v)$, we first derive an expansion
\[K^{-1}(v)=K^{-1}(0)+[K^{-1}]^{(1)}(v)+[K^{-1}]^{(2)}(v)+O(v^3).\]
We therefore write
\begin{align*}
     K(v)=\begin{bmatrix} 
9c_*^{3/2}&0\\
9 &-\tfrac{9}{2}c_*^{1/2}
\end{bmatrix}+\begin{bmatrix} 
b_1&b_2\\
b_3&b_4
\end{bmatrix},
\end{align*}
where $b_1,b_2,b_3$ and $b_4$ denote the functionals
\begin{align*}
    b_1=&\ \bigl\langle (x\partial_x+2)v,\phi_{c_*}\bigr\rangle_{L^2},\\
    b_2=&\ \bigl\langle  \partial_x v, \phi_{c_*}\bigr\rangle_{L^2},\\
    b_3=&\ \bigl\langle (x\partial_x+2)v,\zeta_{c_*}\bigr\rangle_{L^2},\\
   b_4=&\ \bigl\langle \partial_x v,\zeta_{c_*}\bigr\rangle_{L^2}.
\end{align*}
We furthermore write 
    \begin{align*}
        K^{-1}=\frac{1}{\det K}\left(
9\begin{bmatrix} 
-\tfrac{1}{2}c_*^{1/2}&0\\
-1 & c_*^{3/2}
\end{bmatrix}+\begin{bmatrix} 
b_4&-b_2\\
-b_3&b_1
\end{bmatrix} \right)
    \end{align*}
and compute
\begin{align*}
    \det K=&-\tfrac{81}{2}c_*^2+9c_*^{3/2}b_4-\tfrac{9}{2}c_*^{1/2}b_1-9 b_2+b_1 b_4-b_3 b_2.
\end{align*}
Using $\frac{1}{a+x}=\frac{1}{a}-\frac{x}{a^2}+\frac{x^2}{a^3}+O(x^3)$, we expand $\tfrac{1}{\det K(v)}$ as
\begin{align*}
    K^{-1}(0)=&\frac{1}{9}\begin{bmatrix}
c_*^{-3/2}& 0\\
2c_*^{-2}& -2c_*^{-1/2}\\
\end{bmatrix},\\
    [K^{-1}]^{(1)}
=&\frac{2}{81}\begin{bmatrix*}[l]
c_*^{-2}b_1+2c_*^{-5/2}b_2& & c_*^{-2}b_2\\
c_*^{-2}b_3+2c_*^{-5/2}b_4-c_*^{-7/2}b_1-2c_*^{-4}b_2 &  &-2c_*^{-1}b_4+2c_*^{-5/2}b_2
\end{bmatrix*}
\end{align*}
and
\begin{align*}
    [K^{-1}]^{(2)}=&\ \tfrac{2}{729}(-2c_*^{-5/2}b_4+c_*^{-7/2}b_1+2c_*^{-4}b_2)\begin{bmatrix}
        b_4 & -b_2\\
        -b_3 & b_1
    \end{bmatrix}\\
&+\tfrac{2}{729}\Bigl(-4c_*^{-3}b_4^2-c_*^{-5}b_1^2-4c_*^{-6}b_2^2+2c_*^{-4}b_1b_4\\
&+8c_*^{-9/2}b_2b_4-4c_*^{-11/2}b_1b_2+2c_0^{-4}b_3b_2\Bigr)\begin{bmatrix} 
-\tfrac{1}{2}c_*^{1/2}&0\\
-1 & c_*^{3/2}
\end{bmatrix}.    
\end{align*}


We then collect that
\begin{align*}
    \begin{bmatrix}
[\overline{\gamma}_d^0]^{(2)}(v)\\
[\overline{\mu}_d^0]^{(2)}(v)
\end{bmatrix}
&=K^{-1}(0)
\begin{bmatrix}
\langle N(v),\phi_{c_*}\rangle_{L^2} \\
\langle N(v),\zeta_{c_*}\rangle_{L^2}
\end{bmatrix}=
\begin{bmatrix}
\tfrac{1}{9}c_*^{-3/2}\langle \label{eqn:gammaexpansion}v,v\partial_x\phi_{c_*}\rangle_{L^2}\\
\tfrac{2}{9}\bigl\langle v,v( c_*^{-2}\partial_x \phi_{c_*}-c_*^{-1/2}\partial_c\phi_{c_*})\bigr\rangle
\end{bmatrix},
\end{align*}
and 
\[[R_0]^{(2)}(v)=-[\overline{\gamma}_d^0]^{(2)}(v) (2+x\partial_x)\phi_{c_*}-[\overline{\mu}_d^0]^{(2)}(v) \partial_x\phi_{c_*}.\]
For the functionals $\overline{\gamma}_s$ and $\overline{\mu}_s$ we obtain
\begin{align*}
\begin{bmatrix}
[\overline{\gamma}_s]^{(1)}(v)\\
[\overline{\mu}_s]^{(1)}(v)
\end{bmatrix}=
 K^{-1}(0)
\begin{bmatrix}
M^*(v)[\phi_{c_*}]\\
 M^*(v)[\zeta_{c_*}]
\end{bmatrix}
+[K^{-1}]^{(1)}(v)
\begin{bmatrix}
M^*(\phi_{c_*})[\phi_{c_*}]\\
 M^*(\phi_{c_*})[\zeta_{c_*}]
\end{bmatrix}
\end{align*}
and
\begin{align*}
\begin{bmatrix}
[\overline{\gamma}_s]^{(2)}(v)\\
[\overline{\mu}_s]^{(2)}(v)
\end{bmatrix}=
 [K^{-1}]^{(1)}(v)
\begin{bmatrix}
M^*(v)[\phi_{c_*}]\\
 M^*(v)[\zeta_{c_*}]
\end{bmatrix}
+[K^{-1}]^{(2)}(v)
\begin{bmatrix}
M^*(\phi_{c_*})[\phi_{c_*}]\\
 M^*(\phi_{c_*})[\zeta_{c_*}]
\end{bmatrix}.
\end{align*}
Using these expressions, we have
\begin{align*}
   [S]^{(1)}(v)[h]=&M(v)[h] - (x \partial_x+2)\phi_{c_*}\bigl\langle h,[\overline{\gamma_s}]^{(1)}(v)\bigr\rangle_{\mathcal{H}}\nonumber\\
   &- (x \partial_x+2)v\bigl\langle h,\overline{\gamma_s}(0)\bigr\rangle_{\mathcal{H}}-  \partial_x\phi_{c_*} \bigl\langle h,[\overline{\mu}_s]^{(1)}(v)\rangle_{\mathcal{H}}-  \partial_x v \langle h,\overline{\mu}_s(0)\bigr\rangle_{\mathcal{H}}.
\end{align*}
Lastly, we have
\begin{align*}
\begin{bmatrix}
[\overline{\gamma}_d]^{(1)}(v,\alpha)\\
[\overline{\mu}_d]^{(1)}(v,\alpha)
\end{bmatrix}=&
-\alpha K^{-1}(0)\sum_{i=1}^5 \begin{bmatrix}
\bigl\langle [R_i]^{(1)} (v,\alpha),\phi_{c_*}\bigr\rangle_{L^2}\\
\bigl\langle [R_i]^{(1)} (v,\alpha),\zeta_{c_*}\bigr\rangle_{L^2}
\end{bmatrix}\\
&-\alpha [K^{-1}]^{(1)}(v)\sum_{i=1}^5 \begin{bmatrix}
\langle R_i (0,\alpha),\phi_{c_*}\rangle_{L^2}\\
\langle R_i (0,\alpha),\zeta_{c_*}\rangle_{L^2}
\end{bmatrix}
\end{align*}
and
\begin{align*}
\begin{bmatrix}
[\overline{\gamma}_d]^{(2)}(v,\alpha)\\
[\overline{\mu}_d]^{(2)}(v,\alpha)
\end{bmatrix}=&
-\alpha K^{-1}(0)\sum_{i=1}^5 \begin{bmatrix}
\bigl\langle [R_i]^{(2)} (v,\alpha),\phi_{c_*}\bigr\rangle_{L^2}\\
\bigl\langle [R_i]^{(2)} (v,\alpha),\zeta_{c_*}\bigr\rangle_{L^2}
\end{bmatrix}\\
&-\alpha [K^{-1}]^{(1)}(v)\sum_{i=1}^5 \begin{bmatrix}
\bigl\langle [R_i]^{(1)} (v,\alpha),\phi_{c_*}\bigr\rangle_{L^2}\\
\bigl\langle [R_i]^{(1)} (v,\alpha),\zeta_{c_*}\bigr\rangle_{L^2}
\end{bmatrix}\\
&-\alpha [K^{-1}]^{(2)}(v)\sum_{i=1}^5 \begin{bmatrix}
\langle R_i (0,\alpha),\phi_{c_*}\rangle_{L^2}\\
\langle R_i (0,\alpha),\zeta_{c_*}\rangle_{L^2}
\end{bmatrix}.
\end{align*}

\section{Numerical schemes}
\label{app:schemes}
Here, we describe the numerical schemes that were employed to simulate the stochastic KdV equation \eqref{eqn:skdvgeneral} and the modulation system \eqref{eqn:modulationv}-\eqref{eqn:modulationxi}, in the cases of scalar noise and space-time white noise. We employ a semi-implicit finite-difference scheme from \cite{lord}, as was also used in \cite{cartwright}.

In space, we use $N+1$ grid points $x_n=n\Delta x-L$ for $n=1,\ldots,N+1$, where $\Delta x=\frac{2L}{N}$ and $L>0$ is the right-boundary of the computational domain $[-L,L]$. We denote the numerical solution to \eqref{eqn:skdvgeneral} at time $j\Delta t$ by $U^j=[U_1^j,U^j_2,\ldots,U_{n}^j,\dots,U_{N+1}^j]^T$, and denote by $D_1,D_2$ and $D_3$ the $(N+1)\times (N+1)$ centered finite difference matrices of second order for the differential operators $\partial_x,\partial_x^2$ and $\partial_x^3$, respectively. We initialize by using an Euler-Maruyama step as
\[U^1=U^0-\Delta t(D_3*U^0+2U^0(D_1*U^0))+\sigma U^0 \Delta W, \]
where the noise is discretized as $\Delta W\sim \sqrt{\Delta t}N(0,1)$ in the case of scalar noise, and as $\Delta W=[W_1,W_2,\ldots,W_n,\ldots,W_{N+1}]^T$ with $W_n\stackrel{i.i.d.}{\sim} \sqrt{\frac{\Delta t}{\Delta x}}N(0,1)$ for $n=1,\ldots,N+1$ in the case of space-time white noise. Thereafter, for $j\geq 1$, the scheme continues with a semi-implicit step and a two-step Adam-Bashforth discretization for the nonlinear term as
\begin{align*}
    U^{j+1}=&\ \bigg(I+\frac{\Delta t}{2}D_3\bigg)^{-1}\bigg[\bigg(I-\frac{\Delta t}{2}D_3\bigg)*U^j+\sigma U^j \Delta W\\
    & \qquad \qquad -3\Delta t U^j(D_1*U^j)+\Delta t U^{j-1}(D_1*U^{j-1})\bigg].
\end{align*}

To simulate \eqref{eqn:modulationv}, we also employ a semi-implicit method. Denote by $V^j=[V_1^j,V^j_2,\ldots,V_{n}^j,\dots,V_{N+1}^j]^T$ and  $A^j=[A_1^j,A^j_2,\ldots,A_{n}^j,\dots,A_{N+1}^j]^T$ the numerical solutions to \eqref{eqn:modulationv} and \eqref{eqn:modulationalpha}, respectively. For $j\geq 0$, the numerical solution $V^{j+1}$ to \eqref{eqn:modulationv} at time $j\Delta t$ is computed as
\begin{align*}
      V^{j+1}=&\ \bigg(I-\frac{\Delta t}{2}L_0\bigg)^{-1}\bigg[\bigg(I+\frac{\Delta t}{2}L_0\bigg)*V^j-2\Delta t(A^j)^{-3} V^j(D_1*V^j) \\
      & \qquad \qquad +\Delta t \sum_{i=0}^2 R^\sigma_{i}( V^j,A^j)+ \sigma S(V^j) \Delta W\bigg],
\end{align*}
where \[L_0=-D_3+c_*D_1-2D_1*\operatorname{Diag}(\Phi_0)\]
and $\operatorname{Diag}(\Phi_0)$ is a diagonal matrix with entries 
\[\phi_{c_*}(-L),\phi_{c_*}(-L+\Delta x),\ldots,\phi_{c_*}(-L+n\Delta x),\ldots,\phi_{c_*}(L)\] 
on the diagonal.

The numerical solutions $A^{j+1}$  and  $X^{j+1}=[X_1^j,X^j_2,\ldots,X_{n}^j,\dots,X_{N+1}^j]^T$ to \eqref{eqn:modulationalpha} and \eqref{eqn:modulationxi} at time $j\Delta t$ , are respectively given for $j\geq0$  by 
\begin{align*}
     A^{j+1} =&\ A^j+[-(A^j)^{-2}\overline{\gamma}_d^0(V^j)+\sigma^2 A^j \overline{\gamma}_{d;I}(V^j)]\Delta t -\sigma A^j \overline{\gamma}_{s;I}(V^j) \Delta W, \\
    X^{j+1} =&\ X^j+[-(A^j)^{-2}\overline{\mu}_d^0(V^j)+\sigma^2 A^j \overline{\mu}_{d;I}(V^j)] \Delta t -\sigma A^j \overline{\mu}_{s;I}(V^j) \Delta W , 
\end{align*}
 in the case of scalar noise, and as
\begin{align*}
     A^{j+1} =&\ A^j+[-(A^j)^{-2}\overline{\gamma}_d^0(V^j)+\sigma^2  \overline{\gamma}_{d;\RomanIII}(V^j)]\Delta t -\sigma A^j \langle  \Delta W,\overline{\gamma}_{\diamond}(V^j)\rangle, \\
    X^{j+1} =&\ X^j+[-(A^j)^{-2}\overline{\mu}_d^0(V^j)+\sigma^2 \overline{\mu}_{d;\RomanIII}(V^j)] \Delta t -\sigma A^j \langle  \Delta W,\overline{\mu}_{\diamond}(V^j)\rangle , 
\end{align*}
 in the case of space-time white noise. We remark that, to obtain path-wise correspondence in this latter case between numerical solutions to \eqref{eqn:skdvgeneral} and the modulation system \eqref{eqn:modulationv}-\eqref{eqn:modulationxi}, realizations of the noise $\Delta W$ should be shifted and rescaled via the map $T_{\alpha,\xi}$.

\pagebreak
\section{Supplementary figures}
\label{app:supplementary}

\begin{figure*}[h!]
    \centering
    \begin{subfigure}[t]{0.5\textwidth}
        \centering
        \includegraphics[height=2.3in]{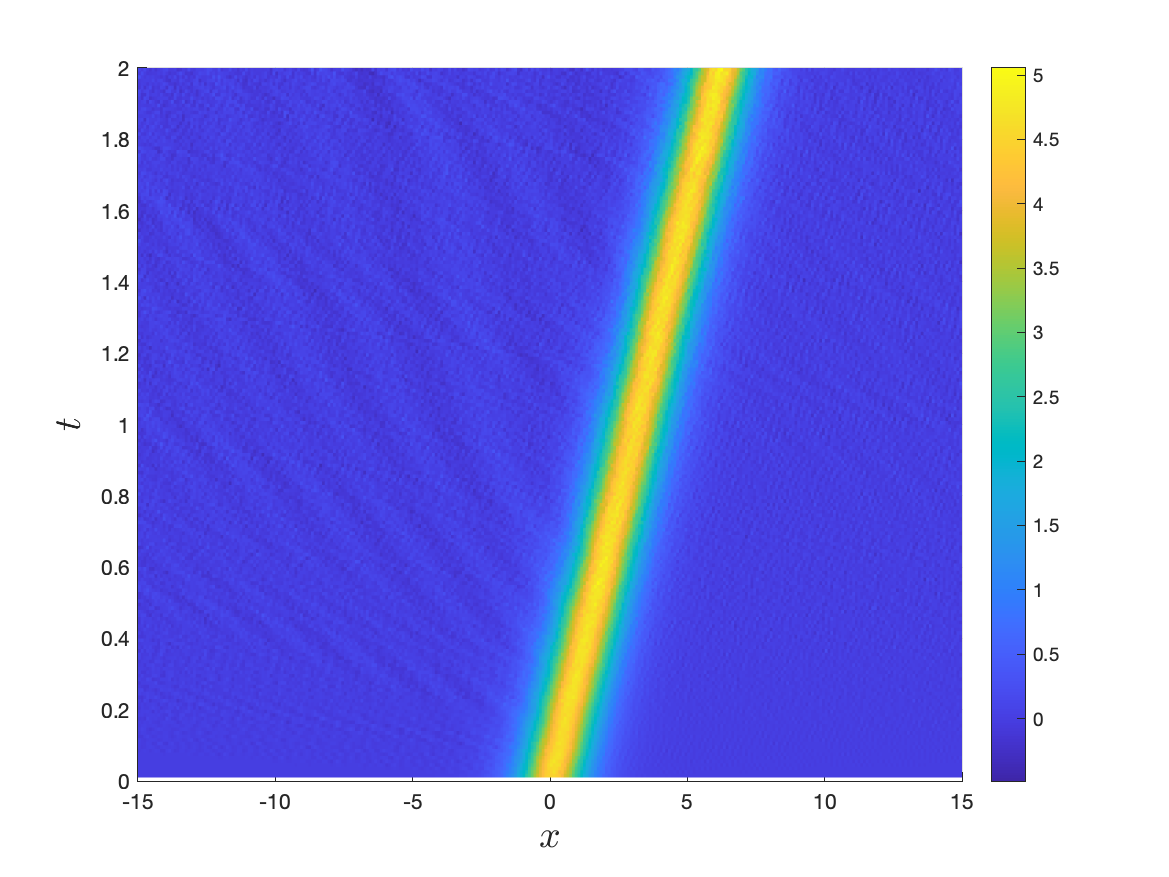}
        \caption{Original frame.}
    \end{subfigure}%
    ~ 
    \begin{subfigure}[t]{0.5\textwidth}
        \centering
        \includegraphics[height=2.3in]{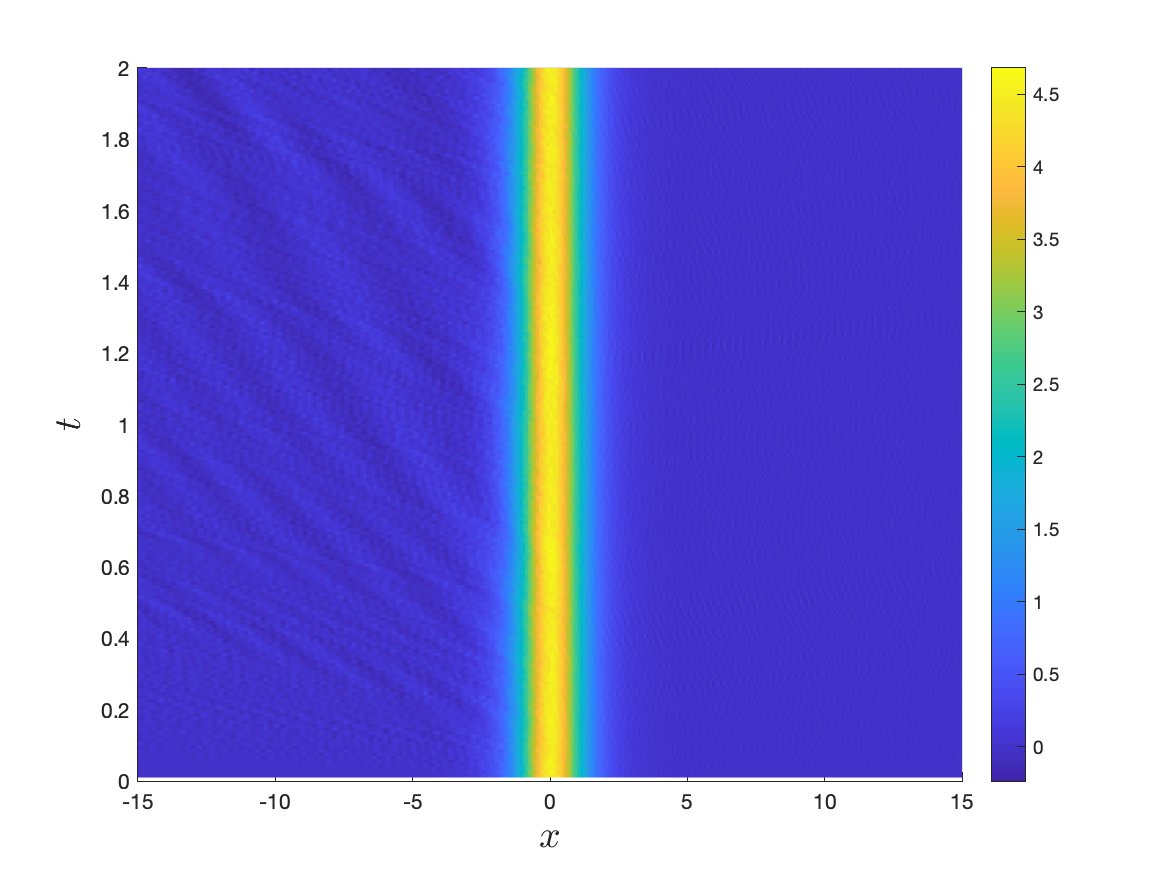}
        \caption{Stochastic co-moving frame.}
    \end{subfigure}
    \begin{subfigure}[t]{0.5\textwidth}
        \centering
        \includegraphics[height=2.3in]{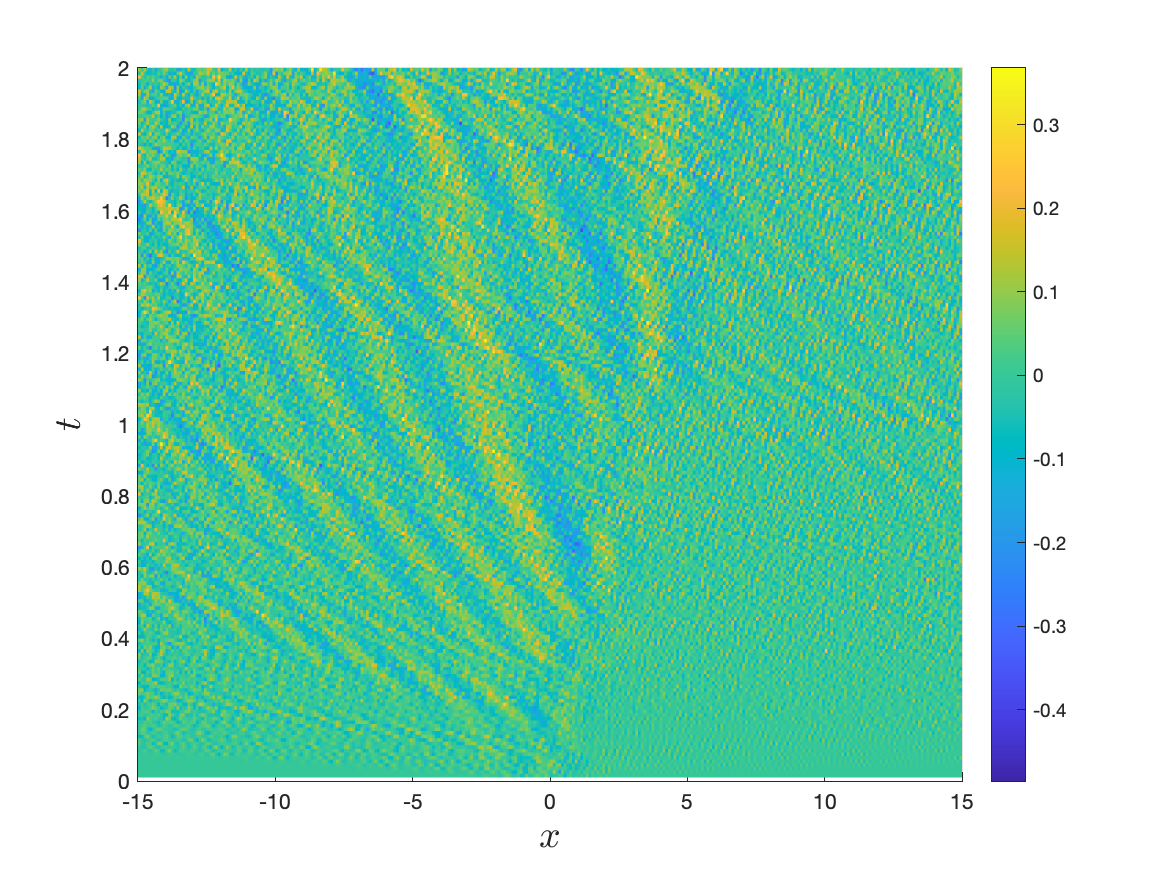}
        \caption{Original frame, soliton removed. }\label{subfig:originalpwhite}
    \end{subfigure}%
    ~ 
    \begin{subfigure}[t]{0.5\textwidth}
        \centering
        \includegraphics[height=2.3in]{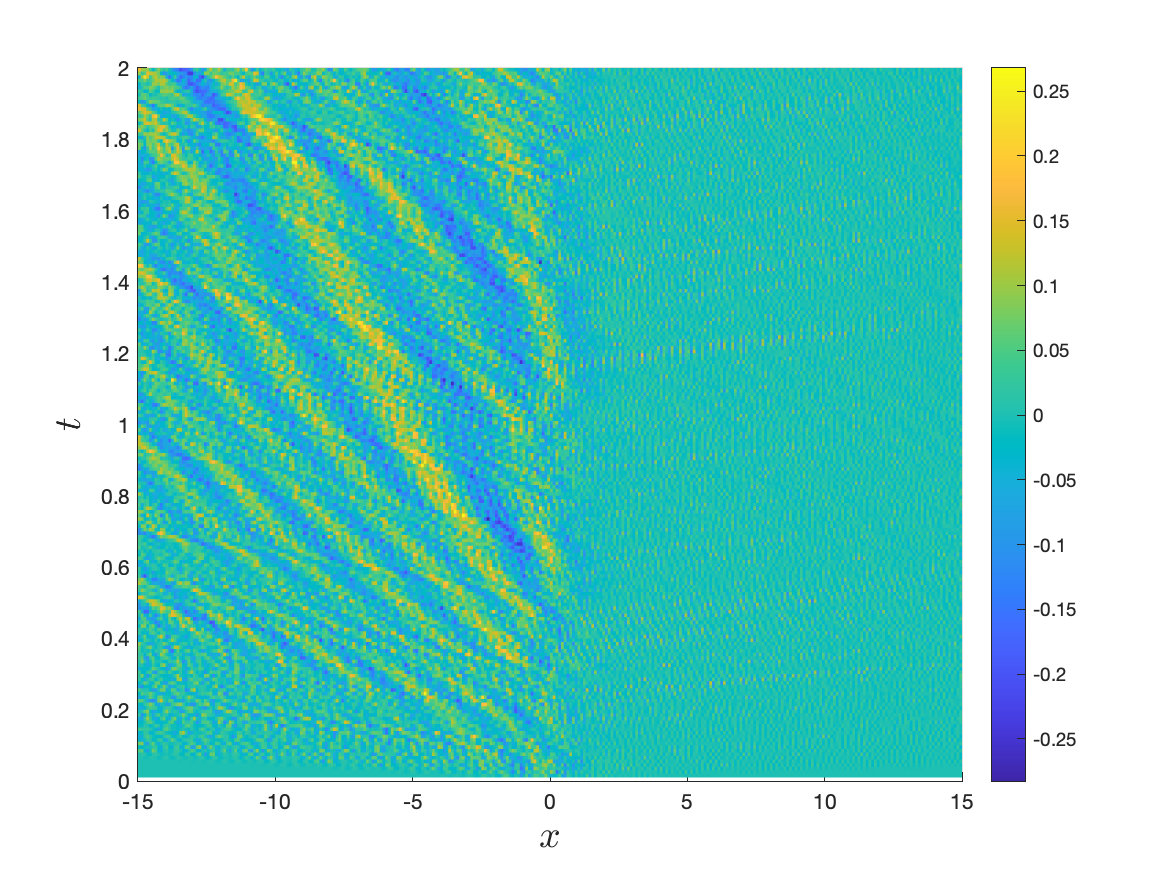}
        \caption{Stochastic co-moving frame, soliton removed.}\label{subfig:frozenpwhite}
    \end{subfigure}
    \caption{Simulation of the KdV equation with space-time white noise of strength $\sigma=0.05$. The original frame realization shows $u(t,x)$, from a simulation of \eqref{eqn:skdvgeneral}. The frozen frame simulation shows $ \phi_{c_*}(x)+v(t,x)$, from simulation of \eqref{eqn:modulationv}-\eqref{eqn:modulationxi} with the same realization of the noise. Figure~\ref{subfig:originalpwhite} and Figure~\ref{subfig:frozenpwhite} show the perturbation with respect to the soliton. The original frame realization shows $u(t,x)-\phi_{c(t)}(x-\xi(t))$ with the phase-definitions \eqref{eqn:phasedefinitions}, and the frozen frame simulation shows $v(t,x)$.}\label{fig:pathwiseperturbation}
\end{figure*}

\begin{figure*}[h!]
\centering
           \begin{subfigure}[t]{0.5\textwidth}
    \centering
        \includegraphics[height=2in]{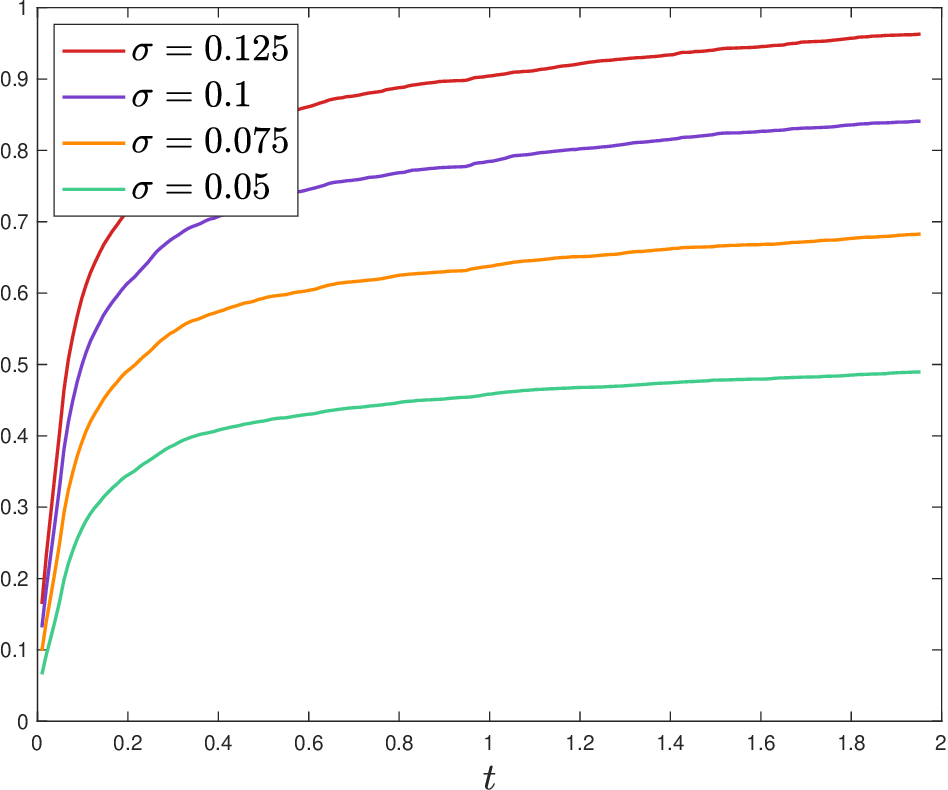}
        \caption{Perturbation size over time.}
    \end{subfigure}%
    ~
    \begin{subfigure}[t]{0.5\textwidth}
    \centering
        \includegraphics[height=2in]{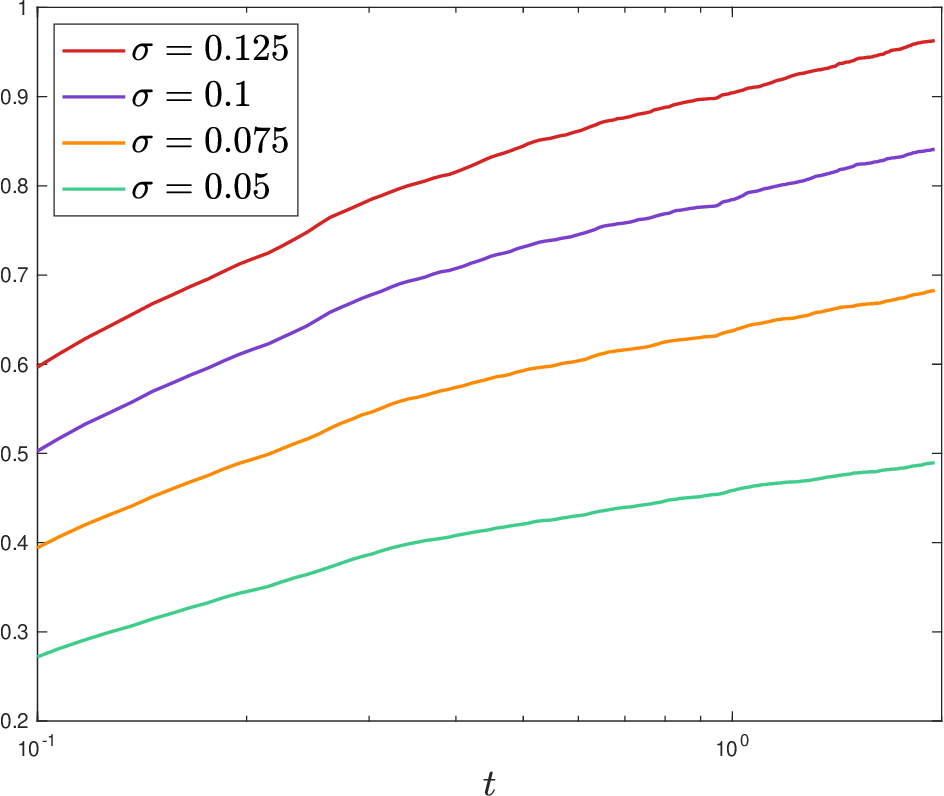}
        \caption{Perturbation size over time, log-scale.}\label{subfig:perturbationlogwhite}
        \end{subfigure}
        \caption{Sample mean of the process $\sup_{s\leq t}\|v(s)\|_{L_a^2([-40,10])}$ for space-time white noise, computed over $200$ realisations for $\sigma\in\{0.05,0.075,0.1,0.125\}.$ This simulation was computed on the computational domain $[-40,40]$, with values $c_*=3$ and $a=0.15$.}\label{fig:perturbationwhite}
\end{figure*}

\begin{figure*}[h!]
    \centering
        \includegraphics[height=2.5in]{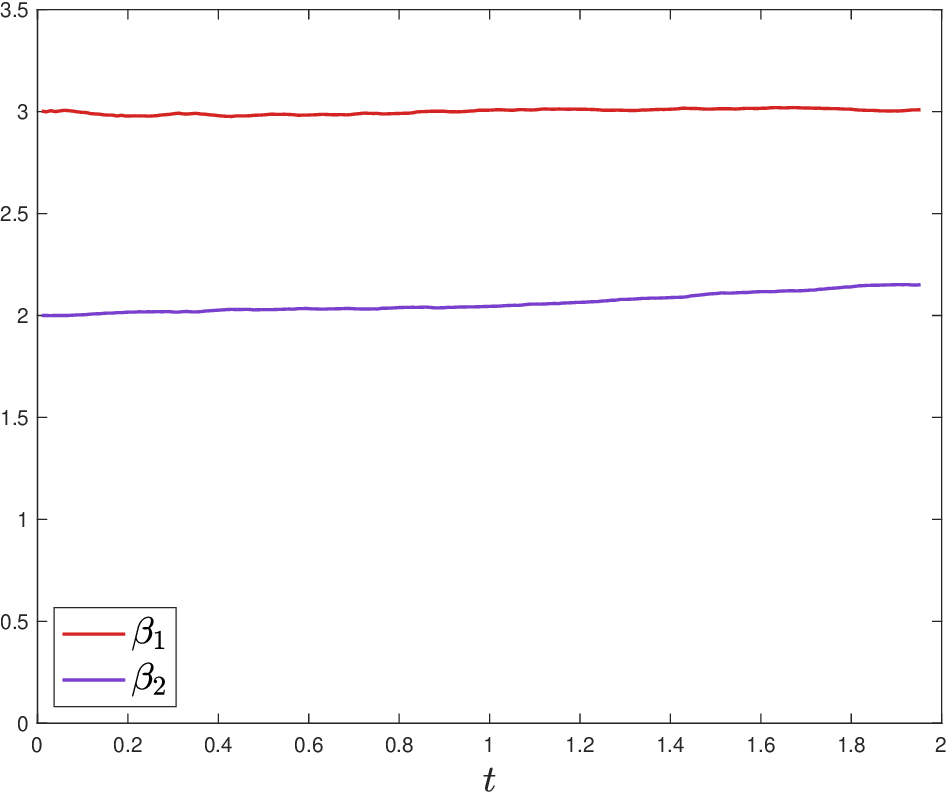}
    \caption{Estimation of the orders $\beta_1$ and $\beta_2$ at which the remainders $|c-c_2|$ and $\|v(t)- v_1(t)\|_{L_a^2}$ depend on the noise strength $\sigma$. Here, $\beta_1(t)$ is obtained from a least squares fit of $\mathbb{E}\sup_{s\leq t}|c(s)-c_2(s)|$ (as in Figure~\ref{fig:errors}) to $k_1(t)\sigma^{\beta_1(t)}$. Similarly, the exponent $\beta_2(t)$ is obtained from a least squares fit of $\mathbb{E}\sup_{s\leq t}\|v(s)- v_1(s)\|_{L_a^2}$ to $k_2(t)\sigma^{\beta_2(t)}$.}\label{fig:errororder}
\end{figure*}
 \bibliographystyle{elsarticle-num} 
 \bibliography{bibliography}





\end{document}